\documentclass[10pt,reqno]{amsart}
\usepackage[utf8]{inputenc} 
\usepackage[T1]{fontenc} 
\usepackage[french,english]{babel} 
\usepackage{lmodern} 
\usepackage{enumerate} 
\usepackage{geometry} 
\usepackage{fancyhdr} 
\usepackage{verbatim} 
\usepackage{amsmath,amssymb, amsthm} 
\usepackage{dsfont} 
\usepackage{graphicx} 
\usepackage{tikz} 
\usepackage{float} 
\usetikzlibrary{matrix,arrows}
\usepackage{url}
\usepackage{xypic} 
\usepackage[all]{xy}
\usepackage{stmaryrd}

\usepackage{hyperref} 
\hypersetup{pdfstartview=XYZ} 

\usepackage[final]{pdfpages} 

\title{Irreducible dual of p-adic $ U(5) $}

\author{Claudia Schoemann}
\date{October 2014}

\usepackage[leqno]{amsmath}     
\usepackage{amssymb}            
\usepackage{mathrsfs}           
\usepackage{dsfont}             
\usepackage{amsxtra}            

\newtheoremstyle{th}{8pt}{8pt}{\itshape}{}{\bfseries}{ ---}{5pt}{\thmname{#1}\thmnumber{ #2} \thmnote{\bfseries (#3).}}
\theoremstyle{th}

\theoremstyle{break}
\newtheorem{Theorem}{Theorem}[section]

\newtheorem{Lemma}[Theorem]{Lemma}

\newtheorem{Corollary}[Theorem]{Corollary}

\newtheorem{Proposition}[Theorem]{Proposition}
\newenvironment{Proof}{\begin{proof}} {\end{proof}}

\newtheoremstyle{def}{8pt}{8pt}{}{}{\bfseries}{}{5pt}{\thmname{#1}\thmnumber{ #2} \thmnote{\bfseries (#3).}}
\theoremstyle{def}
\newtheorem{Definition}[Theorem]{Definition}

\newtheorem{Remark}[Theorem]{Remark}

\theoremstyle{nonumberplain}

\newcommand{\field}[1]{\mathds{#1}}
\newcommand{\Z}{\field{Z}}              
\newcommand{\N}{\field{N}}              
\newcommand{\R}{\field{R}}
\newcommand{\Q}{\field{Q}}              
\newcommand{\C}{\field{C}}              

\newcommand{\trans}[1]{{^t\!{#1}}}

\DeclareMathOperator{\GL}{GL}

\DeclareMathOperator{\Cent}{Cent}

\DeclareMathOperator{\Hom}{Hom}

\DeclareMathOperator{\G}{G}

\DeclareMathOperator{\id}{id}

\DeclareMathOperator{\Ind}{Ind}

\DeclareMathOperator{\F}{F}

\DeclareMathOperator{\Gal}{Gal}

\DeclareMathOperator{\Lg}{Lg}

\DeclareMathOperator{\cusp}{cusp}

\DeclareMathOperator{\St}{St}

\DeclareMathOperator{\supp}{supp}

\DeclareMathOperator{\E}{E}

\DeclareMathOperator{\Normu}{Norm}

\usepackage[scaled=.90]{helvet}

\begin{document}

\begin{abstract}

We study the parabolically induced complex representations of the unitary group in 5 variables, $ U(5), $ defined over a
p-adic field. 

Let $ F $ be a p-adic field. Let $ E : F $ be a field extension of degree two. Let $ \Gal(E : F ) = \{ \id, \sigma \}. $
We write $ \sigma(x) = \overline{x} \; \forall x \in E. $ Let $ E^* := E 
\setminus \{ 0 \} $ and let $ E^1 := \{x \in E \mid x \overline{x} = 1 \}. $

$ U(5) $ has three proper standard Levi subgroups, the minimal Levi subgroup $ M_0 \cong E^* \times E^* \times
E^1 $ and the two maximal Levi subgroups $ M_1 \cong \GL(2, E) \times E^1 $ and $ M_2 \cong E^* \times U(3). $

We consider representations induced from $ M_0 $ and from non-cuspidal, not fully-induced representations of $ M_1 $ and
$ M_2. $ We determine the points and lines of reducibility and the irreducible
subquotients of these representations.

\end{abstract}

\maketitle

\tableofcontents
\section{Introduction}

We study the parabolically induced complex representations of the unitary group in 5 variables - $ U(5) $ - defined over a 
non-archimedean local field of characteristic $0$. This is $ \Q_p $ or a finite extension of $ \Q_p, $ where
$ p $ is a prime number. We speak of a 'p-adic field'.

\medskip

Let $ F $ be a p-adic field. Let $ E:F $ be a field extension of degree two. Let $ \Gal(E:F) = \{ \id, \sigma \} $ be the
Galois group. We write $ \sigma(x) = \overline{x} \; \forall x \in E. $ Let $ E^* := E \setminus \{ 0 \} $
and let $ E^1 := \{ x \in E \mid  x \overline{x} = 1 \}. $ Let $ \mid \; \mid_E $ denote the p-adic norm on 
$ E. $

\medskip

$ U(5) $ has three proper parabolic subgroups. Let $ P_0 $ denote the minimal standard parabolic
subgroup and $ P_1 $ and $ P_2 $ the two maximal standard parabolic subgroups. One has the Levi decomposition
$ P_i = M_i N_i, \; i \in \{ 0, 1, 2 \}, $ where $ M_i $ denote the standard Levi subgroups and $ N_i $ are the corresponding
unipotent subgroups of $ U(5). $

\medskip

$ M_0 \cong E^* \times E^* \times E^1 $ is the minimal Levi subgroup, $ M_1 \cong \GL(2,E) \times E^1 $ and
$ M_2 \cong E^* \times U(3) $ are the two maximal Levi subgroups.

\medskip

We consider representations of the Levi subgroups, extend them trivially to the unipotent subgroups to obtain 
representations of the parabolic groups. We now perform normalised parabolic induction to obtain representations of
$ U(5). $

\medskip

We consider representations of $ M_0, $ further we consider non-cuspidal, not fully-induced representations of $ M_1 $ and $ M_2. $
For $ M_1 $ this means that the representation of the $ \GL(2,E) - $ part is a proper subquotient of a representation
induced from $ E^* \times E^* $ to $ \GL(2,E) $. For $ M_2 $ this means that the representation of the $ U(3) - $ part of 
$ M_2 $ is a proper subquotient of a representation induced from $ E^* \times E^1 $ to $ U(3). $

\smallskip

As an example for $ M_1, $ take $ \mid\det\mid_p^{\alpha} \chi(\det) \St_{\GL_2} \rtimes \lambda', $ where $ \alpha \in \R, \;
\chi $ is a unitary character of $ E^*, \; \St_{\GL_2} $
is the Steinberg 
representation of $ \GL(2,E) $ and $ \lambda' $ is a character of $ E^1. $ As an example for $ M_2, $ take 
$ \mid \; \mid^{\alpha} \chi \rtimes \lambda'(\det) \St_{U(3)}, $ where $ \alpha \in \R,
\; \chi $ is a unitary character of $ E^*, \; \lambda' $ is a character of $ E^1 $ and $ \St_{U(3)} $ is the Steinberg
representation of $ U(3). $ Note that $ \lambda' $ is unitary.

\bigskip

We determine the points and lines of reducibility of the representations of $ U(5), $ and
we determine the irreducible subquotients.

\bigskip

The irreducible complex representations of $ U(3) $ over a p-adic field obtained as subquotients of parabolically induced
representations have been classified by Charles 
David Keys in \cite{Ky}, the irreducible complex representations of $ U(4) $ over a p-adic field obtained as subquotients 
of parabolically induced representations have
been classified by Kazuko Konno in \cite{Ko}.

\bigskip

We start with some basic definitions. In section 3 we give the classification of the irreducible non-cuspidal
representations of $ U(3), $ as has been done by Charles David Keys \cite{Ky}.
We reassemble the results for the irreducible unitary representations.

\medskip

In section 4 we determine when the induced representations to $ U(5) $ are irreducible. It is done for representations
induced from $ M_0 $ and for non-cuspidal, not fully-induced representations of the two maximal Levi subgroups  $ M_1 $
and $ M_2. $

\smallskip

Representations of $ M_0 $ are of the form $ \mid \; \mid_p^{\alpha_1}
\chi_1 \otimes \mid \; \mid_p^{\alpha_2} \chi_2 \otimes \lambda', $ where $ \alpha_1, \alpha_2 \in \R, \; \chi_1, \chi_2 $
are unitary characters of $ E^* $ and $ \lambda' $ is a unitary character of $ E^1. $ Reducibility of the induced
representation $ \mid \; \mid_p^{\alpha_1}
\chi_1 \times \mid \; \mid_p^{\alpha_2} \chi_2 \rtimes \lambda' $ depends on $ \alpha_1, \alpha_2 $ and on the two unitary
characters $ \chi_1 $ and $ \chi_2. $

\medskip

Let $ N_{E/F}(E) $ denote the norm map of $ E $ with respect to the field extension $ E : F, $ i.e.
$ N_{E/F}(x) = x \overline{x} \; \forall x \in E. $

\medskip

In Theorems 4.2, 4.3 and 4.5 we show that for $ \alpha_1, \alpha_2 \in \R_+ \; \mid \; \mid_p^{\alpha_1}
\chi_1 \times \mid \; \mid_p^{\alpha_2} \chi_2 \rtimes \lambda' $ is reducible if and only if at least 
one of the following
cases holds:

\begin{enumerate}
 \item $ \mid \alpha_1 - \alpha_2 \mid = 1 $ and $ \chi_1 = \chi_2, $

\item $ \mid \alpha_1 + \alpha_2 \mid = 1 $ and $ \chi_1(x) = \chi_2^{-1}(\overline{x}) \; \forall x \in E^*, $

\item $ \exists i \in \{ 1, 2 \} $ s.t. $ \alpha_i = 1 $ and $ \chi_i = 1, $

\item $ \exists i \in \{1, 2 \} $ s.t. $ \alpha_i = 1/2 $ and $ \chi_i \mid F^* \neq 1, $ but $ \chi_i \mid N_{E/F}(E^*) =
1, $

\item $ \exists i \in \{1, 2 \} $ s.t. $ \alpha_i = 0 $ and $ \chi_i \neq 1, $ but $ \chi_i \mid F^* = 1. $

\end{enumerate}

Let $ \chi $ be a unitary character of $ E^*. $ The condition that $ \chi(x) = \chi^{-1}(\overline{x}) \; \forall x \in E^* $
is equivalent to the condition that $ \chi \mid N_{E/F}(E^*) = 1 $
and to the fact that $ \chi $ is a character of a type as in 3. 4. or 5. of the list above.

\bigskip

In 4.3 we consider representations induced from irreducible non-cuspidal representations of $ M_1 $ and $ M_2 $
that are not fully-induced.

\smallskip

We consider $ \mid \det \mid_p^{\alpha} \chi \St_{\GL_2} \rtimes \lambda' $ and 
$ \mid \det \mid_p^{\alpha} \chi 1_{\GL_2} \rtimes \lambda', $ where $ \alpha \in \R_+, \; \chi $ is a unitary character of
$ E^*, \; \St_{\GL_2} $ is the Steinberg representation of $ \GL(2,E), \; \lambda' $ is a unitary character of $ E^1 $ and
$ 1_{\GL_2} $ is the trivial representation of $ \GL(2,E). $

\smallskip

In Theorem 4.7 and Proposition 4.8 we show that for $ \alpha \in \R_+, \; \mid \; \mid_p^{\alpha} \chi \St_{\GL_2} \rtimes 
\lambda' $ and $ \mid \; \mid_p^{\alpha} \chi 1_{\GL_2} \rtimes \lambda' $ are irreducible unless one of the following
cases holds:

\begin{enumerate}
 \item $ \alpha = 1/2 $ or $ \alpha = 3/2 $ and $ \chi = 1, $
 \item $ \alpha = 0, \; \alpha = 1/2 $ or $ \alpha = 1 $ and $ \chi \mid F^* \neq 1, $ but $ \chi \mid N_{E/F}(E^*) = 1, $
\item $ \alpha = 1/2 $ and $ \chi \neq 1 $ but $ \chi \mid F^* = 1. $
\end{enumerate}

We consider $ \mid \; \mid^{\alpha} \chi \rtimes \tau, $ where $ \alpha \in \R_+, \; \chi $ is a unitary character of $ E^* $
and $ \tau $ is an irreducible non-cuspidal unitary representation of $ U(3) $ that is not fully-induced. We consider all
irreducible proper subquotients $ \tau $ of representations induced to $ U(3) $ from its unique proper Levi-subgroup 
$ M \cong E^* \times E^1, $ as classified by Charles David Keys in \cite{Ky}.

\smallskip

In Theorems 4.9, 4.11, 4.13 and in Propositions 4.10, 4.12 and Remark 4.14 we show that these representations
are irreducible unless one has a certain combination of $ \alpha \in \{0, 1/2, 1, 3/2, 2 \} $ and $ \chi = 1,$ or 
$ \chi \neq 1 $ but $ \chi \mid \F^* = 1, $
or $ \chi \mid F^* \neq 1 $ but $ \chi \mid N_{E/F}(E^*) = 1. $

\bigskip
  
In the case of reduciblility the irreducible subquotients are determined in the course of section 4. In
several cases the irreducible suquotients are determined in section 5.

\medskip

In section 5 we treat several 'special' reducibility points of representations induced from the minimal parabolic
subgroup $ P_0 \cong M_0 N_0, $ where $ M_0 \cong E^* \times E^* \times E^1. $ 
In some cases the induced representation $ \mid \; \mid_p^{\alpha_1}
\chi_1 \times \mid \; \mid_p^{\alpha_2} \chi_2 \rtimes \lambda' $ has more than two irreducible subquotients, more precisely 
it has four irreducible subquotients.

\medskip

If this is the case, then $ \alpha_1, \alpha_2 \in \{0, 1/2, 1, 3/2, 2 \} $ and $ \chi_i(x) = \chi_i^{-1}(\overline{x})
$ for $ i =1,2 $ and $ \forall x \in E^*. $ I.e. $ \chi_i = 1, $ or $ \chi_i \neq 1 $ but $ \chi_i \mid F^* = 1, $ or $ \chi_i \mid F^* \neq 1 $ but 
$ \chi_i \mid N_{E/F}(E^*) = 1, $ for $ i = 1,2. $

\medskip

We determine the irreducible subquotients in terms of Langlands-quotients, and we determine whether these 
Langlands-quotients are unitary.

\section{Definitions}

Let $ \G $ be a connected reductive algebraic group, defined over a non-archimedean local field of characteristic $ 0. $
Let $ V $ be a vector space, defined over the complex numbers. Let $ \pi $ be a representation of
 $ \G $ on $ V. $ We denote it by $ (\pi,V) $ and sometimes by $ \pi $ or $ V. $
\smallskip

Let $ \mathcal{K} $ denote the set of open compact subgroups of $ \G. $

\begin{Definition}
 A representation $ (\pi, V) $ of $ \G $ is said to be \textbf{smooth} if every $ v \in V $ is fixed by a neighbourhood of the
unity in $ \G. $ This is equivalent to saying that $ \exists K \in \mathcal{K} $ such that $ \pi(k)v = v \; \forall k \in 
K. $ 
\end{Definition}

\begin{Definition}
 A representation $ ( \pi, V) $ of $ \G $ is said to be \textbf{admissible} if it is smooth and for every $ K \in \mathcal{K}
$ the space
$ V^K $ of fixed vectors under $ K $ is finite dimensional.
\end{Definition}

From now on let $ (\pi, V) $ be an admissible representation of $ \G. $

\begin{Definition} Let $ (\pi, V) $ be a smooth representation of $ \G. $ Let $ V^* = \Hom(V, \C) $ be the space of
linear forms on $ V. $ One defines a representation $ (\pi^*, V^*)$ of $ \G: $ if $ v^* \in V^* $ and $ g \in G $ then 
$ \pi^*(g)(v^*) $ is defined by $ v \mapsto v^*(\pi(g^{-1})(v)). $ Let $ \overset{\sim}{V} $ be the subspace of smooth
vectors in $ V^* ,$ i.e. the subset 
$ \overset{\sim}{V} \subset 
V^* $ of elements $ \overset{\sim}{v} $ 
such that the stabiliser of $ \overset{\sim}{v} $ is an open subgroup of $ \G. $ One shows that $ \overset{\sim}{V} $
 is invariant under $ \G $ for the action of $ \pi^*. $ So $ \pi^* $ induces a representation
$ \overset{\sim}{\pi} $ on $ \overset{\sim}{V}; \; (\overset{\sim}{\pi}, \overset{\sim}{V}) $ is
called the \textbf{dual representation} of $ (\pi, V). $  

\end{Definition}

\begin{Definition} A representation $ \pi  $ is called \textbf{hermitian} if 
 $  \pi \cong  \overline{\overset{\sim}{\pi}}. $
\end{Definition}

\begin{Definition} A representation $ (\pi,V) $ of a group $ \G $ is called a \textbf{unitary representation} if and only if on the
vector space
 $ V $ there exists a positive definite hermitean form $ \langle \; , \; \rangle: V \times V \rightarrow \C $ that is
 invariant under the action of $ \G: $
\end{Definition}

$$ \langle \pi(g) v, \pi(g) w \rangle = \langle v,w \rangle \; \forall g \in \G, \forall v,w \in V.  $$

\begin{Definition} 
A \textbf{matrix coefficient} $ f_{v,\overset{\sim}{v}} $ of a representation $ (\pi,V) $ is a (locally constant) function: $ f_{v,
\overset{\sim}{v}}: \G \rightarrow \C: 
g \mapsto \overset{\sim}{v}(\pi(g)(v)), $ where $ v \in V $ and $ \overset{\sim}{v} \in \overset{\sim}{V}. $
\end{Definition}

\begin{Definition} An irreducible representation $ \pi $ of $ \G $ is called \textbf{cuspidal} if $ \pi $ has a non-zero
matrix
coefficient 
$ f: \G \rightarrow \C $ that is
compactly supported modulo the center of G.
\end{Definition} 

\begin{Definition} An irreducible representation $ \pi $ of $ \G $  is called \textbf{square-integrable} if $ \pi $ is unitary and if
$ \pi $ has a non-zero matrix coefficient
 $ f: \G \rightarrow \C $ that is square-integrable
modulo the center $ Z $ of $ \G: \; f \in L^2(\G / Z). $ It follows that all matrix 
coefficients of $ \pi $ are square-integrable.
\end{Definition}

\begin{Definition}  An irreducible representation $ \pi $ of a group $ \G $ is \textbf{tempered} if it is unitary and if
$ \pi $ has a
non-zero matrix coefficient
$ f: \G \rightarrow \C $ that verifies $ f \in L^{2 + \epsilon}(\G / Z) \; \forall \epsilon > 0. $
\end{Definition}

\begin{Remark} Square-integrable representations are tempered.
\end{Remark}

\bigskip

Let $ F $ be a non-archimedean local field of characteristic $ 0. $ i.e. $ \Q_p $ or a finite extension of $ \Q_p, $ where
$ p $ is a prime number.

Let $ E:F $ be a field extension of degree two, hence a Galois extension. Let $ \Gal(E:F) = \{ \id, \sigma \}. $
The action of the non-trivial element $ \sigma $ of the Galois group is called the conjugation of elements in $ E $ corresponding
to the extension $ E: F. $
We write $ \sigma(x) = \overline{x} \; \forall \;x \in E $ and extend $ \overset{-}{\;} $ to
matrices with entries in $ E. $

\medskip

Let $ \Phi \in \GL(n,E) $ be a hermitian matrix, i.e. $ \overline{\Phi}^t = \Phi, $ let $ U_{\Phi} $ be the unitary group
definded by $ \Phi: $

$$ U_{\Phi} = \{ g \in \GL(n,E): g \Phi \overline{g}^t = \Phi \}. $$

\bigskip

Let $ \Phi_n = (\Phi_{ij}), $ where $ \Phi_{ij} = (-1)^{i-1} \delta_{i,n+1-j} $ and $ \delta_{ab} $ is the Cronecker delta.

Let $ \zeta \in E^* $ be an element of trace $ 0, $ i.e. $ tr(\zeta) = \zeta + \overline{\zeta} = 0. $

\smallskip
If  n is odd, then
$ \Phi_n = \left(
\begin{smallmatrix}
&&&&&&1\\
&&&&&\cdotp&\\
&&&&\cdotp&&\\
&&&\cdotp&&&\\
&&1&&&&\\
&-1&&&&&\\
1&&&&&&
\end{smallmatrix}
\right) $ is hermitian. If n is even, 
$ \zeta \Phi_n = \left(
\begin{smallmatrix}
 &&&&&&\zeta\\
&&&&&\cdotp&\\
&&&&\cdotp&&\\
&&&\cdotp&&&\\
&& - \zeta &&&&\\
&\zeta&&&&&\\
- \zeta&&&&&&
\end{smallmatrix}
\right) $ is hermitian.

Denote by \textbf{$ U(n) $} the \textbf{unitary group} corresponding to $ \Phi_n $ if $ n $ is odd or to $ \zeta \Phi_n $ if $ n $ is even,
 respectively.
It is quasi split.

\bigskip

Let $ n $ be an even positive integer. We will call Levi subgroup of $ U(n) $ a subgroup of block diagonal matrices

$ L: = \{ \left(
\begin{smallmatrix}
 A_1&&&&&&-&0\\
&A_2&&&&&&\mid\\
&&\ddots&&&&&\\
&&&A_k&&&&\\
&&&&B&&&\\
&&&&&\trans{\overline{A}}_k^{-1}&&\\
\mid &&&&&&\ddots&\\
0&-&&&&&&\trans{\overline{A}}_1^{-1}
\end{smallmatrix}
\right), \; A_i \in \GL_{n_i}(E) \; \text{ for } \;  1 \leq i \leq k \; \text {and } B \in U(m) \}, $

where $ m, n_1, \ldots, n_k $ are strictly positive integers such that $ m + 2 \underset{i = 1}{\overset{k}{\sum}} n_i = n. $
(If $ k = 0, $ then there are no $ n_i $ and $ L = U(n). $)

It is canonically isomorphic to the product
$ \GL(n_1,E) \times \ldots \times \GL(n_k, E) \times U(m).$
We chose the corresponding parabolic subgroup $ P $ such that it contains $ L $ and the subgroup of upper triangular matrices
in $ U(n). $ We call a parabolic subgroup $ P $ that contains the subgroup of upper triangular matrices standard.
Let $ N $ be the unipotent subgroup with identity matrices for the block diagonal matrices of $ L, $ arbitray entries in 
$ E $ above and $ 0's $ below. Then one has the Levi decomposition $ P = L N. $ 

\smallskip

Let $ \pi_i, \; i = 1, \ldots k, $ be smooth admissible
representations of $ \GL(n_i,E) $ and $ \sigma $ a smooth admissible representation of $ U(m). $ Let 
$ \pi_1 \otimes \ldots \otimes \pi_k \otimes \sigma $ denote the representation of 
$ L = \GL(n_1,E) \times \ldots \times \GL(n_k,E)
\times U(m) $ and denote by
$ \pi: = \Ind_P^{U(n)} (\pi_1 \otimes \ldots \pi_k \otimes \sigma) = \pi_1 \times \ldots \pi_k \rtimes \sigma $ the
normalized parabolically induced representation, where $ P $ is the corresponding standard parabolic subgroup containing
$ L $.    

\begin{Definition} Let $ \pi $ be an irreducible representation of $ \GL(n,E). $ Then there exist irreducible
cuspidal representations $ \rho_1, \rho_2, \ldots, \rho_k $ of general linear groups that are, up to isomorphism,
uniquely defined by $ \pi, $ such that $ \pi  $ is
isomorphic to a subquotient of $ \rho_1 \times \ldots \times \rho_k. $ The multiset of equivalence classess $ (\rho_1,
\ldots, \rho_k) $ is called the \textbf{cuspidal support} of $ \pi. $ It is denoted by $ \supp(\pi).$
\end{Definition}

\begin{Definition}
Let $ n \in \N $ and let $ \tau $ be an irreducible representation of $ U(n). $ Then there exist irreducible cuspidal
representations $ \rho_1, \ldots, \rho_k $ of general linear groups and an irreducible cuspidal representation $ \sigma $
of some $ U(m) $ that are, up to isomorphism and replacement of $ \rho_i $ by $ \rho_i^{-1}(\overset{-}{\;}) $ for some
$ i \in \{ 1, \ldots, k \},  $ uniquely defined by $ \tau, $ s. t. $ \tau $ is isomorphic to a subquotient of
$ \rho_1 \times \ldots \times \rho_k \rtimes \sigma. $ The representation $ \sigma $ is called the \textbf{partial cuspidal
support} of $ \tau $ and is denoted by $ \tau_{\cusp}. $
\end{Definition}

\begin{Definition}
Let $ \pi $ be a smooth representation of finite length of $ G. $ Then  $\hat{\pi} $ denotes the \textbf{Aubert dual} of 
$ \pi, $ as
defined in \cite{MR1285969}.
\end{Definition}

\bigskip

\section{The representations of $ U(3) $}
 
\medskip

\subsection{The irreducible representations of $ U(3) $}

\label{irru3}

$ U(3) $ has a unique standard proper parabolic subgroup, denoted by $ P. $ Let $ M $ be the standard Levi subgroup and 
$ N $ the
unipotent radical corresponding to $ P. $

\medskip

Then
$ M = \{
\left( \begin{smallmatrix}
 x&0&0\\
0&k&0\\
0&0&\overline{x}^{-1}
\end{smallmatrix} \right), x \in E^*, k \in E^1 \}, $

where $ E^* $ is the group of invertible elements of $ E $ and $ E^1 := \{ x \in E :  x \overline{x} = 1 \}. $

\bigskip

$ N = \{
\left( \begin{smallmatrix}
 1&\alpha&\beta\\
0&1&\overline{\alpha}\\
0&0&1
\end{smallmatrix} \right), \alpha, \beta \in E, \alpha \overline{\alpha} = \beta + \overline{\beta} \}. $
Set

$$ P = M N = \{
\left( \begin{smallmatrix}
 x&0&0\\
0&k&0\\
0&0&\overline{x}^{-1}
\end{smallmatrix} \right)
\left( \begin{smallmatrix}
 1&\alpha&\beta\\
0&1&\overline{\alpha}\\
0&0&1
\end{smallmatrix} \right) =
\left( \begin{smallmatrix}
 x&x\alpha&x\beta\\
0&k&k\overline{\alpha}\\
0&0&\overline{x}^{-1}
\end{smallmatrix} \right), x \in E^*, k \in E^1, \alpha, \beta \in E, \alpha \overline{\alpha} = \beta + \overline{\beta} \}. $$

\bigskip

For a smooth character $ \lambda \in \Hom(M,\C^*) $ there exist unique smooth characters $ \lambda_1 \in \Hom(E^*, \C^*) $ and
$ \lambda' \in \Hom(E^1, \C^*) $ such that

$$ \lambda \big(
\left( \begin{smallmatrix}
 x&0&0\\
0&k&0\\
0&0&\overline{x}^{-1}
\end{smallmatrix} \right) \big) =
\lambda_1(x) \lambda'(x\overline{x}^{-1}k), \; \forall x \in E^*, \forall k \in E^1.
$$

\bigskip

Every smooth character of $ E^* $ can be written in the form

$ \lambda_1(x) = \mid x \mid^{\alpha_1} \chi_1(x), $ with $ \alpha_1 \in \R $ and $ \chi_1 $ a unitary character. $ \lambda':
E^1 \rightarrow \C^* $ is smooth and $ E^1 $ is a compact group, hence $ \lambda' $ is unitary.

\bigskip

These are all characters of $ M. $ We extend $ \lambda $ from $ M $ to $ P, $ by taking $ \lambda \mid N = 1 $.
 
\medskip

We induce parabolically from $ P $ to $ U(3) $ and obtain

$$ \pi: = \Ind_P^{U(3)}(\lambda) = \Ind_P^{U(3)}(\lambda_1 \otimes \lambda') =: \lambda_1 \rtimes \lambda'. $$

The complex vector space V of the representation $ \pi $ is defined as follows:

$$ V: = \{ f: U(3) \rightarrow  \C : f \; \text{smooth and} f(mng) = \delta_P^{1/2} (m) \lambda (m) f(g) \;\forall m \in M,
\forall n \in N, \forall
g \in U(3) \}. $$

Here $ \delta_P^{1/2} $ is the modulus character. $ \pi $
acts on $ V $ by right translations.

\bigskip

Let $ \alpha \in \R_+^* $ and $ \chi $ be a unitary character of $ E^*.$ Let $ \lambda' $ be a character of $ E^1. $ 
Let $ \lambda = \mid \; \mid^{\alpha} \chi \otimes \lambda' $ be a character of the Levi subgroup $ M $ and 
$ \mid \; \mid^{\alpha} \chi \rtimes \lambda' $ the parabolically induced representation to $ U(3). $ 
Then $ \mid \; \mid^{\alpha} \chi \rtimes \lambda' $ 
has a unique irreducible quotient denoted by $ \Lg(\mid \; \mid^{\alpha} \chi \rtimes \lambda'), $ the Langlands
quotient.

\bigskip

Let $ N_{E/F}(\;) $ denote the norm on $ E $ corresponding to the field extension $ E/F $ of degree 2: $ \; N_{E/F}(x) = x
\overline{x} \; \forall x \in E. \; N_{E/F}(E^*) \subset F^* $ and $ \mid F^* / N_{E/F}(E^*) \mid = 2. $

Let $ \omega_{E/F}: F^* \rightarrow \C^* $ be the unique non-trivial smooth character with $ \omega_{E/F} \mid N_{E/F}(E^*) = 1. \;
\omega_{E/F} $ is determined by local class field theory.

Let $ X_{\omega_{E/F}} $ be the set of characters $ \chi $ of $ E^* $ such that 
$ \chi \mid F^* = \omega_{E/F}. $ Characters in $ X_{\omega_{E/F}} $ are unitary.

\medskip

Let $ X_{1_{F^*}} $ be the set of characters $ \chi $ of $ E^* $ that are non-trivial and whose restriction to
 $ F^* $ is trivial: $ \chi \neq 1, \; \chi \mid F^* = 1. $

\medskip

Let $ X_{N_{E/F}(E^*)} = \{ 1 \} \cup X_{\omega_{E/F}} \cup X_{1_{F^*}}. \; X_{N_{E/F}(E^*)} $ exhausts all characters 
$ \chi $ of $ E^* $ that are trivial on $ N_{E/F}(E^*), $ i.e. that
verify $ \chi(x) = \chi^{-1}(\overline{x}) \; \forall x \in E^*. $

\bigskip

Let
$ T := \{
\left(
\begin{smallmatrix}
 x&&\\
&1&\\
&&\overline{x}^{-1}
\end{smallmatrix}
\right), x \in F^* \} $
be the maximal split torus over $ F. $

Let $ \Normu_{U(3)}(T) $ be the normaliser of $ T $ in $ U(3), \; $ let $ \Cent_{U(3)}(T) $ be the centraliser of $ T $ in
$ U(3). $
The Weyl group is $ W := \Normu_{U(3)}(T)/\Cent_{U(3)}(T) \cong \Z/2 \Z \cong
\{ \id, \; 
\left(
\begin{smallmatrix}
 &&1\\
&1&\\
1&&
\end{smallmatrix}
\right) \}. $

\bigskip
By David Keys \cite{Ky} the induced representation $ \Ind_P^{U(3)}(\lambda) $ is irreducible except in the following cases:

\bigskip

1. \; $\lambda_1 = \mid \; \mid_E^{\pm 1} $

\smallskip

2. \; $ \lambda_1 = \mid \; \mid_E^{\pm 1/2} \chi_{\omega_{E/F}}, $ where $ \chi_{\omega_{E/F}} \in X_{\omega_{E/F}}. $

\smallskip

3. \; $ \lambda_1 = \chi_{1_{F^*}}, $ where $ \chi_{1_{F^*}} \in X_{1_{F^*}}. $

\bigskip

Note that in 1. and 2. changing the sign of the exponent is equivalent to replacing $ \lambda $ by $ w \lambda, $ where
$ w \in W $ is the non-trivial element of the Weyl group. Thus, the sign of the exponent does not affect the set of
irreducible constituents. We give the classification for positive exponent, for negative exponent the irreducible 
constituents exchange roles.

The classification does not depend on $ \lambda'. $
\medskip

In the first case $ \Ind_P^{U(3)}(\lambda) $ has exactly two constituents, the character
$ \psi := \lambda' \circ \det = \Lg(\lambda_1; \lambda') $ and the square-integrable subrepresentation
$ \St_{U(3)} \cdot \psi. $ \\
$ \psi $ and $ \St_{U(3)} \cdot \psi $ are unitary.

\medskip

In the second case $ \Ind_P^{U(3)}(\lambda) $ has exactly two constituents, a square-integrable (and hence unitary) 
representation
 $ \pi_{1,\chi_{\omega_{E/F}}} $ and a non-tempered unitary representation $ \pi_{2, \chi_{\omega_{E/F}}} = 
\Lg(\lambda_1; \lambda'). $

\medskip

In the third case $ \Ind_P^{U(3)}(\lambda) $ decomposes into the direct sum
$ \sigma_{1, \chi_{1_{F^*}}} \oplus \sigma_{2, \chi_{1_{F^*}}}. $ The two constituents $ \sigma_{1, \chi_{1_{F^*}}} $ and
 $ \sigma_{2, \chi_{1_{F^*}}} $ are tempered, hence unitary.

\bigskip

\begin{Remark} $ \pi_{2, \chi_{\omega_{E/F}}} $ is unitary:

Let $ \chi_{\omega_{E/F}} \in X_{\omega_{E/F}}. \; \chi_{\omega_{E/F}} \rtimes \lambda' $ is irreducible and unitary,
 $ \mid \; \mid^{\alpha} \chi_{\omega_{E/F}} \rtimes \lambda' $ is 
irreducible and unitary for $ \alpha \in (0, 1/2), $ by Theorem \ref{unitaryu3}, (1.3). By \cite{MR0324429}
the irreducible subquotients $ \pi_{\chi_{1,\chi_{\omega_{E/F}}}} $ and 
$ \pi_{2, \chi_{\omega_{E/F}}} $ of $ \mid \; \mid^{1/2} \chi_{\omega_{E/F}} \rtimes \lambda' $ are
unitary. See Theorem \ref{unitaryu3}.
\end{Remark}

\medskip

\begin{Remark} $ \sigma_{1, \chi_{1_{F^*}}} $ and $ \sigma_{2, \chi_{1_{F^*}}} $ are tempered:

In the third case  $ \lambda_1 =: \chi_{1_{F^*}} \in X_{1_{F^*}}. \; \chi_{1_{F^*}} $ is square-integrable. Hence 
$ \chi_{1_{F^*}} \rtimes \lambda' $ is tempered and so are its constituents.
\end{Remark}

\medskip

We obtain the following

\smallskip

\begin{Corollary} If $ \lambda_1 \rtimes \lambda' $ is reducible there are always two distinct irreducible subquotients. They are unitary.
\end{Corollary}

\bigskip
\bigskip

\subsection{The irreducible unitary representations of $ U(3) $}

\medskip

\begin{Proposition} Let $ \alpha > 0, $ and let $ \chi $ be a smooth unitary character of $ E^*. $
The following list exhausts all irreducible hermitian representations of $ U(3) $ supported in its parabolic subgroup 
$ P. $ 
\end{Proposition}

\medskip

0. $ \chi \rtimes \lambda', \; \chi \notin X_{1_{F^*}}, \; \sigma_{1, \chi_{1_{F^*}}}, \sigma_{2, \chi_{1_{F^*}}} $ as introduced above, tempered,

1. $ \lambda'(\det) = \Lg(\mid \; \mid 1; \lambda'), \pi_{2, \chi_{\omega_{E/F}}} = \Lg(\mid \; \mid^{1/2} \chi ;
 \lambda'), \; \chi \in X_{\omega_{E/F}} $ non-tempered, unitary,

2. $ \lambda'(\det) \St_{U(3)}, \pi_{1, \chi_{\omega_{E/F}}} $ square-integrable,

3. $  \mid \; \mid^{\alpha} 1 \rtimes \lambda', \alpha \neq 1; \; \mid \; \mid^{\alpha} \chi \rtimes
\lambda', \; \alpha \neq 1/2, \; \chi \in X_{\omega_{E/F}}; \mid \; \mid^{\alpha} \chi \rtimes \lambda',
 \; \chi \in X_{1_{F^*}}. $

\bigskip

\begin{Proof} Representations of  0., 1. and 2. are unitary, hence hermitian.

If for $ \alpha > 0 \; \mid \; \mid^{\alpha} \chi \rtimes
\lambda' $ is reducible, all subquotients are hermitian and part of the list.
\medskip

3. Let $ \mid \; \mid^{\alpha} \chi \rtimes
\lambda', \; \alpha > 0, $ be irreducible. By \cite{Ca}, 3.1.2, $ \mid \; \mid^{\alpha} \chi \rtimes
\lambda' \cong \overset{\sim}{\overline{\mid \; \mid^{\alpha} \chi \rtimes
\lambda'}} $ iff $ w (\mid \; \mid^{\alpha} \chi \otimes
\lambda') \cong \overset{\sim}{\overline{\mid \; \mid^{\alpha} \chi \otimes
\lambda'}} $ for the non-trivial element $ w $ of $ W. $

We have $ \overset{\sim}{\overline{\mid \; \mid^{\alpha} \chi \otimes
\lambda'}} = \mid \; \mid^{-\alpha} \chi \otimes \lambda' = w (\mid \; \mid^{\alpha} \chi \otimes \lambda') =
 \mid \; \mid^{-\alpha} \chi^{-1}(\overset{-}{\;}) \otimes \lambda' \Leftrightarrow \chi \cong \chi^{-1}(\overset{-}{\;}), $
i.e. $ \chi \in X_{N_{E/F}}. $
\end{Proof}

Let $ \alpha \in \R $ and $ \chi $ be a smooth unitary character of $ E^*. $
Like before we set $ \lambda = \lambda_1 \otimes \lambda', $ where $ \lambda_1 = \mid \; \mid^{\alpha} \chi. $ 

If $ \Ind_P^{U(3)}(\lambda) $ reduces we have seen that all subquotients are unitary.

\begin{Theorem}

\label{unitaryu3}

\medskip

1. \;$ \Ind_P^{U(3)}(\lambda) $ is irreducible and unitary if and only if

\smallskip

1.1  \; $ \chi \notin X_{1_{F^*}} $ and $ \alpha = 0, $ 

1.2 \; $ \chi = 1 $ and $ \alpha \in ]-1,0[ \cup ]0, 1[, $ 

1.3  \; $ \chi \in X_{\omega_{E/F}} $ and $ \alpha \in ] - 1/2, 0 [ \cup ] 0, 1/2 [. $

\bigskip

2. \; $ \Ind_P^{U(3)}(\lambda) $ is irreducible and non-unitary if and only if

\smallskip

2.1 \; $\chi_1 \neq 1, \; \chi_1 \notin X_{\omega_{E/F}} \; \forall \alpha \in \R^*. $ 

2.2 \; $ \chi_1 = 1 $ and $ \alpha \in ] - \infty, -1[ \cup ] 1, \infty [, $

2.3 \; $ \chi_1 \in X_{\chi_{\omega_{E/F}}} $ and $ \alpha \in 
] - \infty, - 1/2 [ \cup ] 1/2, \infty [. $

\bigskip

\end{Theorem}

\begin{Proof}

We use the fact that every character of $ E^* $ can be written as
$ \lambda_1 = \mid \; \mid^{\alpha} \chi $, where $ \alpha \in \R $ and $ \chi $ is a unitary character.

\bigskip

1.1. $ \Leftarrow $ If $ \alpha = 0, $ then $ \lambda = \chi \otimes \lambda' $ is unitary, hence
$ \Ind_P^{U(3)}(\lambda) $ is unitary. The representation $ \Ind_P^{U(3)}(\lambda) $ is irreducible unless
$ \chi \in X_{1_{F^*}} $ \cite{Ky}.

\bigskip

For the rest of the proof we need some preparatory theory.

For $ \alpha \in \R_+^*, $ let $ \pi_{\alpha} = \mid \; \mid^{\alpha} \chi \rtimes \lambda' $
 be a representation of $ U(3) $ and $ V $ be the corresponding vector space.
We give, on the same vector space $ V, $ a family of $ U(3) - $ invariant hermitian forms, parametrised by $ \alpha \in 
\R_+^*: $ 

$ \langle \; , \; \rangle_{\alpha}: V \times V \rightarrow \C, (f,h) \mapsto 
\underset{U(3, \mathcal{O})}{\int} A(w, \lambda) f(k) \overline{h(k)} dk.$

$ w $ is the non-trivial and the longest element of $ W, $ and 
$ A(w, \lambda): \; \mid \; \mid^{\alpha} \chi \rtimes \lambda' \rightarrow \mid \; \mid^{- \alpha} 
\chi^{-1}(\overset{-}{\;})
 \rtimes \lambda' $ is the
 intertwining
operator for
$ \mid \; \mid^{\alpha} \chi \rtimes \lambda' $ corresponding to  $ w. \; \mathcal{O} $ is the ring of integers of $ E. $
For $ \alpha \in \R^*_{-}, $ one can equivalently define such an
intertwining operator.

\bigskip
 
Let $ \G $ be a connected reductive group defined over a p-adic field. Let $ (\pi,V) $ be a representation of $ \G, $
for a finite dimensional vector space $ V. $
One has the following construction:

\begin{Lemma}
\label{cont}  Assume one has, on the same vector space $ V, $ a continuous family of induced irreducible representations
$ (\pi_{\alpha},  V), \alpha \in X, $ where $ X $ is a connected set, that posses non-trivial hermitian forms
(invariant under $ \G $). Suppose that some $ \pi_{\alpha} $ is unitary. If a family of non-degenerate
hermitian forms on a finite dimensional space, parametrised by $ X, $ is positive definite at one point of $ X, $ it is
positive definite everywhere. Hence $ \pi_{\alpha} $ is unitary $ \forall \alpha \in X. $
\end{Lemma}

\begin{Remark}
One may reduce the argument to finite dimensional spaces by considering spaces $ \oplus V(\delta), $ where $ \delta $
runs over fixed finite subsets of the irreducible unitary representations of the maximal compact subgroup $ U(3, \mathcal{O}). $
\end{Remark}

\medskip
We continue the proof of \ref{unitaryu3}.

1.2 and 1.3 $ \Leftarrow $
For $ \alpha = 0 $ and $ \chi = 1 $ or $ \chi \in X_{\omega_{E/F}} \; \Ind_P^{U(3)}(\lambda) $ is irreducible and 
unitary, hence $ \Ind_P^{U(3)}(\lambda) $ is unitary for $ \chi = 1 $ and $ \alpha \in ]-1,1[ $ and for $ \chi \in
X_{\omega_{E/F}} $ and $ \alpha \in ] - 1/2, 1/2 [. $ The hermitian forms are given above.

\bigskip

If $ \alpha \neq 0 $ and $ \chi \notin X_{N_{E/F}(E^*)}, \; \Ind_P^{U(3)}(\lambda) $ is irreducible and not hermitian and
hence not unitarisable.

It remains to show that $ \Ind_P^{U(3)}(\lambda) $ is non-unitary if $ \alpha \in \R^* $ and $ \chi \in 
X_{1_{F^*}} $ (if $ \alpha = 0 $ and $ \chi \in X_{1_{F^*}} $ then $ \Ind_P^{U(3)}(\lambda) $ is reducible).
Further it remains to show that $ \Ind_P^{U(3)}(\lambda) $ is non-unitary if $ \chi = 1 $ and 
$ \alpha \in ]- \infty, -1 [ \cup
 ] 1 , \infty [  $ and if $ \chi \in X_{\omega_{E/F}} $ and $ \alpha \in ] - \infty, - 1/2 [ \cup ] 1/2 , \infty [. $

This will show 1. $ \Rightarrow $ and 2. $ \Leftarrow; \; 2. \Rightarrow $ is shown by $ 1. \Leftarrow. $

\bigskip

Let $ F $ be a non-archimedean local field of characteristic $ 0 $ and let $ \G $ be an arbitrary reductive group defined over
 $ F. $
Let $ \sigma \in \Hom(M, \C^*)^{n.r.}, $ where $  \Hom(M, \C^*)^{n.r.} $ is the group of non-ramified characters of $M.$ It is
canonically isomorphic (as a topological group) to a direct product of a finite number of copies of $ \C^*. $

Let $ \pi $ be an irreducible representation of $ M. $

\smallskip

We will use Lemma 5.1 (i) of \cite{MR1480543} that is a special case of Theorem 4.5 in \cite{MR963153}:

\smallskip

\begin{Lemma}
\label{comp} The set of all $ \sigma \in  \Hom(M, \C^*)^{n.r.} $ such that $ \Ind_P^{\G}(\sigma \otimes \pi) $
has an irreducible unitary subquotient is compact.
\end{Lemma}

\medskip

Here $ \mid \; \mid^{\alpha} \otimes \; 1 \; \in  \Hom(M, \C^*)^{n.r.} $ and $ \chi \otimes \lambda' $ is an irreducible representation
of $ M. $

\medskip

Proof of 1. $ \Rightarrow $ and  2. $ \Leftarrow: $
 $ \Ind_P^{U(3)}(\lambda) $ is irreducible
for $ \chi = 1 $ and $ \alpha \in ] 1, \infty [ $ (or $ \alpha \in ] - \infty , -1 [, $ or for $  \chi \in X_{\omega_{E/F}}
$ and $ \alpha \in ] - \infty, -1/2 [ \cup
 ] 1/2, \infty [, $  or for $ \chi \in X_{1_{F^*}} $ and $ \alpha \in \R^*,$
 respectively). If there existed $ \alpha \in ] 1 , \infty [ $ (or in one
of the
other intervalls
or in $ \R^*, $ respectively) s.t. $ \mid \; \mid^{\alpha} \chi \rtimes \lambda' $ is unitary, with $ \chi $ chosen
appropriately, then by Lemma \ref{cont} all 
representations $ \mid \; \mid^{\alpha} \chi \rtimes \lambda' $ with $ \alpha \in ] 1, \infty [ $ (or in one of the 
other intervalls or
in $ \R^* $)
would be unitary, in contradiction to Lemma \ref{comp}.
\end{Proof} 

The induced representations of $ U(4) $ over a p-adic field have been classified by Kazuko Konno \cite{Ko}.

\bigskip

\section{The irreducible representations of $ U(5) $}

\begin{subsection}{Levi decomposition for $ U(5) $}

\bigskip

Recall the Levi decompositon $ P = MN, $ where $ P $ is a standard parabolic subgroup, $ M $ is the standard Levi-subgroup
corresponding to $ P $ and $ N $ is 
the
unipotent subgroup corresponding to $ P $ and to $ M. $

\medskip

The standard Levi subgroups of $ U(5) $ are the following three:

\medskip

$ M_0 :=  E^* \times E^* \times E^1 $ (the Levi-group corresponding to the minimal parabolic subgroup),

\medskip

$ M_1 := \GL(2,E) \times E^1 $ and

\medskip

$ M_2:= E^* \times U(3) $ (the two Levi-groups corresponding to the maximal parabolic subgroups).

\bigskip

In matrix form we can write:

$ M_0 = \{
\left( \begin{smallmatrix}
x&&&&0\\
&y&&&\\
&&k&&\\
&&&\overline{y}^{-1}&\\
0&&&&\overline{x}^{-1}
\end{smallmatrix} \right), x,y \in E^*, k \in E^1 \}, $

$ M_1 = \{
\left( \begin{smallmatrix}
a&&0\\
&k&\\
0&&\overline{a}^{-1}
\end{smallmatrix} \right), \; a \in \GL(2,E), k \in E^1 \cong U(1) \} $ and

$ M_2 = \{
\left( \begin{smallmatrix}
 x&&&&0\\
&u&\\
0&&&&\overline{x}^{-1}
\end{smallmatrix} \right), x \in E^* \cong \GL(1,E), u \in U(3) \}. $

\bigskip

The unipotent subgroups are

$ N_0 = \{
\left(
\begin{smallmatrix}
1&&&-&*\\
&1&&&\mid\\
&&1&&\\
\mid&&&1&\\
0&-&&&1
\end{smallmatrix}
\right) \} \cap U(5), $
$ N_1 = \{
\left(
\begin{smallmatrix}
 1&0&&-&*\\
0&1&&&\mid\\
&&1&&\\
\mid&&&1&0\\
0&-&&0&1
\end{smallmatrix}
\right) \} \cap U(5), $ and
$ N_2 = \{
\left(
\begin{smallmatrix}
1&&&-&*\\
&1&0&0&\mid\\
&0&1&0&\\
\mid&0&0&1&\\
0&-&&&1
\end{smallmatrix}
\right) \} \cap U(5), $
where $ * $ are entries in $ E. $

\bigskip

We obtain the parabolic subgroups

$ P_0 = M_0 N_0 = \{
\left(
\begin{smallmatrix}
 x&&&&*\\
&y&&&\\
&&k&&\\
&&&\overline{y}^{-1}&\\
0&&&&\overline{x}^{-1}
\end{smallmatrix}
\right), \; x,y, \in E^*, k \in E^1, * \in E  \} \cap U(5), $

$ P_1 = M_1 N_1 = \{
\left(
\begin{smallmatrix}
a&&&*\\
&&k&&\\
0&&&\overline{a}^{-1}
\end{smallmatrix}
\right), \; a \in \GL(2,E), k \in E^1, * \in E \} \cap U(5), $ and

$ P_2  = M_2 N_2 = \{
\left(
\begin{smallmatrix}
 x&&&*\\
&u&\\
0&&&\overline{x}^{-1}
\end{smallmatrix}
\right), x \in E^*, u \in U(3), * \in E \} \cap U(5). $

\bigskip

We consider representations of the Levi-subgroups. One extends them to representations of $ P, $ by extending them
trivially to the unipotent subgroup $ N. $ Then one performs normalized parabolic induction to the whole group $ U(5). $
\end{subsection}

\begin{subsection}{Representations with cuspidal support in $ M_0, $ fully-induced}

The irreducible representations of $ M_0 $ are characters.
Let $ \lambda_1, \lambda_2
\in \Hom(E^*, \C^*) $ and $ \lambda' \in \Hom(E^1, \C^*) $ be smooth characters.
One may write $ \lambda_i = \mid \; \mid_E^{\alpha_i} \chi_i, \; i = 1,2, $ where $ \alpha_i \in \R $ and $ \chi_i $
is a unitary character of $ E^*. \; \lambda' $ is unitary.

\bigskip

Then each character $ \lambda $ of $ M_0 $ can be written as

$$ \lambda(m) = \mid x \mid_E^{\alpha_1} \chi_1(x) \mid y \mid_E^{\alpha_2} \chi_2(y) 
\lambda'(x\overline{x}^{-1}y \overline{y}^{-1}
k), \; m =
\left(
\begin{smallmatrix}
 x&&&&0\\
&y&&&\\
&&k&&\\
&&&\overline{y}^{-1}&\\
0&&&&\overline{x}^{-1}
\end{smallmatrix}
\right), \; x,y \in E^*, \; k \in E^1.
 $$

\bigskip

By $ \lambda: = \lambda_1 \otimes \lambda_2
 \otimes \lambda' $ we denote the characters of $ M_0 $ and by $ \lambda_1 \times \lambda_2 \rtimes \lambda': = 
\Ind_P^{U(5)}(\lambda_1 \otimes \lambda_2 \otimes \lambda') $ the induced representations to $ U(5). $

\bigskip

We start with the case where $ \lambda_1 = \chi_1 $ and $ \lambda_2 = \chi_2 $ are unitary characters, i.e.
$ \alpha_1 = \alpha_2 = 0. $

\bigskip
\smallskip

\begin{center}{4.2.1 \textbf{Irreducible subquotients of} $ \chi_1 \times \chi_2 \rtimes \lambda' $} \end{center}

\bigskip

Let $ P_0 $ be the minimal parabolic subgroup of $ U(5) $ (the upper triangular matrices in $ U(5) $) with Levi subgroup
$ M_0 $ and unipotent subgroup $ N_0, $ such that $ P_0 = M_0 N_0). $

\bigskip
 
\begin{Theorem}
\label{chi12}
Let $ \chi_1, \chi_2 $ be unitary characters of $ E^* $ and let $ \lambda' $ be a (unitary) character 
of $ E^1. $

The induced representation $ \chi_1 \times \chi_2 \rtimes \lambda' $ is reducible if and only if

\bigskip

$ \exists \;  i \in \{ 1,2 \} \; \text{s.t.} \;  \chi_i \in X_{1_{F^*}}. $

\end{Theorem}

\begin{Proof} For the proof we use the theory of $ R  - $ groups, these are subgroups of the Weyl group W of $ U(5). $

\bigskip

We have $ W \cong S_2 \times (\Z/2 \Z)^2, $ where $ S_2 \cong \Z / 2 \Z $ is the symmetric group in 2 letters.

\medskip

Let $ \lambda $ be a character of $ M_0, $ and $ W_{\lambda} := \{ w \in W: w \lambda = \lambda \}. $

Let $ a(w, \lambda): \Ind_P^{U(5)}(\lambda) \rightarrow \Ind_P^{U(5)}(w\lambda) $ be the intertwining operator of 
$ \Ind_P^{\G}(\lambda) $ corresponding to $ w, $ where
$ w \lambda(m) := \lambda(w^{-1}mw). $

Let $ W':= \{  w \in W_{\lambda} : a(w, \lambda) \; \text{is scalar} \}. $

We have $ W_{\lambda} = R \ltimes W'. $
By a result of D. Keys, the number of inequivalent irreducible components of $ \Ind_{P_0}^{U(5)}(\lambda) $ equals the
number of conjugacy classes in 
$ R, $ see Corollary 1 in \cite{MR675406}.

\bigskip

In order to make use of the Theorem 3.4 in \cite{1995} we
introduce
some notation:

\bigskip
Let $ \G: = \left[^{U(2n)}_{U(2n+1)} \right. $ be the unitary group in $ 2n $ or $ 2n + 1 $
variables, respectively.

For $ m \leq n $ let $ \G(m) := \left[^{U(2m) \text{\; if } \G = U(2n)}_{ U(2m+1) \text{ \; if \; } \G = U(2n+1)} 
\right.. $

By convention $ \G(0) = U(1). $
Let $ \epsilon_2(\G) $ be the set of equivalence classes of irreducible square-integrable representations of $ \G. $
Let $ \sigma_i \in \epsilon_2(\GL_{n_i}(E)), \; i = 1,2, \ldots $ and $ \rho \in \epsilon_2(\G(m)). $

\medskip

\begin{Theorem}[\cite{1995}] Let $ \G = U(2n) $ or $ U(2n+1). $
Let
$ P=MN $ be a parabolic subgroup of $ \G. $ Suppose that $ M \cong
\GL(n_1,E) \times \ldots \times \GL(n_r,E) \times \G(m). $
Let $ \sigma \in \epsilon_2(M), $ with $ \sigma \cong \sigma_1 \otimes \ldots \otimes \sigma_r
\otimes \rho. $
Let $ d $ be the number of inequivalent $ \sigma_i, $ such that
$ \Ind_{\GL_{n_i} \times \G(m)}^{\G(m+n_i)}(\sigma_i
\otimes \rho) $ reduces. Then R  $ \cong (\Z/2 \Z)^d. $
\end{Theorem}

\bigskip

In the present case this translates to:

\medskip

$ \G = U(5) $, $ P = P_0 $ is the minimal parabolic subgroup, $ M = M_0 \cong \GL_1(E) \times \GL_1(E) \times G(0)
\cong E^* \times E^* \times E^1. $

\smallskip

$ \sigma_1 = \chi_1, \; \sigma_2 = \chi_2, $ and $ \rho = \lambda'. $

\medskip

Recall that for a unitary character $ \chi $ of $ E^*, \; \chi \rtimes \lambda' $ is 
reducible if and only if
$ \chi \in X_{1_{F^*}}. $ Then $ \chi \rtimes \lambda' = \sigma_{1, \chi} \oplus \sigma_{2, \chi}, $ where $ \sigma_{1, \chi} $
and $ \sigma_{2, \chi} $ are tempered.

\bigskip

We apply the theorem:
for $ \lambda = \chi_1 \otimes \chi_2 \otimes \lambda' $ and $ W \cong S_2 \times (\Z/2 \Z)^2 $ the integer $ d $ may equal $ 0,1 $ or $ 2.$

\medskip

0. Let $ d = 0. \; \chi_i \rtimes \lambda', \; i \in \{ 1,2 \} \; $ is irreducible for $ i \in \{ 1,2 \}, $ and
$ R \cong \{ 1 \}. $
$ \chi_1 \times \chi_2 \rtimes \lambda' $ is irreducible and unitary.

\bigskip

1. Let $ d = 1. $ Then $ \exists \; i,j \in \{ 1,2 \}, i \neq j, $ s.t. $ \chi_i \in X_{1_{F^*}} $ and $ \chi_j \notin
X_{1_{F^*}} $ or $ \chi_i \in X_{1_{F^*}} $ and $ \chi_j \cong \chi_i. $
Hence $ R \cong \Z/2\Z, $ and $  \chi_1 \times \chi_2 \rtimes \lambda' $ has two irreducible 
inequivalent constituents: 
$ \chi_j \rtimes \sigma_{1, \chi_i} $ and 
$ \chi_j \rtimes \sigma_{2, \chi_i}. $  They are tempered and hence unitary.

\bigskip

2. Let $ d = 2.  \;  \chi_1 $ and $ \chi_2 $ are two inequivalent characters and $ \chi_i \in X_{1_{F^*}} $ for 
$  i = 1,2. $

$ R \cong (\Z / 2 \Z)^2, $ and $ \chi_1 \times \chi_2 \rtimes \chi' $ has four irreducible 
inequivalent unitary constituents. By \cite[Theorem 4.3]{1995} they are tempered and elliptic.
\end{Proof}

\bigskip
\smallskip

\begin{center}{4.2.2 \textbf{Irreducible subquotients of}
 $ \mid \; \mid^{\alpha_1} \chi_1 \times \mid \; \mid^{\alpha_2} \chi_2
\rtimes
 \lambda', \; \alpha_1, \alpha_2 > 0 $ \textbf{and of} $ \mid \; \mid^{\alpha} \chi_1 \times \chi_2 \rtimes \lambda', \; 
\alpha > 0 $}\end{center}

\bigskip

Let $ M_0 \cong E^* \times E^* \times E^1 $ be the minimal Levi subgroup, let
$ P_0 = M_0 N_0 $ be the corresponding parabolic subgroup, i. e. the minimal parabolic subgroup.
$ N_0 $ is the
unipotent radical of $ P_0. $

Let $ \lambda := \lambda_1 \otimes \lambda_2 \otimes \lambda' = \mid \; \mid^{\alpha_1} \chi_1 \otimes 
\mid \; \mid^{\alpha_2} \chi_2
\otimes \lambda' $ be a character of $ M_0, $ where $ \alpha_1, \alpha_2 \in \R $ and $ \chi_1, \chi_2 $ are unitary
characters of $ E^*. $

\medskip

W. l . o. g. let $ \alpha_1 \geq  \alpha_2 > 0. $ If $ \alpha_2 > \alpha_1 $ we change parameters.
The case $ \alpha_2 = 0 $ is treated
seperately.

Recall that $ X_{\omega_{E/F}} $ is the set of characters of $ E^* $ whose restriction to $ F^* $ is the character 
$ \omega_{E/F}, $
i.e. whose restriction to $ F^* $ is non-trivial on $ F^* $ but trivial on $ N_{E/F}(E^*). \; X_{1_{F^*}} $ is the set 
of non -
trivial characters of $ E^* $ whose restriction to $ F^* $ is trivial. $ X_{N_{E/F}(E^*)} =  {1} \cup X_{\omega_{E/F}}
\cup X_{1_{F^*}}. \; X_{N_{E/F}(E^*)} $ contains all characters $ \chi $ of $ E^* $ satisfying 
$ \chi(x) = \chi^{-1}(
\overline{x}). $ They are unitary.

\medskip

From now on, the lack of an entry at position $ ij $ in a matrix means that the entry equals 0.

\bigskip

\begin{Theorem}
\label{alpha12chi12}

Let $ \chi_1, \chi_2 $ be unitary characters of $ E^* $  and let $ \lambda' $ be a character of $ E^1. $
Let $ \alpha_1, \alpha_2 \in \R_+^* $ s.t. $ \alpha_1 \geq \alpha_2, $ then

$$  \mid \; \mid^{\alpha_1}  \chi_1 \times  \mid \; \mid^{\alpha_2} \chi_2 \rtimes \lambda' $$ is reducible if and only if
at least one of the following conditions holds:

\bigskip

1. $ \alpha_1 - \alpha_2 =  1 $ and $ \chi_1 = \chi_2,  $ 

2. $ \alpha_1 + \alpha_2 =  1 $ and $ \chi_1(x) = \chi_2^{-1}(\overline{x}), $

3. $ \alpha_1 =  1 $ and $ \chi_1 = 1 $ or $ \alpha_1 =  1/2 $ and $ \chi_1 \in X_{\omega_{E/F}}, $

4. $ \alpha_2 = 1 $ and $ \chi_2 = 1 $ or $ \alpha_2 =  1/2 $ and $ \chi_2  \in X_{\omega_{E/F}}. $

\end{Theorem}

\bigskip

\begin{Proof}
Let $ \lambda := \mid \; \mid^{\alpha_1} \chi_1 \otimes \mid \; \mid^{\alpha_2} \chi_2 \otimes \lambda' $ be a 
character of $ M_0. $

Let $ A(w, \lambda): \; \Ind_{P_0}^{U(5)}(\lambda) = \mid \; \mid^{\alpha_1} \chi_1 \times \mid \; \mid^{\alpha_2} 
\chi_2 \rtimes \lambda' \rightarrow 
\Ind_{P_0}^{U(5)}(w\lambda) = \mid \; \mid^{-\alpha_1} \chi_1^{-1}(\overset{-}{\;}) \times \mid \; 
\mid^{-\alpha_2} 
\chi_2^{-1}(\overset{-}{\;}) \rtimes \lambda'
$ be a standard long intertwining operator for the representation $ \mid \; \mid^{\alpha_1} \chi_1 \times 
\mid \; \mid^{\alpha_2} \chi_2 \rtimes \lambda'. $ 

\bigskip

\begin{Remark}

$  w =
\left(
\begin{smallmatrix}
&&&&1\\
&&&1&\\
&&1&&\\
&1&&&\\
1&&&&
\end{smallmatrix}
\right) $
is the longest element of the Weyl group, and for $ m \in M_0 $ it is

$ w \lambda(m) :=
\lambda(\left(
\begin{smallmatrix}
&&&&1\\
&&&1&\\
&&1&&\\
&1&&&\\
1&&&&
\end{smallmatrix}
\right)
\left(
\begin{smallmatrix}
 x&&&&\\
&y&&&\\
&&k&&\\
&&&\overline{y}^{-1}&\\
&&&&\overline{x}^{-1}
\end{smallmatrix}
\right)
\left(
\begin{smallmatrix}
&&&&1\\
&&&1&\\
&&1&&\\
&1&&&\\
1&&&&
\end{smallmatrix}
\right)) = \lambda(
\left(
\begin{smallmatrix}
\overline{x}^{-1}&&&&\\
&\overline{y}^{-1}&&&\\
&&k&&\\
&&&y&\\
&&&&x
\end{smallmatrix}
\right)). $
\end{Remark}

Hence $ \Ind_{P_0}^{U(5)}(w\lambda) $ equals $ \mid \; \mid^{- \alpha_1} \chi_1^{-1}(\overset{-}{\;}) \times \mid \; \mid^{- \alpha_2}
\chi_2^{-1}(\overset{-}{\;}) \rtimes \lambda'. $

\bigskip
\medskip

If $ A(w, \lambda) $ is either not injective or not surjective it follows that $ \Ind_{P_0}^{U(5)}(\lambda) $ is 
reducible.
The decomposition of the long intertwining operator into short operators shows for which $ \alpha_1, \alpha_2 \;
 $ 
and unitary characters $ \chi_1 $ and $ \chi_2 $ the long intertwining operator is not an isomorphism.

\bigskip

Let
$ w_1 :=
\left(
\begin{smallmatrix}
 0&1&&&\\
1&0&&&\\
&&1&&\\
&&&0&1\\
&&&1&0
\end{smallmatrix}
\right) $
and
$ w_2 :=
\left(
\begin{smallmatrix}
1&&&&\\
&&&1&\\
&&1&&\\
&1&&&\\
&&&&1
\end{smallmatrix}
\right). $

\bigskip

We have $ w = 
\left(
\begin{smallmatrix}
&&&&1\\
&&&1&\\
&&1&&\\
&1&&&\\
1&&&&
\end{smallmatrix}
\right) =
\left(
\begin{smallmatrix}
0&1&&&\\
1&0&&&\\
&&1&&\\
&&&0&1\\
&&&1&0
\end{smallmatrix}
\right) 
\left(
\begin{smallmatrix}
1&&&&\\
&&&1&\\
&&1&&\\
&1&&&\\
&&&&1
\end{smallmatrix}
\right)
\left(
\begin{smallmatrix}
0&1&&&\\
1&0&&&\\
&&1&&\\
&&&0&1\\
&&&1&0
\end{smallmatrix}
\right)
\left(
\begin{smallmatrix}
1&&&&\\
&&&1&\\
&&1&&\\
&1&&&\\
&&&&1
\end{smallmatrix}
\right) = w_1 w_2 w_1 w_2
$.

\bigskip

The following diagram gives the decomposition of $ A(w, \lambda). $

\medskip

$
\begin{smallmatrix}
 \mid \; \mid^{\alpha_1} \chi_1 \times \mid \; \mid^{\alpha_2} \chi_2 \rtimes \lambda' & \overset{\cong}\rightarrow & 
\mid \; \mid^{\alpha_1} \chi_1 \times \mid \; \mid^{\alpha_2} \chi_2 \rtimes \lambda'\\
&& A(w_2, \lambda) \downarrow w_2\\
&&\mid \; \mid^{\alpha_1} \chi_1 \times \mid \; \mid^{-\alpha_2} \chi_2^{-1}(\overset{-}{\;}) \rtimes \lambda'\\
A(w, \lambda) \downarrow w && A(w_1, w_2 \lambda) \downarrow w_1\\
&&\mid \; \mid^{-\alpha_2} \chi_2^{-1}(\overset{-}{\;}) \times \mid \; \mid^{\alpha_1} \chi_1 \rtimes \lambda'\\
&& A(w_2, w_1 w_2 \lambda) \downarrow w_2\\
&&\mid \; \mid^{-\alpha_2} \chi_2^{-1}(\overset{-}{\;}) \times \mid \; \mid^{-\alpha_1} \chi_1^{-1}(\overset{-}{\;}) \rtimes \lambda'\\
&& A(w_1, w_2 w_1 w_2 \lambda) \downarrow w_1\\
\mid \; \mid^{-\alpha_1} \chi_1^{-1}(\overset{-}{\;}) \times \mid \; \mid^{-\alpha_2} \chi_2^{-1}(\overset{-}{\;}) \rtimes 
\lambda'&\overset{\cong}{\rightarrow}&\mid \; \mid^{-\alpha_1} \chi_1^{-1}(\overset{-}{\;}) \times \mid \; \mid^{-\alpha_2}
 \chi_2^{-1}(\overset{-}{\;}) \rtimes \lambda'
\end{smallmatrix} $

\bigskip

If $ A(w, \lambda) $ is not an isomorphism, then at least one of the operators $ A(w_2 \lambda),
A(w_1, w_2 \lambda),
A(w_2, w_1 w_2 \lambda) $ or $ A(w_1, w_2 w_1 w_2 \lambda) $ is not an isomorphism.

$ A(w_1, \lambda) $ is no isomorphism if and only if the induced representation 
$ \mid \; \mid^{\alpha_2} \chi_2 \rtimes
\lambda' $ is reducible. This is the case if and only if
$ \alpha_2 =  1 $ and $ \chi_2 = 1 $ or $ \alpha_2 =  1/2 $ and $  \chi_2  \in X_{\omega_{E/F}}. $

$ A(w_1, w_2 \lambda) $ is no isomorphism if and only if the corresponding representation $ \mid \; \mid^{\alpha_1}
\chi_1 \times \mid \; \mid^{-\alpha_2} \chi_2^{-1}(\overset{-}{\;}) $ is reducible if and only if $ \alpha_1 + \alpha_2 = 1 $ and 
$ \chi_1(x) = \chi_2^{-1}(\overline{x}) \; \forall x \in E^*.$ 

$ A(w_2, w_1 w_2 \lambda) $ is no isomorphism if and only if $
 \mid \; \mid^{\alpha_1} \chi_1 \rtimes \lambda' $ is
reducible if and only if $ \alpha_1 =  1 $ and $ \chi_1 = 1 $ or 
$ \alpha_1 =  1/2 $ and $ \chi_1 \in X_{\omega_{E/F}}, $ and

$ A(w_1, w_2 w_1 w_2 \lambda) $ is no isomorphism if and only if  $
 \mid \; \mid^{-\alpha_2} \chi_2^{-1}(\overset{-}{\;}) \times \mid \; \mid^{-\alpha_1} \chi_1^{-1}(\overset{-}{\;}) $ is reducible
if and only if $ \alpha_1 - \alpha_2 =  1 $ and $\chi_1 = \chi_2. $

\bigskip

In all other cases the short intertwining operators are holomorphic and isomorphisms, hence $ A(w, \lambda) $
is an isomorphism and the representation
$ \mid \; \mid^{\alpha_1} \chi_1 \times \mid \; \mid^{\alpha_2} \chi_2 \rtimes \lambda' $
is irreducible.

\bigskip

On the other hand, if at least one of the short intertwining operators is no isomorphism,
$ \mid \; \mid^{\alpha_1} \chi_1 \times \mid \; \mid^{\alpha_2} \chi_2 \rtimes \lambda' $ is reducible by indcution in
stages; in these cases 
we determine the irreducible constituents in Theorems \ref{alphachist}, \ref{alphachistu3},
\ref{alphachipi1}, \ref{2111}, \ref{1/21/2chiomega}, \ref{3/21/2chiomega} and \ref{11/21chiomega}.
\end{Proof}

Let $ \alpha_1 > \alpha_2 = 0. $

\begin{Theorem}
\label{alpha1chi12} Let $ \chi_1, \chi_2 $ be unitary characters of $ E^*, $ let $ \lambda' $ be a (unitary) character of
$ E^1. $ Let $ \alpha_1 \in \R_+^*. $ The induced representation

$$ \; \mid \; \mid^{\alpha_1} \chi_1 \times \chi_2 \rtimes \lambda' $$ is reducible if and only if

\bigskip
1. $ \alpha_1 = 1 $ and $ \chi_1 = \chi_2, $

2. $ \alpha_1 = 1 $ and $ \chi_1(x) = \chi_2^{-1}(\overline{x}{\;}), $

3. $ \alpha_1 = 1 $ and $ \chi_1 = 1 $ or $ \alpha_1 = 1/2 $ and $ \chi_1 \in X_{\omega_{E/F}}, $

4. $ \chi_2 \in X_{1_{F^*}}. $

\end{Theorem}

\bigskip

\begin{Proof}
1. $ \Rightarrow $
We use Lemma 2.1 of \cite{MR2504024}.

\smallskip
\begin{Lemma}[\cite{MR2504024}]
\label{Tadicalphachi12}

 Let $ \pi $ be an irreducible representation of $ U(m) $ and let $  \rho $ be an irreducible cuspidal
representation of a general linear group $ \GL(p,F). $ Suppose

(1) $ \rho \neq \overset{\sim}{\rho}(\overset{-}{\;}). $

(2) $ \rho \rtimes \pi_{\cusp} $ is irreducible.

(3) $ \rho \times \rho' $ and $ \overset{\sim}{\rho}(\overset{-}{\;}) \times \rho' $ are irreducible for any factor $ \rho' $ of $ \pi. $

(4) Neither $ \rho $ nor $ \overset{\sim}{\rho}(\overset{-}{\;}) $ is a factor of $ \pi. $

Then $$ \rho \rtimes \pi $$ is irreducible.
\end{Lemma}

Here $ \pi \cong \chi_2 \rtimes \lambda' $ and $ \; \rho \cong \; \mid \; \mid^{\alpha_1} \chi_1. $

If none of the four cases in Theorem \ref{alpha1chi12} holds we are in the position to apply Lemma \ref{Tadicalphachi12},
hence $ \mid \; \mid^{\alpha_1} \chi_1
 \times \chi_2 \rtimes \lambda' $
is irreducible.

$ \Leftarrow $ If at least one of the four cases holds, clearly $ \mid \; \mid^{\alpha_1} \chi_1 \times \chi_2 \rtimes
\lambda' $ is reducible.
\end{Proof}

In those cases where $ \mid \; \mid^{\alpha_1} \chi_1 \times \chi_2 \rtimes \lambda', \; \alpha > 0, $ is reducible, the
irreducible constituents will be investigated \ref{alphachist}, \ref{chistu3}, \ref{chipi1} and in \ref{alphachisigma12}.

\end{subsection}

\subsection{Representations induced from $ M_1 $ and $ M_2, $ with cuspidal support in $ M_0 $}
 
We consider representations induced from the maximal parabolic subgroups with Levi-groups $ M_1 $ and $ M_2, $ whose
cuspidal support is in $ M_0 $ and that are not fully induced. We begin with $ M_1 = \GL(2,E) \times E^1. $ 

\bigskip

\begin{center}{4.3.1 \textbf{Representations} $ \mid \; \mid^{\alpha} \chi \St_{\GL_2} \rtimes \lambda' $ \textbf{and} $ 
 \mid \; \mid^{\alpha} \chi 1_{\GL_2} \rtimes \lambda', \; \alpha > 0 $}\end{center}

Let $ \alpha \in \R_+^* $ and $ \chi $  be a unitary character of $ E^*. $
We study $ \mid \; \mid^{\alpha} \chi
\St_{\GL_2} \rtimes \lambda' $  that is a subrepresentation
of the induced
representation $ \mid \; \mid^{\alpha + 1/2} \chi \times \mid \; \mid^{\alpha - 1/2} \chi \rtimes \lambda', $ and its
Aubert dual
$ \mid \; \mid^{\alpha} \chi 1_{\GL_2} \rtimes \lambda'. $ 
 
\bigskip

\begin{Theorem}
\label{alphachist}

Let $ \alpha \in \R_+^* $ and $ \chi $ be a unitary character of $ E^*. $ 
The representations $ \mid \; \mid^{\alpha} \chi(\det) \St_{\GL_2} \rtimes \lambda' $ and
$ \mid \; \mid^{\alpha} \chi(\det) 1_{\GL_2} \rtimes \lambda' $ are irreducible, unless one of the following cases holds:

1. $ \alpha = 1/2 $ and  $ \chi \in X_{N_{E/F}(E^*)}, $

2. $ \alpha = 3/2 $ and $ \chi = 1, $

3. $ \alpha = 1 $ and $ \chi \in X_{\omega_{E/F}}. $

\end{Theorem}

\begin{Proof}
Let $ R(U(n)) $ be the Grothendieck group of the
category of admissible representations of finite length of $ U(n) $ and let
$ R(U) := \underset{n \geq 0}{\oplus} R(U(n)). $

In $ R(U) $ we have $ \mid \; \mid^{\alpha + 1/2} \chi \times \mid \; \mid^{\alpha - 1/2} \chi \rtimes \lambda' =
\mid \; \mid^{\alpha} \chi \St_{\GL_2} \rtimes \lambda' + \mid \; \mid^{\alpha} \chi 1_{\GL_2} \rtimes \lambda'. $ 

$ \mid \; \mid^{\alpha} \chi \St_{\GL_2} \rtimes \lambda' $ and 
$ \mid \; \mid^{\alpha} \chi 1_{\GL_2} \rtimes \lambda' $ are dual in the sense of Aubert and 
have the same number of irreducible constituents.
We give the proof for $ \mid \; \mid^{\alpha} \chi \St_{\GL_2} \rtimes \lambda' $ as
subrepresentation of
$ \mid \; \mid^{\alpha + 1/2} \chi \times \mid \; \mid^{\alpha - 1/2} \chi \rtimes \lambda'. $

\medskip

Let $ \lambda: = \mid \; \mid^{\alpha + 1/2} \chi \otimes \mid \; \mid^{\alpha -1/2} \chi \otimes \lambda'. $
Let $ A(w', \lambda): \; \mid \; \mid^{\alpha} \chi  \St_{\GL_2} \rtimes \lambda' \rightarrow
\mid \; \mid^{-\alpha} \chi^{-1}(\overline{\det}) \St_{\GL_2} \rtimes \lambda' $ be the long intertwining operator
for the representation $ \mid \; \mid^{\alpha} \chi \St_{\GL_2} \rtimes \lambda', $ where $ w' $ is the
longest element of $ W $ respecting $ M_1 \cong \GL(2,E) \times E^1. $

\medskip

We have
$ w' :=
\left(
\begin{smallmatrix}
&&&1&0\\
&&&0&1\\
&&1&&\\
1&0&&&\\
0&1&&&
\end{smallmatrix}
\right) =
\left(
\begin{smallmatrix}
1&&&&\\
&&&1&\\
&&1&&\\
&1&&&\\
&&&&1
\end{smallmatrix}
\right)
\left(
\begin{smallmatrix}
0&1&&&\\
1&0&&&\\
&&1&&\\
&&&0&1\\
&&&1&0
\end{smallmatrix}
\right)
\left(
\begin{smallmatrix}
1&&&&\\
&&&1&\\
&&1&&\\
&1&&&\\
&&&&1
\end{smallmatrix}
\right) =
w_2 w_1 w_2.
$ 

\medskip

The decomposition of $ A(w',\lambda) $ into short intertwining operators 
gives information for which $ \alpha > 0 $ and unitary characters $ \chi $ of $ E^* $
this operator is an isomorphism.
The following diagram shows the decomposition of $ A(w', \lambda). $
 $ i_1 $ and $ i_2 $ are inclusions that depend holomorphically on $ \alpha. $

\bigskip

$\begin{smallmatrix}
  \mid \; \mid^{\alpha} \chi (\det) \St_{\GL_2} \rtimes \lambda'& \overset{i_1}{\hookrightarrow}& \mid \; \mid^{\alpha + 1/2} \chi 
\times \mid \; \mid^{\alpha - 1/2} \chi \rtimes \lambda'\\
&& A(w_2, \lambda) \downarrow w_2\\
&& \mid \; \mid^{\alpha + 1/2} \chi \times \mid \; \mid^{-\alpha + 1/2} \chi^{-1}(\overset{-}{\;}) \rtimes \lambda'\\
A(w', \lambda) \downarrow w' && A(w_1, w_2 \lambda) \downarrow w_1\\
&& \mid \; \mid^{-\alpha + 1/2} \chi^{-1}( \overset{-}{\;}) \times \mid \; \mid^{\alpha + 1/2} \chi \rtimes \lambda'\\
&& A(w_2, w_1 w_2 \lambda) \downarrow w_2 \\
\mid \; \mid^{-\alpha} \chi^{-1}(\overline{\det}) \St_{\GL_2} \rtimes \lambda' & \overset{i_2}{\hookrightarrow} & 
\mid \; \mid^{-\alpha + 1/2} \chi^{-1}(\overset{-}{\;}) \times \mid \; \mid^{-\alpha - 1/2} \chi^{-1} (\overset{-}{\;})
 \rtimes \lambda'\\
 \end{smallmatrix}
$

\bigskip

If $ \alpha \neq 1/2 \; A(w_2, \lambda) $ is no isomorphism if and only if $ \mid \; \mid^{\alpha - 1/2} \chi \rtimes 
\lambda' $ reduces,
if and only if
$ \alpha = 3/2 $ and $ \chi = 1 $ or $ \alpha = 1 $ and $ \chi \in  X_{\omega_{E/F}}. $

If $ \alpha = 1/2 $ and $ \chi \in X_{1_{F^*}}, $ then $ \chi \rtimes \lambda' $ reduces.

\medskip

$ A(w_1, w_2 \lambda) $ is no isomorphism if and only if $ \mid \; \mid^{\alpha + 1/2} \chi \times 
\mid \; \mid^{-\alpha + 1/2} \chi^{-1}(\overset{-}{\;})
$ reduces, if and only if $ \alpha = 1/2 $ and $ \chi \in X_{N_{E/F}(E^*)}.$ 

$  A(w_2, w_1 w_2 \lambda) $ is no isomorphism if and only if $ \mid \; \mid^{\alpha + 1/2} \chi \rtimes \lambda' $
reduces if and only if $ \alpha = 1/2 $ and $ \chi = 1. $

\bigskip

In all other cases $ A(w_2, \lambda), A(w_1, w_2 \lambda) $ and $ A(w_2, w_1 w_2 \lambda) $ are holomorphic and
 isomorphisms and 
$ A(w', \lambda) $  is also an isomorphism. Hence the representations
$ \mid \; \mid^{\alpha} \chi(\det) \St_{\GL_2} \rtimes \lambda' $ and
$ \mid \; \mid^{\alpha} \chi(\det) 1_{\GL_2} \rtimes \lambda' $ are irreducible.

If one of the three cases in Theorem 1.5 holds, reducibility of 
$ \mid \; \mid^{\alpha} \chi(\det) \St_{\GL_2} \rtimes \lambda' $ and
$ \mid \; \mid^{\alpha} \chi(\det) 1_{\GL_2} \rtimes \lambda' $ has to be investigated. This is done in \ref{chistu3},
\ref{2111}, \ref{3/21/2chiomega}, \ref{1chiomega} and in \ref{1chi1F*}.
\end{Proof}

\bigskip

\begin{center}{4.3.2 \textbf{Representations} $ \chi \St_{\GL_2} \rtimes \lambda' $ \textbf{and}
$ \chi 1_{\GL_2} \rtimes \lambda' $}
\end{center}

Let $ 0 < \alpha_2 \leq \alpha_1, \; \alpha > 0. $ Let $ \chi_1, \chi_2, \chi $ and $ \chi' $ be unitary characters of $ E^*. $ Let
$ \lambda' $ be a unitary character of $ E^1. $ Let $ \chi \notin X_{1_{F^*}} $ (hence $ \chi \rtimes \lambda' $ is irreducible by
\cite{Ky}). Let $ \tau_1 $ be a tempered representation of $ \GL(2,E), $ let $ \tau_2 $ be a
tempered representation of $ U(3) $ and let $ \tau $ be a tempered representation of $ U(5). $

The representations $ \mid \; \mid^{\alpha_1} \chi_1 \times \mid \; \mid^{\alpha_2} \chi_2 \rtimes \lambda', \;
\mid \; \mid^{\alpha} \chi_1 \times \chi \rtimes \lambda', \; \mid \; \mid^{\alpha} \tau_1 \rtimes \lambda', \;
\mid \; \mid^{\alpha} \chi' \rtimes \tau_2 $ and $ \tau $ have a unique irreducible quotient, the Langlands quotient, denoted
by
$ \Lg(\mid \; \mid^{\alpha_1} \chi_1 ; \mid \; \mid^{\alpha_2} \chi_2 \; \lambda'), \;
\Lg(\mid \; \mid^{\alpha} \chi_1 ; \chi \rtimes \lambda'), \; \Lg(\mid \; \mid^{\alpha} \tau_1 ; \lambda'), \;
\Lg(\mid \; \mid^{\alpha} \chi' ; \tau_2) $ and $ \tau, $ respectively.

\begin{Proposition}

\label{chist}

 Let $ \chi $ be a unitary character of $ E^*, $ let $ \lambda' $ be a (unitary) character of $ E^1. $
The representations $ \chi \St_{\GL_2} \rtimes \lambda' $ and $ \chi 1_{\GL_2} \rtimes \lambda' $ are reducible if
and only if $ \chi \in X_{\omega_{E/F}}. $

\medskip

Let $ \chi =: \chi_{\omega_{E/F}} \in X_{\omega_{E/F}}, $ Let 
$ \pi_{1,\chi_{\omega_{E/F}}} $ be the unique irreducible square-integrable subquotient
of $ \mid \; \mid^{1/2} \chi_{\omega_{E/F}} \rtimes \lambda'  $ \cite{Ky}. Then

\smallskip

$ \chi_{\omega_{E/F}} 1_{\GL_2} \rtimes \lambda' =
 \Lg(\mid \; \mid^{1/2}  \chi_{\omega_{E/F}}; \mid \; \mid^{1/2} \chi_{\omega_{E/F}} ; \lambda') +
 \Lg(\mid \; \mid^{1/2} \chi_{\omega_{E/F}} ; \pi_{1,\chi_{\omega_{E/F}}}), $

$ \chi_{\omega_{E/F}} \St_{\GL_2} \rtimes \lambda' = \tau_1 + \tau_2, $

\smallskip

where $ \tau_1 $ and $ \tau_2 $ are tempered such that $ \tau_1 = \widehat{\Lg(\mid \; \mid^{1/2}  \chi_{\omega_{E/F}}; 
\mid \; \mid^{1/2} \chi_{\omega_{E/F}} ; \lambda')} $ and $ \tau_2 =
 \widehat{\Lg(\mid \; \mid^{1/2} \chi_{\omega_{E/F}} ; \pi_{1,\chi_{\omega_{E/F}}})}. $
All subquotients are unitary.

\end{Proposition}

\bigskip

\begin{Proof}
We consider Jaquet-restriction of $ \chi 1_{\GL_2} \rtimes \lambda' $ to the minimal parabolic subgroup:

\smallskip

$ s_{\min} (\chi(\det) 1_{\GL_2} \rtimes \lambda') = \mid \; \mid^{-1/2} \chi \otimes \mid \; \mid^{1/2} \chi \otimes \lambda'
 + \mid \; \mid^{-1/2} \chi^{-1}
(\overset{-}{\;}) \otimes \mid \; \mid^{1/2} \chi^{-1}(\overset{-}{\;}) \otimes \lambda' + \mid \; \mid^{-1/2} \chi \otimes
\mid \; \mid^{-1/2} \chi^{-1}(\overset{-}{\;}) \otimes \lambda' +  \mid \; \mid^{-1/2} \chi^{-1}
(\overset{-}{\;}) \otimes \mid \; \mid^{-1/2} \chi \otimes \lambda'. $

\medskip

Hence all subquotients of $ \chi 1_{\GL_2} \rtimes \lambda' $ are non-tempered.

\bigskip

$ \chi \St_{\GL_2} \rtimes \lambda' $ and $ \chi 1_{\GL_2} \rtimes \lambda' $ are subquotients of $ \mid \; \mid^{1/2} \chi
\times \mid \; \mid^{1/2} \chi \rtimes \lambda'. $

\medskip
 
For $ w =
\left(
\begin{smallmatrix}
 1&&&&\\
&&&&1&\\
&&1&&\\
&1&&&\\
&&&&1
\end{smallmatrix}
\right) $ we have $ w (\mid \; \mid^{1/2} \chi \otimes \mid \; \mid^{-1/2} \chi \otimes \lambda') =
\mid \; \mid^{1/2} \chi \otimes \mid \; \mid^{1/2} \chi^{-1}(\overset{-}{\;}) \otimes \lambda', $ and 
$ \mid \; \mid^{1/2} \chi \times \mid \; \mid^{-1/2} \chi \rtimes \lambda' $ and
 $ \mid \; \mid^{1/2} \chi \times \mid \; \mid^{1/2} \chi^{-1}(\overset{-}{\;}) \rtimes \lambda'$ have the same
irreducible constituents. Therefore we consider the reducibility of $
 \mid \; \mid^{1/2} \chi \times \mid \; \mid^{1/2} \chi^{-1}(\overset{-}{\;}) \rtimes \lambda'. $

\medskip

If $ \chi \notin X_{\omega_{E/F}}, $ then $ \Lg(\mid \; \mid^{1/2} \chi ; \mid \; \mid^{1/2}
 \chi^{-1}(\overset{-}{\;}) ; \lambda')$ is the only non-tempered Langlands quotient supported in $
\mid \; \mid^{1/2} \chi \otimes \mid \; \mid^{1/2}
 \chi^{-1}(\overset{-}{\;}) \otimes \lambda'. $ Hence  $ \chi 1_{\GL_2} \rtimes \lambda' =
 \Lg(\mid \; \mid^{1/2} \chi ; \mid \; \mid^{1/2} \chi^{-1}(\overset{-}{\;}) ; \lambda')$ is irreducible.
$ \chi \St_{\GL_2} \rtimes \lambda' $ is irreducible by Aubert duality, it is tempered.

\medskip

Let $ \chi =: \chi_{\omega_{E/F}} \in X_{\omega_{E/F}}. \; 
\Lg(\mid \; \mid^{1/2}  \chi_{\omega_{E/F}}; \mid \; \mid^{1/2} \chi_{\omega_{E/F}} ; \lambda') $ and $
 \Lg(\mid \; \mid^{1/2} \chi_{\omega_{E/F}} ; \pi_{1,\chi_{\omega_{E/F}}})$ are the only non-tempered Langlands quotients
supported in $ \mid \; \mid^{1/2} \chi_{\omega_{E/F}} \otimes \mid \; \mid^{1/2} \chi_{\omega_{E/F}} \otimes \lambda'. \:
 \chi_{\omega_{E/F}} \St_{\GL_2} \rtimes \lambda' $ is tempered and so are its subquotients. Hence 
$ \Lg(\mid \; \mid^{1/2}  \chi_{\omega_{E/F}}; \mid \; \mid^{1/2} \chi_{\omega_{E/F}} ; \lambda') $ and $
 \Lg(\mid \; \mid^{1/2} \chi_{\omega_{E/F}} ; \pi_{1,\chi_{\omega_{E/F}}}) $ are the subquotients of  
$ \chi_{\omega_{E/F}} 1_{\GL_2} \rtimes \lambda'. $ By Aubert duality  $\chi_{\omega_{E/F}} \St_{\GL_2} \rtimes \lambda' $
has the two irreducible subquotients $ \tau_1 := \widehat{\Lg(\mid \; \mid^{1/2}  \chi_{\omega_{E/F}}; \mid \; \mid^{1/2} 
\chi_{\omega_{E/F}} ; \lambda')} $ and
$ \tau_2 := \widehat{\Lg(\mid \; \mid^{1/2} \chi_{\omega_{E/F}} ; \pi_{1,\chi_{\omega_{E/F}}})}. $

\medskip

We consider restriction to the parabolic subgroup $ P_1: $

$ s_{P_1}(\chi_{\omega_{E/F}} \St_{\GL_2} \rtimes \lambda') = \chi_{\omega_{E/F}} \St_{\GL_2} \otimes 
\lambda' + \chi_{\omega_{E/F}} \St_{\GL_2} \otimes \lambda' + \mid \; \mid^{1/2} \chi_{\omega_{E/F}} \times
\mid \; \mid^{1/2} \chi_{\omega_{E/F}} \otimes \lambda'. $

\smallskip

$ \chi_{\omega_{E/F}} \St_{\GL_2} \rtimes \lambda' $ is unitary, hence $ \tau_1 \hookrightarrow 
\chi_{\omega_{E/F}} \St_{\GL_2} \rtimes \lambda' $ and $ \tau_2 \hookrightarrow \chi_{\omega_{E/F}} \St_{\GL_2} 
\rtimes \lambda'. $
By Frobenius reciprocity $ s_{P_1}(\tau_1) \twoheadrightarrow  \chi_{\omega_{E/F}} \St_{\GL_2} \otimes \lambda' $ and
$ s_{P_1}(\tau_2) \twoheadrightarrow  \chi_{\omega_{E/F}} \St_{\GL_2} \otimes \lambda' . \; 
\chi_{\omega_{E/F}} \St_{\GL_2} \otimes \lambda' $ is irreducible and has multiplicity 2 in
$ s_{P_1}(\chi_{\omega_{E/F}} \St_{\GL_2} \rtimes \lambda'). $ Hence $ \tau_1 $ and $ \tau_2 $ have multiplicities 1
 and
$ \chi_{\omega_{E/F}} \St_{\GL_2} \rtimes \lambda' $ is a representation of length 2. By Aubert duality
$ \Lg(\mid \; \mid^{1/2}  \chi_{\omega_{E/F}}; \mid \; \mid^{1/2} \chi_{\omega_{E/F}} ; \lambda') $ and $
 \Lg(\mid \; \mid^{1/2} \chi_{\omega_{E/F}} ; \pi_{1,\chi_{\omega_{E/F}}}) $ have multiplicities 1 and 
$ \chi_{\omega_{E/F}}
1_{\GL_2} \rtimes \lambda' $ is of length 2.

\medskip

$ \chi_{\omega_{E/F}}
\St_{\GL_2} \rtimes \lambda' $ and $ \chi_{\omega_{E/F}} 1_{\GL_2} \rtimes \lambda' $ are unitary, hence all
subquotients are unitary.
\end{Proof}

\bigskip

\begin{center}{4.3.3 \textbf{Representations} $ \mid \; \mid^{\alpha} \chi \rtimes \tau $ \textbf{and}
$ \chi \rtimes \tau, \; \alpha > 0, \;
 \tau $ \textbf{irreducible non-cuspidal of
$ U(3), $ not fully-induced}} \end{center}

We now look at representations induced from the maximal parabolic subgroup $ P_2, $ whose cuspidal support
is in $ M_0 $ and that are not fully induced.

\bigskip

Recall that $ P_2 = M_2 N_2, $ where $ M_2 \cong E^* \times U(3) $ is a maximal standard
Levi subgroup and $ N_2 $ the unipotent subgroup corresponding to $ P_2 $ and $ M_2. $

\bigskip

Let $ \chi $ be a unitary character of $ E^*. $

Let $ \beta \in \R_+. $ Recall from \ref{irru3} the irreducible subquotients of the induced representations to $ U(3) $ in the
cases that
$ \mid \; \mid^{\beta} \chi \rtimes \lambda' $ is reducible:
$ \lambda'(\det) \St_{U(3)}, \lambda'(\det) 1_{U(3)}, \pi_{1, \chi_{\omega_{E/F}}}, \pi_{2, \chi_{\omega_{E/F}}}, 
\sigma_{1, \chi_{1_{F^*}}}, \sigma_{2, \chi_{1_{F^*}}}. $ All irreducible subquotients are unitary.

\medskip

Let $ \alpha \in \R_+^*. $ We study the representations $ \mid \; \mid^{\alpha} \chi \rtimes \lambda'(\det) \St_{U(3)}, \;
 \mid \; \mid^{\alpha} \chi \rtimes 
\lambda'(\det)
1_{U(3)}, \; \mid \; \mid^{\alpha} \chi \rtimes \pi_{1, \chi_{\omega_{E/F}}}, \; \mid \; \mid^{\alpha} \chi \rtimes 
\pi_{2, \chi_{\omega_{E/F}}}, \mid \; \mid^{\alpha} \chi \rtimes
\sigma_{1, \chi_{1_{F^*}}} $ and $ \mid \; \mid^{\alpha} \chi \rtimes \sigma_{2, \chi_{1_{F^*}}}. $ Further we study
representations $  \chi \rtimes \lambda'(\det) \St_{U(3)}, \;
 \chi \rtimes 
\lambda'(\det)
1_{U(3)}, \;  \chi \rtimes \pi_{1, \chi_{\omega_{E/F}}}, \; \chi \rtimes 
\pi_{2, \chi_{\omega_{E/F}}}. \;  \chi \rtimes
\sigma_{1, \chi_{1_{F^*}}} $ and $  \chi \rtimes \sigma_{2, \chi_{1_{F^*}}}. $

\bigskip

\begin{center}{4.3.3.1 Representations $ \mid \; \mid^{\alpha} \chi \rtimes \lambda'(\det) \St_{U(3)} $ and
 $ \mid \; \mid^{\alpha} \chi \rtimes \lambda'(\det) 1_{U(3)}, \; \alpha > 0 $}\end{center}

\begin{Theorem} 

\label{alphachistu3}

Let $ \alpha \in R_+^* $ and $ \chi $ be a unitary character of $ E^*. $
The representations $ \mid \; \mid^{\alpha} \chi \rtimes \lambda'(\det) \St_{U(3)} $ and
$ \mid \; \mid^{\alpha} \chi 
\rtimes
\lambda'(\det) 1_{U(3)} $ are irreducible unless one of the following conditions holds:

\bigskip

1. $ \alpha = 2 $ and $ \chi = 1, $

2. $ \alpha =  1 $ and $ \chi = 1, $

3. $ \alpha = 1/2 $ and  $ \chi \in X_{\omega_{E/F}}. $

\end{Theorem}

\begin{Proof} In $ R(U) $ we have  $\mid \; \mid^{\alpha} \chi \times \mid \; \mid^1 1 \rtimes \lambda' =
\mid \; \mid^{\alpha} \chi \rtimes \lambda'
\det \St_{U(3)} + \mid \; \mid^{\alpha} \chi \rtimes \lambda' \det 1_{U(3)}. $
We give the proof for $ \mid \; \mid^{\alpha} \chi \rtimes \lambda'(\det) \St_{U(3)}, $ its Aubert dual
$ \mid \; \mid^{\alpha} \chi \rtimes \lambda'(\det) 1_{U(3)} $ has the same points of reducibility.

Let $ \lambda := \mid \; \mid^{\alpha} \chi \otimes \mid \; \mid 1 \otimes \lambda'. $
Let $ A(w'', \lambda): \mid \; \mid^{\alpha} \chi \rtimes \lambda'(\det) \St_{U(3)} \rightarrow
\mid \; \mid^{-\alpha} \chi^{-1}(\overset{-}{\;}) \rtimes \lambda'(\det) \St_{U(3)} $ be the long intertwining operator for the
representation $ \mid \; \mid^{\alpha} \chi \rtimes \lambda'(\det) \St_{U(3)} $ corresponding to $ w'', $ the longest element of
$ W $ respecting the Levi subgroup $ M_2 \cong E^* \times U(3). $
We analyse the decomposition of $ A(w'', \lambda) $ into short intertwining operators.
As in the previous case this will show when  $ \mid \; \mid^{\alpha} \chi \rtimes \lambda' \det
\St_{U(3)} $ and $ \mid \; \mid^{\alpha} \chi \rtimes
\lambda'(\det) 1_{U(3)} $ are irreducible.

\bigskip

Let
$ w'' :=
\left(
\begin{smallmatrix}
&&&&1\\
&1&&&\\
&&1&&\\
&&&1&\\
1&&&&
\end{smallmatrix}
\right) =
\left(
\begin{smallmatrix}
0&1&&&\\
1&0&&&\\
&&1&&\\
&&&0&1\\
&&&1&0
\end{smallmatrix}
\right)
\left(
\begin{smallmatrix}
1&&&&\\
&&&1&\\
&&1&&\\
&1&&&\\
&&&&1
\end{smallmatrix}
\right)
\left(
\begin{smallmatrix}
0&1&&&\\
1&0&&&\\
&&1&&\\
&&&0&1\\
&&&1&0
\end{smallmatrix}
\right) = w_1w_2w_1
$

\bigskip

The following diagram gives the decomposition of $ A(w'', \lambda) $ into short intertwining operators.
$ i_1 $ and $ i_2 $ are inclusions and depend holomorphically on $ \alpha. $

\bigskip

$
\begin{smallmatrix}
 \mid \; \mid^{\alpha} \chi \rtimes \lambda'(\det) \St_{U(3)} & \overset{i_1}{\hookrightarrow} &\mid \;
 \mid^{\alpha} \chi \times \mid \; \mid 1 \rtimes \lambda' \\
&& A(w_1, \lambda) \downarrow w_1 \\
&& \mid \; \mid 1 \times \mid \; \mid^{\alpha} \chi \rtimes \lambda' \\
A(w'', \lambda) \downarrow w'' &&  A(w_2, w_1 \lambda) \downarrow w_2 \\
&& \mid \; \mid 1 \times \mid \; \mid^{ - \alpha} \chi^{-1}(\overset{-}{\;}) \rtimes \lambda' \\
&&  A(w_1, w_2 w_1 \lambda) \downarrow w_1 \\
\mid \; \mid^{-\alpha} \chi^{-1}(\overset{-}{\;}) \rtimes \lambda'(\det) \St_{U(3)}&\overset{i_2}{\hookrightarrow} 
& \mid \; \mid^{-\alpha} \chi^{-1}(\overset{-}{\;}) \times \mid \; \mid 1 \rtimes \lambda' \\
\end{smallmatrix}
$

 \bigskip
$ A(w_1, \lambda) $ is no isomorphism if and only if $ \mid \; \mid^{\alpha} \chi \times \mid \; \mid 1 $
reduces if and only if $ \alpha = 2 $ and $ \chi = 1. $

\medskip

$ A(w_2, w_1 \lambda) $ is no isomorphism if and only if $ \mid \; \mid^{\alpha} \chi \rtimes \lambda' $ reduces
if and only of $ \alpha = 1 $ and $ \chi = 1 $ or $ \alpha = 1/2 $ and $ \chi \in X_{\omega_{E/F}}. $

\medskip

$ A(w_1, w_2w_1\lambda) $ is an isomorphism as $\mid \; \mid 1 \times \mid \; \mid^{ - \alpha} \chi^{-1}(\overset{-}{\;}) $
is irreducible for all $ \alpha \in \R_+^*. $

\bigskip
 
In those cases where $ A(w_1, \lambda), \; A(w_2, w_1 \lambda)  $ and $ A(w_2, w_1 w_2 \lambda) $ are isomorphisms, 
$ A(w'', \lambda) $ is an 
isomorphism, and the representations
$ \mid \; \mid^{\alpha} \chi \rtimes \lambda'(\det) \St_{U(3)} $ and $ \mid \; \mid^{\alpha} \chi \rtimes
\lambda'(\det) 1_{U(3)} $ are irreducible.

On the other hand, if $ A(w_1, \lambda), \; A(w_2, w_1 \lambda)  $ or $ A(w_2, w_1 w_2 \lambda) $ is no isomorphism, then
reducibility of $ \mid \; \mid^{\alpha} \chi \rtimes \lambda'(\det) \St_{U(3)} $ and
 $ \mid \; \mid^{\alpha} \chi \rtimes
\lambda'(\det) 1_{U(3)} $ is left to be investigated. It is done in \ref{2111} and in \ref{11/21chiomega}.
\end{Proof}

\bigskip

\begin{center}{4.3.3.2 Representations $ \chi \rtimes \lambda'(\det) \St_{U(3)} $ and 
$ \chi \rtimes \lambda'(\det) 1_{U(3)} $}\end{center}

Let $ \chi_{1_{F^*}} \in X_{1_{F^*}}. $ Recall that $ \chi_{1_{F^*}} \rtimes \lambda' = \sigma_{1, \chi_{1_{F^*}}} \oplus
 \sigma_{2, \chi_{1_{F^*}}}, $ where
$ \sigma_{1, \chi_{1_{F^*}}} $ and $ \sigma_{2, \chi_{1_{F^*}}} $ are tempered \cite{Ky}.

\medskip

\begin{Proposition}

\label{chistu3}
 Let $ \chi $ be a unitary character of $ E^*, $ let $ \lambda' $ be a (unitary) character of $ E^1. $
The representations $ \chi \rtimes \lambda'(\det) \St_{U(3)} $ and 
$ \chi \rtimes \lambda'(\det) 1_{U(3)} $ are reducible if and only if $ \chi = 1 $ or $ \chi \in X_{1_{F^*}}. $

\medskip

Let $ \chi = 1. $ 

\smallskip
$ 1 \rtimes \lambda'(\det) 1_{U(3)} = \Lg(\mid \; \mid 1 ; 1 \rtimes \lambda') + \Lg( \mid \; \mid^{1/2} \St_{\GL_2}
; \lambda') $ and

$ 1 \rtimes \lambda'(\det) \St_{U(3)} = \tau_3 + \tau_4, $

where $ \tau_3 $ and $ \tau_4 $ are tempered such that $ \tau_3 = \widehat{\Lg(\mid \; \mid 1 ; 1 \rtimes \lambda')} $
and $ \tau_4 = \widehat{\Lg( \mid \; \mid^{1/2} \St_{\GL_2}; \lambda')}. $

\medskip

Let $ \chi =: \chi_{1_{F^*}} \in X_{1_{F^*}}. $

\smallskip
$ \chi_{1_{F^*}} \rtimes \lambda'(\det) 1_{U(3)} = \Lg(\mid \; \mid 1 ; \sigma_{1, \chi_{1_{F^*}}}) + \Lg( \mid \; \mid 1 ;
\sigma_{2, \chi_{1_{F^*}}}) $ and

$ \chi_{1_{F^*}} \rtimes \lambda'(\det) \St_{U(3)} = \tau_5 + \tau_6. $

$ \tau_5 $ and $ \tau_6 $ are tempered, such that $ \tau_5 = \widehat{\Lg(\mid \; \mid 1 ; 
\sigma_{1, \chi_{1_{F^*}}})}, \;
 \tau_6 =  \widehat{\Lg(\mid \; \mid 1; \sigma_{2 ; \chi_{1_{F^*}}})}. $
\end{Proposition}

\begin{Proof}
We consider the Jaquet restriction of $ \chi \rtimes \lambda'(\det) 1_{U(3)} $ to the minimal parabolic subgroup:

\smallskip

$ s_{\min}(\chi \rtimes \lambda'(\det) 1_{U(3)}) =  \chi \otimes \mid \; \mid^{-1} 1 \otimes \lambda' +
\mid \; \mid^{-1} 1
\otimes \chi \otimes \lambda' + \chi^{-1}(\overset{-}{\;}) \otimes \mid \; \mid^{-1} 1 \otimes \lambda' +
\mid \; \mid^{-1} 1 \otimes \chi^{-1}(\overset{-}{\;}) \otimes \lambda'. $

Hence all subquotionts of $ \chi \rtimes \lambda'(\det) 1_{U(3)} $ are non-tempered.

\smallskip

$ \chi \rtimes \lambda'(\det) \St_{U(3)} $ is tempered, hence all its subquotients are tempered.

\medskip

If $ \chi \neq 1 $ and $ \chi \notin X_{1_{F^*}}, $ then

$ \Lg(\mid \; \mid 1 ; \chi \rtimes \lambda') $ is the only non-tempered Langlands quotient supported in $ \mid \; \mid 1 
\otimes \chi \otimes \lambda', $ and $ \chi \rtimes \lambda'(\det) 1_{U(3)} = \Lg(\mid \; \mid 1 ; \chi \rtimes \lambda') $
is irreducible. 

$ \chi \rtimes \lambda'(\det) \St_{U(3)} $ is irreducible by Aubert duality.

\medskip

Let $ \chi = 1. \; \Lg(\mid \; \mid 1 ; 1 \rtimes \lambda') $ and $  \Lg( \mid \; \mid^{1/2} \St_{\GL_2}
; \lambda') $ are the only non-tempered Langlands-quotients supported in $ \mid \; \mid 1 \otimes 1 \otimes \lambda'. $
 $ 1 \rtimes \lambda'(\det) \St_{U(3)} $ is tempered and so are its subquotients. Hence 
$ \Lg(\mid \; \mid 1 ; 1 \rtimes \lambda') $ and $  \Lg( \mid \; \mid^{1/2} \St_{\GL_2}
; \lambda') $ are the irreducible subquotients of 
$ 1 \rtimes \lambda'(\det) 1_{U(3)}.$ By Aubert duality $ 1 \rtimes \lambda'(\det) \St_{U(3)} $ has two irreducible 
subquotients, the tempered representations
$ \tau_3 :=  \widehat{\Lg(\mid \; \mid 1 ; 1 \rtimes \lambda')} $ and
$ \tau_4 := \widehat{\Lg( \mid \; \mid^{1/2} \St_{\GL_2}; \lambda')}. $

\medskip

We consider restriction to the parabolic subgroup $ P_2. $

$ s_{P_2}(1 \rtimes \lambda'(\det) \St_{U(3)}) = 1 \otimes \lambda'(\det) \St_{U(3)} + 1 \otimes \lambda'(\det) \St_{U(3)}
+  \mid \; \mid 1 \otimes 1 \rtimes \lambda'. $ 

\smallskip

$ 1 \otimes \lambda'(\det) \St_{U(3)} $ is irreducible and appears with multiplicity 2 in 
$ s_{P_2}(1 \rtimes \lambda'(\det) \St_{U(3)}). $

$  1 \rtimes \lambda'(\det) \St_{U(3)} $ is unitary, hence $ \tau_3 \hookrightarrow 
 1 \rtimes \lambda'(\det) \St_{U(3)} $ and $ \tau_4 \hookrightarrow
 1 \rtimes \lambda'(\det) \St_{U(3)}. $ By Frobenius reciprocity 
$ s_{P_2}(\tau_3) \twoheadrightarrow  1 \otimes \lambda'(\det) \St_{U(3)}, $ and $ s_{P_2}(\tau_4) \twoheadrightarrow
1 \otimes \lambda'(\det) \St_{U(3)}. $

Hence multiplicities of $ \tau_3 $ and $ \tau_4 $ equal 1 and $ 1 \rtimes \lambda'(\det) \St_{U(3)} $
 is a representation of length 2. By Aubert duality multiplicities of the irreducible subquotients of
$ 1 \rtimes \lambda'(\det) 1_{U(3)} $ equal 1 and it is a representation of
length 2.

\medskip

If $ \chi =: \chi_{1_{F^*}} \in X_{1_{F^*}}, $ then $ \Lg(\mid \; \mid 1 ; \sigma_{1, \chi_{1_{F^*}}}) $ and 
$ \Lg( \mid \; \mid 1 ;
\sigma_{2, \chi_{1_{F^*}}}) $ are the only non-tempered Langlands-quotients supported in $ \mid \; \mid 1 \otimes \chi_{1_{F^*}}
\otimes \lambda'. \; \chi_{1_{F^*}} \rtimes \lambda'(\det) \St_{U(3)} $  is tempered and so are its subquotients. Hence
$ \Lg(\mid \; \mid 1 ; \sigma_{1, \chi_{1_{F^*}}}) $ and 
$ \Lg( \mid \; \mid 1 ;
\sigma_{2, \chi_{1_{F^*}}}) $ are the irreducible subquotients of $ \chi_{1_ {F^*}}
\rtimes \lambda'(\det) 1_{U(3)}. $ By Aubert duality $ \chi_{1_{F^*}} \rtimes \lambda'(\det) \St_{U(3)} $ has two
irreducible
subquotients, the tempered representations $ \tau_5 := \widehat{\Lg(\mid \; \mid 1 ; 
\sigma_{1, \chi_{1_{F^*}}})} $ and $
 \tau_6 :=  \widehat{\Lg(\mid \; \mid 1; \sigma_{2 ; \chi_{1_{F^*}}})}. $

\medskip

We consider the restriction to the parabolic subgroup $ P_2: $

$ s_{P_2}(\chi_{1_{F^*}} \rtimes \lambda'(\det) \St_{U(3)}) = \chi_{1_{F^*}} \otimes \lambda'(\det) \St_{U(3)} +
 \chi_{1_{F^*}} \otimes \lambda'(\det) \St_{U(3)}
+  \mid \; \mid 1 \otimes \chi_{1_{F^*}} \rtimes \lambda'. $ 

\smallskip

$ \chi_{1_{F^*}} \otimes \lambda'(\det) \St_{U(3)} $ is irreducible and appears with multiplicity 2 in 
$ s_{P_2}( \chi_{1_{F^*}} \rtimes \lambda'(\det) \St_{U(3)}). $
$  \chi_{1_{F^*}} \rtimes \lambda'(\det) \St_{U(3)} $ is unitary, hence $ \tau_5 \hookrightarrow 
 \chi_{1_{F^*}} \rtimes \lambda'(\det) \St_{U(3)} $ and $ \tau_6 \hookrightarrow
 \chi_{1_{F^*}} \rtimes \lambda'(\det) \St_{U(3)}. $ By Frobenius reciprocity
$ s_{P_2}(\tau_5) \twoheadrightarrow  \chi_{1_{F^*}} \otimes \lambda'(\det) \St_{U(3)}, $ and $ s_{P_2}(\tau_6)
\twoheadrightarrow
\chi_{1_{F^*}} \otimes \lambda'(\det) \St_{U(3)}. $

Hence multiplicities of $ \tau_5 $ and $ \tau_6 $ equal 1 and $ \chi_{1_{F^*}} \rtimes \lambda'(\det) \St_{U(3)} $
 is a representation of length 2. By Aubert duality multiplicities of the irreducible subquotients of
$ \chi_{1_{F^*}} \rtimes \lambda'(\det) 1_{U(3)} $ equal 1 and it is a representation of
length 2.
\end{Proof}

\bigskip

\begin{center}{4.3.3.3 Representations $ \mid \; \mid^{\alpha} \chi \rtimes \pi_{1, \chi_{\omega_{E/F}}} $ and
 $ \mid \; \mid^{\alpha} \chi \rtimes \pi_{2, \chi_{\omega_{E/F}}}, \; \alpha > 0 $}\end{center}

Let $ \alpha \in \R_+^* $ and let $ \chi $ be a unitary character of $ E^*. $
Let $ \chi_{\omega_{E/F}} \in X_{\omega_{E/F}}, $ i.e. $ \chi_{\omega_{E/F}} $ is a (unitary) character of $ E^* $ whose
restriction to $ F^* $ equals the character $ \omega_{E/F}. $

Let $ \pi_{1, \chi_{\omega_{E/F}}} $ be the unique square-integrable subquotient and let
$ \pi_{2, \chi_{\omega_{E/F}}} $ be the unique 
irreducible non-tempered subquotient of $ \mid \; \mid^{1/2} \chi_{\omega_{E/F}} \rtimes \lambda' $ \cite{Ky}.

\bigskip

\begin{Theorem}

\label{alphachipi1}

 Let $ \alpha \in \R_+^* $ and $ \chi $ be a unitary character of $ E^*. $
 The representations $ \mid \; \mid^{\alpha} \chi \rtimes \pi_{1, \chi_{\omega_{E/F}}} $ and $ \mid \; \mid^{\alpha} \chi \rtimes
\pi_{2, \chi_{\omega_{E/F}}} $ are irreducible unless

\smallskip
1. $ \alpha = 1/2 $ or $ \alpha = 3/2 $ and $ \chi = \chi_{\omega_{E/F}}, $

2. $ \alpha =  1 $ and $ \chi = 1, $

3. $ \alpha = 1/2 $ and $ \chi \in X_{\omega_{E/F}}. $
\end{Theorem}

\begin{Proof} Let $ \alpha \in \R_+^*. $
Let $ \lambda := \mid \; \mid^{\alpha} \chi \otimes \mid \; \mid^{1/2} \chi_{\omega_{E/F}} \otimes \lambda'. $

Let $ A(w'', \lambda): \; \mid \; \mid^{\alpha} \chi \rtimes \pi_{1, \chi_{\omega_{E/F}}} \rightarrow
\mid \; \mid^{-\alpha} \chi^{-1}(\overset{-}{\;}) \rtimes \pi_{1, \chi_{\omega_{E/F}}} $ be the long intertwining operator for the
representation $ \mid \; \mid^{\alpha} \chi \rtimes \pi_{1, \chi_{\omega_{E/F}}} $ corresponding to $ w'', $ the longest element of $ W $
respecting $ M_2 \cong E^* \times U(3). $
We analyse the decomposition of $ A(w'', \lambda) $ into short intertwining operators.
As in the previous case this shows when  $ \mid \; \mid^{\alpha} \chi \rtimes 
\pi_{1, \chi_{\omega_{E/F}}} $ and $ \mid \; \mid^{\alpha} \chi \rtimes
\pi_{2, \chi_{\omega_{E/F}}} $ are irreducible.

\bigskip

The following diagram gives the decomposition of $ A(w'', \lambda) $ into short intertwining operators.
$ w_1 $ and $ w_2 $ are like before:

$
\begin{smallmatrix}
 \mid \; \mid^{\alpha} \chi \rtimes \pi_{1, \chi_{\omega_{E/F}}} &\overset{i_1}{\hookrightarrow} & \mid \; \mid^{\alpha} \chi \times 
\mid \; \mid^{1/2} \chi_{\omega_{E/F}} \rtimes \lambda' \\
&& A(w_1, \lambda) \downarrow w_1 \\
&& \mid \; \mid^{1/2} \chi_{\omega_{E/F}} \times \mid \; \mid^{\alpha} \chi \rtimes \lambda' \\
A(w'', \lambda) \downarrow w'' && A(w_2, w_1 \lambda) \downarrow w_2 \\
&& \mid \; \mid^{ 1/2} \chi_{\omega_{E/F}} \times \mid \; \mid^{-\alpha} \chi^{-1} (\overset{-}{\;}) \rtimes \lambda'\\
&& A(w_1 , w_2 w_1 \lambda) \downarrow w_1\\
\mid \; \mid^{- \alpha} \chi^{-1}(\overset{-}{\;}) \rtimes \pi_{1, \chi_{\omega_{E/F}}} & \overset{i_2}{\hookrightarrow} &
 \mid \; \mid^{- \alpha} \chi^{-1}(\overset{-}{\;}) \times \mid \; \mid^{1/2} \chi_{\omega_{E/F}} \rtimes \lambda'
\end{smallmatrix}
$

\bigskip

$ A (w_1, \lambda) $ is not an isomorphism if and only if
$ \alpha = 3/2 $ and $ \chi = \chi_{\omega_{E/F}}.$

$ A(w_2, w_1 \lambda) $ is not an isomorphism if and only if $ \alpha = 1 $ and $ \chi = 1 $ or $ \alpha = 1/2 $
and $ \chi \in X_{\omega_{E/F}}. $

\medskip

$ A(w_1, w_2 w_1 \lambda) $ is not an isomorphism if and only if $ \alpha = 1/2 $ and
 $ \chi = \chi_{\omega_{E/F}}. $

\bigskip

In those cases where $ A(w_1, \lambda), A(w_2, w_1 \lambda), A(w_1, w_2 w_1 \lambda) $ and hence $ A(w'', \lambda) $
 are isomorphisms,
$ \mid \; \mid^{\alpha} \chi \rtimes \pi_{1, \chi_{\omega_{E/F}}} $ and $ \mid \; \mid^{\alpha} \chi \rtimes
 \pi_{2, \chi_{\omega_{E/F}}} $
are irreducible.

On the other hand, if $ A(w_1, \lambda), A(w_2, w_1 \lambda) $ or $ A(w_1, w_2 w_1 \lambda) $ is no isomorphism, 
reducibility of
 $ \mid \; \mid^{\alpha} \chi \rtimes \pi_{1, \chi_{\omega_{E/F}}} $ and
$ \mid \; \mid^{\alpha} \chi \rtimes \pi_{2, \chi_{\omega_{E/F}}} $ needs to be investigated. This is done in 
\ref{chist} and \ref{1/21/2chiomega} (for $ \mid \; \mid^{1/2} \chi_{\omega_{E/F}} \rtimes 
\pi_{1, \chi_{\omega_{E/F}}} $ and $ \mid \; \mid^{1/2} \chi_{\omega_{E/F}} \rtimes 
\pi_{2, \chi_{\omega_{E/F}}}), $ in \ref{3/21/2chiomega}, \ref{11/21chiomega} and in \ref{1/21/2chiomega12}. 
\end{Proof}

\bigskip

\begin{center}{4.3.3.4 Representations $ \chi \rtimes \pi_{1, \chi_{\omega_{E/F}}} $ and
$ \chi \rtimes \pi_{2, \chi_{\omega_{E/F}}} $}\end{center}

Let $ \chi_{\omega_{E/F}} \in X_{\omega_{E/F}}. $ Let
$ \pi_{1,\chi_{\omega_{E/F}}} $ be the unique square-integrable subquotient and $ \pi_{2, \chi_{\omega_{E/F}}} $ the
unique non-tempered irreducible subquotient of $ \mid \; \mid^{1/2} \chi_{\omega_{E/F}} \rtimes \lambda'. $
Let $ \chi_{1_{F^*}} \in X_{1_{F^*}}. $ Recall that $ \chi_{1_{F^*}} \rtimes \lambda' = \sigma_{1, \chi_{1_{F^*}}} \oplus
\sigma_{2, \chi_{1_{F^*}}}, $ where $ \sigma_{1, \chi_{1_{F^*}}} $ and $ \sigma_{2, \chi_{1_{F^*}}} $ are tempered
\cite{Ky}.

\medskip

\begin{Proposition}
\label{chipi1}

Let $ \chi $ be a unitary character of $ E^* $ and let $ \lambda' $ be a (unitary) character of $ E^1. $ 
The representations $ \chi \rtimes \pi_{1, \chi_{\omega_{E/F}}} $ and
$ \chi \rtimes \pi_{2, \chi_{\omega_{E/F}}} $ are reducible if and only if $ \chi \in X_{1_{F^*}}. $

\medskip

Let $ \chi =: \chi_{1_{F^*}} \in X_{1_{F^*}}. $ Then

\smallskip

$ \chi_{1_{F^*}} \rtimes \pi_{2, \chi_{\omega_{E/F}}} = \Lg(\mid \; \mid^{1/2} \chi_{\omega_{E/F}} ; \sigma_{1, \chi_{1_{F^*}}}) 
+ \Lg( \mid \; \mid^{1/2} \chi_{\omega_{E/F}} ;
\sigma_{2, \chi_{1_{F^*}}}), $

$ \chi_{1_{F^*}} \rtimes \pi_{1, \chi_{\omega_{E/F}}} = \tau_7 + \tau_8, $

where $ \tau_7 $ and $ \tau_8 $ are tempered representations such that 
$ \tau_7 = \widehat{\Lg(\mid \; \mid^{1/2} \chi_{\omega_{E/F}} ; 
\sigma_{1, \chi_{1_{F^*}}})} $ and $ \tau_8 = \widehat{\Lg( \mid \; \mid^{1/2} \chi_{\omega_{E/F}} ;
 \sigma_{2, \chi_{1_{F^*}}})}. $ 
\end{Proposition}

\begin{Proof}
We consider Jaquet restriction to the minimal parabolic subgroup:

\smallskip

$ s_{\min} (\chi \rtimes \pi_{2, \chi_{\omega_{E/F}}}) =  \chi \otimes \mid \; \mid^{-1/2} \chi_{\omega_{E/F}}
 \otimes \lambda' +
 \mid \; \mid^{-1/2} \chi_{\omega_{E/F}} \otimes
\chi \otimes \lambda' + \chi^{-1}(\overset{-}{\;}) \otimes \mid \; \mid^{-1/2} \chi_{\omega_{E/F}} \otimes
 \lambda' +  \mid \; \mid^{-1/2} \chi_{\omega_{E/F}} \otimes \chi^{-1}(\overset{-}{\;}) \otimes \lambda'. $

\smallskip

Hence all irreducible subquotients of $ \chi \rtimes \pi_{2,\chi_{\omega_{E/F}}} $ are non-tempered.

\smallskip

If $ \chi \notin X_{1_{F^*}}, $ then $ \Lg( \mid \; \mid^{1/2} \chi_{\omega_{E/F}} ; \chi \rtimes \lambda') $ is the only
non-tempered Langlands quotient supported in $ \mid \; \mid^{1/2} \chi_{\omega_{E/F}} \otimes \chi \otimes \lambda', $ and
$ \chi \rtimes \pi_{2,\chi_{\omega_{E/F}}} = \Lg( \mid \; \mid^{1/2} \chi_{\omega_{E/F}} ; \chi \rtimes \lambda') $
is irreducible. By Aubert duality 
$ \chi \rtimes \pi_{1,\chi_{\omega_{E/F}}} $ is irreducible, it is tempered.

\medskip

Let $ \chi =: \chi_{1_{F^*}} \in X_{1_{F^*}}. \; \Lg(\mid \; \mid^{1/2} \chi_{\omega_{E/F}} ; \sigma_{1, \chi_{1_{F^*}}}) $ and 
$ \Lg( \mid \; \mid^{1/2} \chi_{\omega_{E/F}} ;
\sigma_{2, \chi_{1_{F^*}}}) $ are the only non-tempered Langlands quotients supported in 
$ \mid \; \mid^{1/2} \chi_{\omega_{E/F}} \otimes \chi_{1_{F^*}} \otimes \lambda'. \; \chi_{1_{F^*}} \rtimes 
\pi_{1,\chi_{\omega_{E/F}}} $ is tempered, hence its subquotients are tempered. Hence 
$ \Lg(\mid \; \mid^{1/2} \chi_{\omega_{E/F}} ; \sigma_{1, \chi_{1_{F^*}}}) $ and 
$ \Lg( \mid \; \mid^{1/2} \chi_{\omega_{E/F}} ;
\sigma_{2, \chi_{1_{F^*}}}) $ are the subquotients of
$ \chi_{1_{F^*}} \rtimes \pi_{2,\chi_{\omega_{E/F}}}. $ By Aubert duality 
$ \chi_{1_{F^*}} \rtimes \pi_{1,\chi_{\omega_{E/F}}} $ has two subquotients, $ \tau_7 := \widehat{\Lg(\mid \; \mid^{1/2} 
\chi_{\omega_{E/F}} ; 
\sigma_{1, \chi_{1_{F^*}}})}$ and $ \tau_8 := \widehat{\Lg( \mid \; \mid^{1/2} \chi_{\omega_{E/F}}; 
\sigma_{2, \chi_{1_{F^*}}})}. $ They are tempered.

\medskip

We consider restriction to the parabolic subgroup $ P_2: $

$ s_{P_2}(\chi_{1_{F^*}} \rtimes \pi_{1,\chi_{\omega_{E/F}}}) = \chi_{1_{F^*}} \otimes \pi_{1,\chi_{\omega_{E/F}}} +
\chi_{1_{F^*}} \otimes \pi_{1,\chi_{\omega_{E/F}}} + \mid \; \mid^{1/2} \chi_{\omega_{E/F}} \otimes \chi_{1_{F^*}}
\rtimes \lambda' = 2 \; \chi_{1_{F^*}} \otimes \pi_{1,\chi_{\omega_{E/F}}} + \mid \; \mid^{1/2} \chi_{\omega_{E/F}}
 \otimes  \sigma_{1,\chi_{1_{F^*}}} + \mid \; \mid^{1/2} \chi_{\omega_{E/F}} \otimes \sigma_{2, \chi_{1_{F^*}}}. $

\smallskip

$ \chi_{1_{F^*}} \rtimes \pi_{1,\chi_{\omega_{E/F}}} $ is unitary, hence $ \tau_7 \hookrightarrow
\chi_{1_{F^*}} \rtimes \pi_{1,\chi_{\omega_{E/F}}} $ and $ \tau_8 \hookrightarrow
\chi_{1_{F^*}} \rtimes \pi_{1,\chi_{\omega_{E/F}}}. $ By Frobenius reciprocity $ s_{P_2}(\tau_7) \twoheadrightarrow 
\chi_{1_{F^*}} \otimes \pi_{1,\chi_{\omega_{E/F}}} $ and $ s_{P_2}(\tau_8) \twoheadrightarrow 
\chi_{1_{F^*}} \otimes \pi_{1,\chi_{\omega_{E/F}}}. \; \; \chi_{1_{F^*}} \otimes \pi_{1,\chi_{\omega_{E/F}}} $ is irreducible
and of multiplicity 2 in $ s_{P_2}(\chi_{1_{F^*}} \rtimes \pi_{1, \chi_{\omega_{E/F}}}), $ hence $ \tau_7 $ and $ \tau_8 
$
are of multiplicity 1 and $ \chi_{1_{F^*}} \times \pi_{1,\chi_{\omega_{E/F}}} $ is a representation of length 2.
By Aubert duality $ \Lg(\mid \; \mid^{1/2} \chi_{\omega_{E/F}} ; \sigma_{1, \chi_{1_{F^*}}}) $ and 
$ \Lg( \mid \; \mid^{1/2} \chi_{\omega_{E/F}} ;
\sigma_{2, \chi_{1_{F^*}}}) $ are of multiplicity 1, and $ \chi_{1_{F^*}} \rtimes \pi_{2,\chi_{\omega_{E/F}}} $ is a
representation of length 2.
\end{Proof}

\bigskip

\begin{center}{4.3.3.5 Representations $ \mid \; \mid^{\alpha} {\chi} \rtimes \sigma_{1, \chi_{1_{F^*}}} $ and $
 \mid \; \mid^{\alpha} \chi \rtimes 
\sigma_{2, \chi_{1_{F^*}}}, \; \alpha > 0 $}\end{center}

Let $ \alpha \in \R_+^* $ and let $ \chi $ be a unitary character of $ E^*. $ Let $ \chi_{1_{F^*}} \in X_{1_{F^*}}, $
 i.e.
$ \chi_{1_{F^*}} $ is a non-trivial unitary character of $ E^* $ whose restriction to $ F^* $ is trivial.
Recall that $ \chi_{1_{F^*}} \rtimes \lambda' = \sigma_{1, \chi_{1_{F^*}}} \oplus \sigma_{2, \chi_{1_{F^*}}}, $ where
$ \sigma_{1, \chi_{1_{F^*}}} $ and $ \sigma_{2, \chi_{1_{F^*}}} $ are tempered \cite{Ky}. 

\smallskip

\begin{Theorem}

\label{alphachisigma12}

Let $ \chi $ be a unitary character of $ E^*. $ Let $ \alpha \in \R_+^*. $
The representations $ \mid \; \mid^{\alpha} \chi \rtimes \sigma_{1, \chi_{1_{F^*}}} $ and 
$  \mid \; \mid^{\alpha} \chi 
\rtimes \sigma_{2, \chi_{1_{F^*}}} $
are
irreducible unless one of the following cases holds:

\bigskip

1. $ \alpha =  1 $ and $ \chi = \chi_{1_{F^*}}, $

2. $ \alpha =  1 $ and $ \chi = 1, $

3. $ \alpha = 1/2 $ and $ \chi \in X_{\omega_{E/F}}. $
\end{Theorem}

\begin{Proof}
Let $ \alpha \in \R_+^*. $ In $ R(U) $ we have $  \mid \; \mid^{\alpha} \chi \times \chi_{1_{F^*}} \rtimes \lambda' =
 \mid \; \mid^{\alpha} \chi \rtimes \sigma_{1, \chi_{1_{F^*}}} + \mid \; \mid^{\alpha} \chi \rtimes 
\sigma_{2, \chi_{1_F^*}}. $
We give the proof for $ \mid \; \mid^{\alpha} \chi \rtimes \sigma_{1, \chi_{1_{F^*}}}, $ the proof for 
$  \mid \; \mid^{\alpha} \chi 
\rtimes \sigma_{2, \chi_{1_{F^*}}} $ is analogous.
Let $ \lambda := \mid \; \mid^{\alpha} \chi \otimes \chi_{1_{F^*}} \otimes \lambda'. $ 
Let $ A(w'', \lambda): \mid \; \mid^{\alpha} \chi \times \chi_{1_{F^*}} \rtimes \lambda' \rightarrow \mid \; \mid^{-\alpha}
\chi^{-1}(\overset{-}{\;}) \times \chi_{1_{F^*}} \rtimes \lambda' $ be the intertwining operator corresponding to the longest
element $ w'' \in W $ respecting $ M_2 \cong \E^* \times U(3). $ The decomposition into short intertwining operators shows
when $ \mid \; \mid^{\alpha} \chi \rtimes \sigma_{1,\chi_{1_{F^*}}} $ is irreducible.

Let $ w'' = 
\left(
\begin{smallmatrix}
 &&&&1\\
&1&&&&\\
&&1&&\\
&&&1&\\
1&&&&
\end{smallmatrix}
\right) $ as before.

\bigskip

$
\begin{smallmatrix}
 \mid \; \mid^{\alpha} \chi \rtimes \sigma_{1, \chi_{1_{F^*}}} & \overset{i_1}{\hookrightarrow} &  \mid \; \mid^{\alpha} \chi \times 
\chi_{1_{F^*}} \rtimes \lambda' \\
&& A(w_1, \lambda) \downarrow w_1\\
&& \chi_{1 \mid F^*} \times \mid \; \mid^{\alpha} \chi \rtimes \lambda'\\
A(w'', \lambda) \downarrow w'' && A(w_2, w_1 \lambda) \downarrow w_2 \\
&& \chi_{1_{F^*}} \times \mid \; \mid^{- \alpha} \chi^{-1}(\overset{-}{\;}) \rtimes \lambda'\\
&& A_3(w_1, w_2 w_1 \lambda) \downarrow w_1\\
\mid \; \mid^{-\alpha} \chi^{-1} (\overset{-}{\;}) \rtimes \sigma_{1, \chi_{1_{F^*}}} & \overset{i_2}{\hookrightarrow} &
 \mid \; \mid^{-\alpha} 
\chi^{-1}(\overset{-}{\;}) \times \chi_{1_{F^*}} \rtimes \lambda'
\end{smallmatrix}
$

\bigskip

$ A(w_1, \lambda) $ is not an isomorphism if and only if $ \mid \; \mid^{\alpha} \chi \times \chi_{1_{F^*}} $ is reducible
if and only if $ \alpha = 1 $ and $ \chi = \chi_{1_{F^*}}.$

$ A(w_2, w_1 \lambda) $ is not an isomorphism if and only if  $ \mid \; \mid^{\alpha} \chi \rtimes \lambda'$ reduces
if and only if $ \alpha = 1 $ and $ \chi = 1 $ or $ \alpha = 1/2 $ and $ \chi \in  X_{\omega_{E/F}}. $

$ A(w_1, w_2 w_1 \lambda) $ is no isomorphism if and only if $ \mid \; \mid^{-\alpha} \chi^{-1}(\overset{-}{\;})
 \times \chi_{1_{F^*}} $
reduces if and only if $ \alpha =  1 $ and $ \chi = \chi_{1_{F^*}}. $

\bigskip

In all other cases $ A(w_1, \lambda), A(w_2, w_1 \lambda), A(w_1, w_2 w_1 \lambda) $ and hence $ A(w'', \lambda) $ are
isomorphisms, and
$ \mid \; \mid^{\alpha} \chi \rtimes \sigma_{1, \chi_{1_{F^*}}}  $ and $ \mid \; \mid^{\alpha} \chi \rtimes 
\sigma_{2, \chi_{1_F^*}} $ are irreducible.

If $ A(w_1, \lambda), A(w_2, w_1 \lambda) $ or $ A(w_1, w_2 w_1 \lambda) $ is no isomorphism, then reducibility of
$ \mid \; \mid^{\alpha} \chi \rtimes \sigma_{1, \chi_{1_{F^*}}} $ 
and $ \mid \; \mid^{\alpha} \chi \rtimes \sigma_{2, \chi_{1_{F^*}}} $ is left to be investigated. It is done in \ref{1chi1F*},
\ref{11chi1F*} and in \ref{1/2chiomegachi1F*}.
\end{Proof}

\bigskip

\begin{center}{4.3.3.6 Representations $ {\chi} \rtimes \sigma_{1, \chi_{1_{F^*}}} $ and $ \chi \rtimes \sigma_{2, \chi_{1_{F^*}}}, 
\; \alpha > 0 $}\end{center}

Let $ \chi $ be a unitary character of $ E^*. $ Let $ \chi_{1_{F^*}} \in X_{1_{F^*}}. $ Recall that $ \chi_{1_{F^*}} \rtimes
\lambda' = \sigma_{1, \chi_{1_{F^*}}} \oplus \sigma_{2, \chi_{1_{F^*}}}, $ where $ \sigma_{1, \chi_{1_{F^*}}} $ and
$ \sigma_{1, \chi_{1_{F^*}}} $ are tempered \cite{Ky}.

\begin{Remark}
By Theorem \ref{chi12} the representations $ \chi \rtimes \sigma_{1, \chi_{1_{F^*}}} $ and $ \chi \rtimes \sigma_{2, \chi_{1_{F^*}}} $ are
reducible if and only if $ \chi \in X_{1_{F^*}} $ such that $ \chi \ncong \chi_{1_{F^*}}. $
\end{Remark}

\bigskip

\section{'Special' Reducibility points of representations of $ U(5) $ with cuspidal support in $ M_0 $}

We determine the irreducible subquotients of the representations whose reducibility has not been examined in Chapter 2.

\bigskip

Let $ \chi_{\omega_{E/F}} \in X_{\omega_{E/F}}. $ Let $ \pi_{1, \chi_{\omega_{E/F}}} $ be the unique irreducible
 square-integrable
subquotient and let $ \pi_{2, \chi_{\omega_{E/F}}} $ be the unique irreducible non-tempered subquotient of
$ \mid \; \mid^{1/2} \chi_{\omega_{E/F}} \rtimes \lambda'. $ Let $ \chi_{1_{F^*}} \in X_{1_{F^*}}. $ Recall that
$ \chi_{1_{F^*}} = \sigma_{1, \chi_{1_{F^*}}} + \sigma_{2, \chi_{1_{F^*}}}, $ where $ \sigma_{1, \chi_{1_{F^*}}} $ and
$ \sigma_{2, \chi_{1_{F^*}}} $ are tempered \cite{Ky}.

\medskip

In Theorem \ref{alphachist} the irreducible subquotients of the following representations are left to be examined:

$  \mid \; \mid^{1/2} \St_{\GL_2} \rtimes \lambda', \; \mid \; \mid^{1/2} 1_{\GL_2} \rtimes \lambda', \;
 \mid \; \mid^{3/2} \St_{\GL_2} \rtimes \lambda' , \;  \mid \; \mid^{3/2} 1_{\GL_2} \rtimes \lambda', \;
\mid \; \mid^{1/2} \chi_{\omega_{E/F}} \St_{\GL_2} \rtimes \lambda', \; \mid \; \mid^{1/2} \chi_{\omega_{E/F}} 1_{\GL_2} 
\rtimes \lambda', \; \mid \; \mid^{1/2} \chi_{1_{F^*}} \St_{\GL_2} \rtimes \lambda', \;  \mid \; \mid^{1/2} 
\chi_{1_{F^*}} 1_{\GL_2} \rtimes \lambda', \; \mid \; \mid \chi_{\omega_{E/F}} \St_{\GL_2} \rtimes \lambda', \; 
\mid \; \mid \chi_{\omega_{E/F}} 1_{\GL_2} \rtimes \lambda'. $

\medskip

In Theorem \ref{alphachistu3} the irreducible subquotients of the following representations are left to be examined:

$ \mid \; \mid^2 1 \rtimes \lambda'(\det) \St_{U(3)}, \mid \; \mid^2 1 \rtimes \lambda'(\det) 1_{U(3)}, \;
 \mid \; \mid 1 \rtimes \lambda'(\det) \St_{U(3)}, \mid \; \mid 1 \rtimes \lambda'(\det) 1_{U(3)},
 \mid \; \mid^{1/2} \chi_{\omega_{E/F}} \rtimes \lambda'(\det) \St_{U(3)}, \mid \; \mid^{1/2} \chi_{\omega_{E/F}}
 \rtimes \lambda'(\det) 1_{U(3)}. $

\medskip

Theorem \ref{alphachipi1} leaves the following representations to be examined:

$ \mid \; \mid^{1/2} \chi_{\omega_{E/F}} \rtimes \pi_{1, \chi_{\omega_{E/F}}}, \mid \; \mid^{1/2} \chi_{\omega_{E/F}}
 \rtimes \pi_{2, \chi_{\omega_{E/F}}}, \; \mid \; \mid^{3/2} \chi_{\omega_{E/F}} \rtimes \pi_{1, \chi_{\omega_{E/F}}},
\;  \mid \; \mid^{3/2} \chi_{\omega_{E/F}}
 \rtimes \pi_{2, \chi_{\omega_{E/F}}}, \;  \mid \; \mid 1 \rtimes \pi_{1, \chi_{\omega_{E/F}}}, \; \mid \; \mid 1
 \rtimes \pi_{2, \chi_{\omega_{E/F}}}, \; \mid \; \mid^{1/2} \chi \rtimes \pi_{1, \chi_{\omega_{E/F}}},  \; \mid \; 
\mid^{1/2} \chi \rtimes \pi_{2, \chi_{\omega_{E/F}}}, \; \chi \in X_{\omega_{E/F}}, \chi \ncong \chi_{\omega_{E/F}}. $

\medskip

Theorem \ref{alphachisigma12} leaves the following representations to be examined:

$ \mid \; \mid \chi_{1_{F^*}} \rtimes \sigma_{1, \chi_{1_{F^*}}}, \;
\mid \; \mid \chi_{1_{F^*}} \rtimes \sigma_{2, \chi_{1_{F^*}}}, \;
 \mid \; \mid 1 \rtimes \sigma_{1, \chi_{1_{F^*}}}, \; \mid \; \mid 1 \rtimes \sigma_{2, \chi_{1_{F^*}}}, \;
\mid \; \mid^{1/2} \chi_{\omega_{E/F}} \rtimes \sigma_{1, \chi_{1_{F^*}}} $ and
$ \mid \; \mid \chi_{\omega_{E/F}} \rtimes \sigma_{2, \chi_{1_{F^*}}}. $

\medskip

All representations are treated in this chapter. We determine whether the irreducible subquotients are unitary.

\subsection{$ \mid \; \mid 1 \times 1 \rtimes \lambda' $}

\medskip

In the Grothendieck group of the category of admissible representations of finite length one has
$  \mid \; \mid 1 \times 1 \rtimes \lambda' = \mid \; \mid^{1/2} \St_{\GL_2} \rtimes \lambda' 
+ \mid \; \mid^{1/2} 1_{\GL_2} \rtimes \lambda' = 1 \rtimes \lambda'(\det) \St_{U(3)} + 1 \rtimes \lambda'(\det) 1_{U(3)}. $

\begin{Theorem}

\label{111}

$ \mid \; \mid^{1/2} \St_{\GL_2} \rtimes \lambda' = \Lg(\mid \; \mid^{1/2} \St_{\GL_2} ; \lambda') + 
\tau_3, $

$ \mid \; \mid^{1/2} 1_{\GL_2} \rtimes \lambda' =  \Lg( \mid \; \mid 1 ; 1 \rtimes \lambda') + \tau_4, $

$ 1 \rtimes \lambda'(\det) \St_{U(3)} = \tau_3 + \tau_4, $

$ 1 \rtimes \lambda'(\det) 1_{U(3)} = \Lg( \mid \; \mid 1 ; 1 \rtimes \lambda') + \Lg( \mid \; \mid^{1/2} \St_{\GL_2} ;
\lambda'), $ where

\medskip

$ \tau_3 = \widehat{\Lg( \mid \; \mid 1 ; 1 \rtimes \lambda')}, \; \tau_4 =
\widehat{\Lg(\mid \; \mid^{1/2} \St_{\GL_2} ; \lambda')^{\textasciicircum}}. \; \tau_3 $ and $ \tau_4 $ are tempered.
All irreducible subquotients are unitary.

\end{Theorem}

\begin{Proof} We have seen in Proposition \ref{chistu3} that

$ 1 \rtimes \lambda'(\det) \St_{U(3)} = \tau_3 + \tau_4, $

$ 1 \rtimes \lambda'(\det) 1_{U(3)} = \Lg( \mid \; \mid 1 ; 1 \rtimes \lambda') + \Lg( \mid \; \mid^{1/2} \St_{\GL_2} ;
\lambda'), $ where

\medskip

$ \tau_3 = \widehat{\Lg( \mid \; \mid 1 ; 1 \rtimes \lambda')}, \; \tau_4 =
\widehat{\Lg(\mid \; \mid^{1/2} \St_{\GL_2} ; \lambda')}, $ and $ \tau_3 $ and $ \tau_4 $ are tempered. 

$ \mid \; \mid^{1/2} \St_{\GL_2} \rtimes \lambda' $ is a subrepresentation of $
 \mid \; \mid 1 \times 1 \rtimes \lambda',  \; \mid \; \mid^{1/2} 1_{\GL_2} \rtimes \lambda' $ is a quotient.
Hence $ \Lg( \mid \; \mid 1 ; 1 \rtimes \lambda') $ is a subquotient of $ \mid \; \mid^{1/2} 1_{\GL_2} \rtimes \lambda'. $ 
$ \mid \; \mid^{1/2} \St_{\GL_2} \rtimes \lambda' $ is the Aubert dual of $ \mid \; \mid^{1/2} 1_{\GL_2} \rtimes \lambda', $
hence $ \tau_3 =  \widehat{\Lg( \mid \; \mid 1 ; 1 \rtimes \lambda')} $ is a subquotient of $
\mid \; \mid^{1/2} \St_{\GL_2} \rtimes \lambda'. \; \Lg(\mid \; \mid^{1/2} \St_{\GL_2} ; \lambda') $ is a subquotient of
$ \mid \; \mid^{1/2} \St_{\GL_2} \rtimes \lambda', $ hence $ \tau_4 = 
\widehat{\Lg(\mid \; \mid^{1/2} \St_{\GL_2} ; \lambda')} $ is a subquotient of 
$ \mid \; \mid^{1/2} 1_{\GL_2} \rtimes \lambda'. $

$ 1 \rtimes \lambda'(\det) \St_{U(3)} $ and $ 1 \rtimes \lambda'(\det) 1_{U(3)} $ are unitary, hence all irreducible
subquotients are
unitary.
\end{Proof}

\bigskip

\subsection{$ \mid \; \mid^2 1 \times \mid \; \mid 1 \rtimes \lambda' $}

In the Grothendieck group of admissible representations of finite length one has

$ \mid \; \mid^2 1 \times \mid \; \mid 1 \rtimes \lambda' =  \mid \; \mid^{3/2} \St_{\GL_2} \rtimes \lambda' +
\mid \; \mid^{3/2} 1_{\GL_2} \rtimes \lambda' = \mid \; \mid^2 1 \rtimes \lambda'(\det) \St_{U(3)} +
\mid \; \mid^2 1 \rtimes \lambda'(\det) 1_{U(3)}. $

\begin{Theorem}

\label{2111}

 The representation $ \mid \; \mid^2 1 \times \mid \; \mid 1 \rtimes \lambda' $ is reducible and we have

$ \mid \; \mid^2 1 \times \mid \; \mid 1 \rtimes \lambda' =  \lambda'(\det) \St_{U(5)} + \Lg(\mid \; \mid^{3/2} \St_{\GL_2} ;
 \lambda') +  \lambda'(\det) 1_{U(5)} + \Lg(\mid \; \mid^2 1 ; \lambda'(\det) \St_{U(3)}). $

\smallskip

$ \lambda'(\det) \St_{U(5)} $ and $ \lambda'(\det) 1_{U(5)} $ are unitary, $ \Lg(\mid \; \mid^{3/2} \St_{\GL_2} ;
 \lambda') $ and $ \Lg(\mid \; \mid^2 1 ; \lambda'(\det) \St_{U(3)}) $ are non-unitary.
\end{Theorem}

\begin{Proof}
By \cite{Ca} 
$ \mid \; \mid^2 1 \times \mid \; \mid 1 \rtimes \lambda' $ is a representation of length 4.

$ \lambda'(\det) 1_{U(5)} = \Lg(\mid \; \mid^2 1 ; \mid \; \mid 1 ; \lambda'), \; \Lg(\mid \; \mid^{3/2} \St_{\GL_2} ; \lambda') $ and 
$ \Lg(\mid \; \mid^2
1 ; \lambda'(\det) \St_{U(3)}) $ are all non-tempered Langlands-quotients supported in 
$ \mid \; \mid^2 1 \otimes \mid \; \mid 1 \otimes \lambda'. $ The subrepresentation $ \lambda'(\det) \St_{U(5)} =
 \widehat{\lambda'(\det) 1_{U(5)}} $ is square-integrable.

\medskip

By results of Casselmann (\cite{MR656064}, page 915)  $ \lambda'(\det) 
\St_{U(5)} $ and its Aubert dual  $ \lambda'(\det) 1_{U(5)} =
\Lg(\mid \; \mid^2 1 ; \mid \; \mid 1 ; \lambda') $ are unitary, $ \Lg(\mid \; \mid^{3/2} \St_{\GL_2} ; \lambda') $ and 
$ \Lg(\mid \; \mid^2
1 ; \lambda'(\det) \St_{U(3)}) $ are not unitary.
\end{Proof}

\subsection{$ \mid \; \mid^{1/2} \chi_{\omega_{E/F}} \times \mid \; \mid^{1/2} \chi_{\omega_{E/F}} \rtimes \lambda' $ }

Let $ \chi_{\omega_{E/F}} \in X_{\omega_{E/F}}. $ Let $ \pi_{1,\chi_{\omega_{E/F}}} $ be the unique irreducible
square-integrable subquotient of $ \mid \; \mid^{1/2} \chi_{\omega_{E/F}} \rtimes \lambda', $ let
$ \pi_{2, \chi_{\omega_{E/F}}} $
be the unique irreducible non-tempered subquotient of $ \mid \; \mid^{1/2} \chi_{\omega_{E/F}} \rtimes \lambda' $
\cite{Ky}.

\medskip

In the Grothendieck group of admissible representations of finite length one has

$ \mid \; \mid^{1/2} \chi_{\omega_{E/F}} \times \mid \; \mid^{1/2} \chi_{\omega_{E/F}} \rtimes \lambda' = 
\chi_{\omega_{E/F}} \St_{\GL_2} \rtimes \lambda' + \chi_{\omega_{E/F}} 1_{\GL_2} \rtimes \lambda' =
\mid \; \mid^{1/2} \chi_{\omega_{E/F}} \rtimes \pi_{1, \chi_{\omega_{E/F}}} + \mid \; \mid^{1/2} \chi_{\omega_{E/F}} 
\rtimes
\pi_{2, \chi_{\omega_{E/F}}}. $

\medskip

\begin{Theorem}

\label{1/21/2chiomega}

The representation $ \mid \; \mid^{1/2} \chi_{\omega_{E/F}} \times \mid \; \mid^{1/2} \chi_{\omega_{E/F}}
 \rtimes \lambda' $
is reducible and we have

 $ \mid \; \mid^{1/2} \chi_{\omega_{E/F}} \times \mid \; \mid^{1/2} \chi_{\omega_{E/F}} \rtimes \lambda' = 
\chi_{\omega_{E/F}} \St_{\GL_2} \rtimes \lambda' + \chi_{\omega_{E/F}} 1_{\GL_2} \rtimes \lambda' =
\mid \; \mid^{1/2} \chi_{\omega_{E/F}} \rtimes \pi_{1, \chi_{\omega_{E/F}}} + \mid \; \mid^{1/2} \chi_{\omega_{E/F}} 
\rtimes
\pi_{2, \chi_{\omega_{E/F}}}. $

Moreover we have

\smallskip

$ \chi_{\omega_{E/F}} \St_{\GL_2} \rtimes \lambda' =  \tau_1 + \tau_2, $

$ \chi_{\omega_{E/F}} 1_{\GL_2} \rtimes \lambda' = \Lg(\mid \; \mid^{1/2} \chi_{\omega_{E/F}} ;
\mid \; \mid^{1/2} \chi_{\omega_{E/F}} ; \lambda') + \Lg(\mid \; \mid^{1/2} \chi_{\omega_{E/F}} ; 
\pi_{1, \chi_{\omega_{E/F}}}) $

$ \mid \; \mid^{1/2} \chi_{\omega_{E/F}} \rtimes \pi_{1, \chi_{\omega_{E/F}}} = 
\Lg(\mid \; \mid^{1/2} \chi_{\omega_{E/F}} ; \pi_{1, \chi_{\omega_{E/F}}}) + \tau_1 $

$ \mid \; \mid^{1/2} \chi_{\omega_{E/F}} \rtimes \pi_{2, \chi_{\omega_{E/F}}} = 
\Lg(\mid \; \mid^{1/2} \chi_{\omega_{E/F}} ; \mid \; \mid^{1/2} \chi_{\omega_{E/F}} ; \lambda') + \tau_2, $

\medskip

where $ \tau_1 $ and $ \tau_2 $ are tempered such that $ \tau_1 = \widehat{\Lg(\mid \; \mid^{1/2} \chi_{\omega_{E/F}} ;
\mid \; \mid^{1/2} \chi_{\omega_{E/F}} ; \lambda')} $ and $ \tau_2 = \widehat{\Lg(\mid \; \mid^{1/2} \chi_{\omega_{E/F}} 
; \pi_{1, \chi_{\omega_{E/F}}})}. $
All irreducible subquotients are unitary.

\end{Theorem}

\begin{Proof}
In Proposition \ref{chist} we have seen that

$ \chi_{\omega_{E/F}} \St_{\GL_2} \rtimes \lambda' =  \tau_1 + \tau_2, $

$ \chi_{\omega_{E/F}} 1_{\GL_2} \rtimes \lambda' = \Lg(\mid \; \mid^{1/2} \chi_{\omega_{E/F}} ;
\mid \; \mid^{1/2} \chi_{\omega_{E/F}} ; \lambda') + \Lg(\mid \; \mid^{1/2} \chi_{\omega_{E/F}} ; 
\pi_{1, \chi_{\omega_{E/F}}}), $

where $ \tau_1 $ and $ \tau_2 $ are tempered such that $ \tau_1 = \widehat{\Lg(\mid \; \mid^{1/2} \chi_{\omega_{E/F}} ;
 \mid \; \mid^{1/2}
\chi_{\omega_{E/F}} ; \lambda')} $ and 
$ \tau_2 = \widehat{\Lg( \mid \; \mid^{1/2} \chi_{\omega_{E/F}} ; \pi_1)}. $

$ \Lg(\mid \; \mid^{1/2} \chi_{\omega_{E/F}} ; \pi_{1, \chi_{\omega_{E/F}}}) $ is a subquotient of $  \mid \; \mid^{1/2} \chi_{\omega_{E/F}}
 \rtimes \pi_{1, \chi_{\omega_{E/F}}}. \; \pi_{2, \chi_{\omega_{E/F}}} $ is a quotient of $ \mid \; \mid^{1/2}
\chi_{\omega_{E/F}} \rtimes \lambda'$ [Ke], hence $ \mid \; \mid^{1/2} \chi_{\omega_{E/F}} \rtimes 
\pi_{2, \chi_{\omega_{E/F}}} $ is a
quotient of $ \mid \; \mid^{1/2} \chi_{\omega_{E/F}} \times \mid \; \mid^{1/2} \chi_{\omega_{E/F}}
 \rtimes \lambda'. \; \Lg(\mid \; \mid^{1/2} \chi_{\omega_{E/F}} ; \mid \; \mid^{1/2} \chi_{\omega_{E/F}} ; \lambda') $
is the unique irreducible Langlands-quotient of $ \mid \; \mid^{1/2} \chi_{\omega_{E/F}} \times \mid \; \mid^{1/2} \chi_{\omega_{E/F}}
 \rtimes \lambda', $ hence $  \Lg(\mid \; \mid^{1/2} \chi_{\omega_{E/F}} ; \mid \; \mid^{1/2} \chi_{\omega_{E/F}} ; \lambda') $
is a quotient of  $ \mid \; \mid^{1/2} \chi_{\omega_{E/F}} \rtimes \pi_{2, \chi_{\omega_{E/F}}}. $
Hence $ \tau_1 =  \widehat{\Lg(\mid \; \mid^{1/2} \chi_{\omega_{E/F}} ;
\mid \; \mid^{1/2} \chi_{\omega_{E/F}} ; \lambda')} $ is a subquotient of $  \mid \; \mid^{1/2} 
\chi_{\omega_{E/F}} \rtimes \pi_{1, \chi_{\omega_{E/F}}} $ and $ \tau_2 = \widehat{\Lg(\mid \; \mid^{1/2} \chi_{\omega_{E/F}}
 ; \pi_{1, \chi_{\omega_{E/F}}})} $ is a subquotient of
$ \mid \; \mid^{1/2} \chi_{\omega_{E/F}} \rtimes \pi_{2, \chi_{\omega_{E/F}}}. $
$ \chi_{\omega_{E/F}} \St_{\GL_2} \rtimes \lambda' $ and $ \chi_{\omega_{E/F}} 1_{\GL_2} \rtimes \lambda' $ are unitary,
hence all irreducible subquotients are unitary.
\end{Proof}

\subsection{$ \mid \; \mid^{3/2} \chi_{\omega_{E/F}} \times \mid \; \mid^{1/2} \chi_{\omega_{E/F}} \rtimes
\lambda' $}

\begin{Theorem}

\label{3/21/2chiomega}

Let $ \chi_{\omega_{E/F}} \in X_{\omega_{E/F}}. $  Let $ \pi_{1, \chi_{\omega_{E/F}}} $ be the unique irreducible 
square-integrable subquotient and
let $ \pi_{2, \chi_{\omega_{E/F}}} $ be the unique irreducible non-tempered subquotient of
$ \mid \; \mid^{1/2} \chi_{\omega_{E/F}} \rtimes 
\lambda'. $
The representation $ \mid \; \mid^{3/2} \chi_{\omega_{E/F}} \times \mid \; \mid^{1/2} \chi_{\omega_{E/F}} \rtimes
\lambda' $ is reducible, and we have $ \mid \; \mid^{3/2} \chi_{\omega_{E/F}} \times \mid \; \mid^{1/2} \chi_{\omega_{E/F}}
\rtimes \lambda' = \mid \; \mid \chi_{\omega_{E/F}} \St_{\GL_2} \rtimes \lambda' + \mid \; \mid \chi_{\omega_{E/F}}
 1_{\GL_2} \rtimes
\lambda' = \mid \; \mid^{3/2} \chi_{\omega_{E/F}} \rtimes \pi_{1, \chi_{\omega_{E/F}}} + \mid \; \mid^{3/2}
 \chi_{\omega_{E/F}} \rtimes \pi_{2, \chi_{\omega_{E/F}}}. $
We have
\medskip

$ \mid \; \mid \chi_{\omega_{E/F}} \St_{\GL_2} \rtimes \lambda'= \Lg(\mid \; \mid \chi_{\omega_{E/F}} \St_{\GL_2} ; 
\lambda') 
+ \delta, $

$ \mid \; \mid \chi_{\omega_{E/F}} 1_{\GL_2} \rtimes \lambda' = \Lg(\mid \; \mid^{3/2} \chi_{\omega_{E/F}} ; \mid \; \mid^{1/2} \chi_{\omega_{E/F}} ;
 \lambda') + \Lg( \mid \; \mid^{3/2} \chi_{\omega_{E/F}} ; \pi_{1, \chi_{\omega_{E/F}}}), $

$ \mid \; \mid^{3/2} \chi_{\omega_{E/F}} \rtimes \pi_{1, \chi_{\omega_{E/F}}} = \Lg(\mid \; \mid^{3/2} \chi_{\omega_{E/F}} ; \pi_1) + 
\delta, $

$ \mid \; \mid^{3/2} \chi_{\omega_{E/F}} \rtimes \pi_{2, \chi_{\omega_{E/F}}} = \Lg(\mid \; \mid^{3/2} \chi_{\omega_{E/F}} ; \mid \; \mid^{1/2} \chi_{\omega_{E/F}} ;
 \lambda') + \Lg(\mid \; \mid \chi_{\omega_{E/F}} \St_{\GL_2} ; \lambda'), $

\medskip

where $ \delta = \widehat{\Lg(\mid \; \mid^{3/2} \chi_{\omega_{E/F}} ; \mid \; \mid^{1/2} \chi_{\omega_{E/F}} ;
 \lambda')} $ is square-integrable. 
$ \Lg(\mid \; \mid \chi_{\omega_{E/F}} \St_{\GL_2} ; \lambda') $ and 
$ \Lg(\mid \; \mid^{3/2} \chi_{\omega_{E/F}} ; \pi_{1, \chi_{\omega_{E/F}}}) $ are not unitary. 

\end{Theorem}

\begin{Proof} $ 
\Lg(\mid \; \mid \chi_{\omega_{E/F}} \St_{\GL_2} ; \lambda'),  \Lg( \mid \; \mid^{3/2} \chi_{\omega_{E/F}} ; 
\pi_{1, \chi_{\omega_{E/F}}}) $ and
$ \Lg(\mid \; \mid^{3/2} \chi_{\omega_{E/F}} ; \mid \; 
\mid^{1/2} \chi_{\omega_{E/F}} ; \lambda') $ are the only non-tempered irreducible subquotients of $ 
 \mid \; \mid^{3/2} \chi_{\omega_{E/F}} \times \mid \; \mid^{1/2} \chi_{\omega_{E/F}} \rtimes
\lambda'. $
\medskip

$ \Lg(\mid \; \mid \chi_{\omega_{E/F}} \St_{\GL_2} ; \lambda') $ is a subquotient of
$  \mid \; \mid \chi_{\omega_{E/F}} \St_{\GL_2} \rtimes \lambda', \; \Lg( \mid \; \mid^{3/2} \chi_{\omega_{E/F}} ; 
\pi_{1, \chi_{\omega_{E/F}}}) $
is a subquotient of $ \mid \; \mid^{3/2} \chi_{\omega_{E/F}} \rtimes \pi_{1, \chi_{\omega_{E/F}}}. $

\medskip

We consider Jaquet-restriction to the minimal parabolic
subgroup:

\smallskip

$ s_{\min}(\mid \; \mid \chi_{\omega_{E/F}} \St_{\GL_2} \rtimes \lambda') = \mid \; \mid^{3/2} \chi_{\omega_{E/F}} \otimes
\mid \; \mid^{1/2} \chi_{\omega_{E/F}} \otimes \lambda' + \mid \; \mid^{3/2} \chi_{\omega_{E/F}} \otimes 
\mid \; \mid^{-1/2} \chi_{\omega_{E/F}} \otimes \lambda' + \mid \; \mid^{-1/2} \chi_{\omega_{E/F}} \otimes \mid \; \mid^{-3/2} 
\chi_{\omega_{E/F}}
 \otimes \lambda' +
\mid \; \mid^{-1/2} \chi_{\omega_{E/F}} \otimes \mid \; \mid^{3/2} \chi_{\omega_{E/F}} \otimes \lambda'. $

$ s_{\min}(\mid \; \mid^{3/2} \chi_{\omega_{E/F}} \rtimes \pi_{1, \chi_{\omega_{E/F}}}) = \mid \; \mid^{3/2} \chi_{\omega_{E/F}} \otimes
\mid \; \mid^{1/2} \chi_{\omega_{E/F}} \otimes \lambda' + \mid \; \mid^{-3/2} \chi_{\omega_{E/F}} \otimes \mid \; \mid^{1/2}
\chi_{\omega_{E/F}} \otimes \lambda' + \mid \; \mid^{1/2} \chi_{\omega_{E/F}} \otimes \mid \; \mid^{3/2} \chi_{\omega_{E/F}} 
\otimes \lambda' +
\mid \; \mid^{1/2} \chi_{\omega_{E/F}} \otimes \mid \; \mid^{-3/2} \chi_{\omega_{E/F}} \otimes \lambda'. $

\medskip

$ \mid \; \mid^{3/2} \chi_{\omega_{E/F}} \otimes \mid \; \mid^{1/2} \chi_{\omega_{E/F}} \otimes \lambda' $ is the
only common irreducible subquotient in the restrictions of $ \mid \; \mid \chi_{\omega_{E/F}} \St_{\GL_2} \rtimes \lambda' $ and of
 $ \mid \; \mid^{3/2} \chi_{\omega_{E/F}} \rtimes \pi_{1, \chi_{\omega_{E/F}}}. $ Hence these representations have exactly
one subquotient
in common, denoted by $ \delta. $ By the Casselmann square-integrability criterion $ \delta $ is square-integrable.

\medskip
 We have $ \mid \; \mid \chi_{\omega_{E/F}} \St_{\GL_2} \rtimes \lambda' + \mid \; \mid \chi_{\omega_{E/F}}
 1_{\GL_2} \rtimes
\lambda' = \mid \; \mid^{3/2} \chi_{\omega_{E/F}} \rtimes \pi_{1, \chi_{\omega_{E/F}}} + \mid \; \mid^{3/2} 
\chi_{\omega_{E/F}} \rtimes \pi_{2, \chi_{\omega_{E/F}}}, $
hence $ \Lg(\mid \; \mid \chi_{\omega_{E/F}} \St_{\GL_2} ; \lambda') $ is a subquotient of $ \mid \; \mid^{3/2} 
\chi_{\omega_{E/F}} \rtimes \pi_{2, \chi_{\omega_{E/F}}}, $ and $ \Lg( \mid \; \mid^{3/2} \chi_{\omega_{E/F}} ; 
\pi_{1, \chi_{\omega_{E/F}}}) $ is a subquotient of
$ \mid \; \mid \chi_{\omega_{E/F}} 1_{\GL_2} \rtimes \lambda'. $

\smallskip

$ \mid \; \mid \chi_{\omega_{E/F}} 1_{\GL_2} $ is the Langlands quotient of $ \mid \; \mid^{3/2} \chi_{\omega_{E/F}}
\times \mid \; \mid^{1/2} \chi_{\omega_{E/F}}. \; \mid \; \mid \chi_{\omega_{E/F}} 1_{\GL_2} \otimes \lambda' $ is a quotient
of $ \mid \; \mid^{3/2} \chi_{\omega_{E/F}} \times \mid \; \mid^{1/2} \chi_{\omega_{E/F}} \otimes \lambda'. $ Hence
$ \mid \; \mid \chi_{\omega_{E/F}} 1_{\GL_2} \rtimes \lambda' $ is a quotient of $  \mid \; \mid^{3/2} \chi_{\omega_{E/F}}
\times \mid \; \mid^{1/2} \chi_{\omega_{E/F}} \rtimes \lambda'. $

$ \Lg(\mid \; \mid^{3/2} \chi_{\omega_{E/F}} ; \mid \; \mid^{1/2} \chi_{\omega_{E/F}} ; \lambda') $ is the unique
irreducible quotient of
$ \mid \; \mid^{3/2} \chi_{\omega_{E/F}} \times \mid \; \mid^{1/2} \chi_{\omega_{E/F}} \rtimes \lambda', $ in particular
it is a quotient of $ \mid \; \mid \chi_{\omega_{E/F}} 1_{\GL_2} \rtimes \lambda'. $ In the same manner
$ \Lg(\mid \; \mid^{3/2} \chi_{\omega_{E/F}} ; \mid \; \mid^{1/2} \chi_{\omega_{E/F}} ; \lambda') $ is a quotient
of $ \mid \; \mid^{3/2} \chi_{\omega_{E/F}} \rtimes \pi_{2, \chi_{\omega_{E/F}}}. $

\medskip

So far we have shown:

\medskip

$ \mid \; \mid \chi_{\omega_{E/F}} \St_{\GL_2} \rtimes \lambda'= \Lg(\mid \; \mid \chi_{\omega_{E/F}} \St_{\GL_2} ;
 \lambda') 
+ \delta + A_1 $

$ \mid \; \mid \chi_{\omega_{E/F}} 1_{\GL_2} \rtimes \lambda' = \Lg(\mid \; \mid^{3/2} \chi_{\omega_{E/F}} ; \mid \; \mid^{1/2} \chi_{\omega_{E/F}} ;
 \lambda') + \Lg( \mid \; \mid^{3/2} \chi_{\omega_{E/F}} ; \pi_{1, \chi_{\omega_{E/F}}}) + A_2 $

$ \mid \; \mid^{3/2} \chi_{\omega_{E/F}} \rtimes \pi_{1, \chi_{\omega_{E/F}}} = \Lg(\mid \; \mid^{3/2}
 \chi_{\omega_{E/F}} ; \pi_{1, \chi_{\omega_{E/F}}}) + 
\delta + A_3 $

$ \mid \; \mid^{3/2} \chi_{\omega_{E/F}} \rtimes \pi_{2, \chi_{\omega_{E/F}}} = \Lg(\mid \; \mid^{3/2} \chi_{\omega_{E/F}} ; \mid \; \mid^{1/2} \chi_{\omega_{E/F}} ;
 \lambda') + \Lg(\mid \; \mid \chi_{\omega_{E/F}} \St_{\GL_2} ; \lambda') + A_4, $

\medskip

where $ A_1, A_2, A_3, A_4 $ are sums of tempered representations. We will prove that $ A_1 = A_2 = A_3 = A_4 = 0. $
\bigskip

A tempered representation is the subquotient of a representation induced from a square-integrable representation.
Here, for each proper Levi subgroup $ M_i, \; i = 0,1,2 $ of $ U(5), \; \Ind_{M_0}^{M_i}(\mid \; \mid^{3/2} \chi_{\omega_{E/F}}
\otimes \mid \; \mid^{1/2} \chi_{\omega_{E/F}} \otimes \lambda') $ does not contain any square-integrable subquotient.
Hence all tempered subquotients of $ \mid \; \mid^{3/2} \chi_{\omega_{E/F}} \times \mid \; \mid^{1/2} \chi_{\omega_{E/F}}
\rtimes \lambda' $ are square-integrable.

\medskip

Assume there existed a square-integrable subquotient $ \beta $ of $ \mid \; \mid \chi_{\omega_{E/F}} 1_{\GL_2} 
\rtimes \lambda'. $

\smallskip

We consider Jaquet restriction to the minimal parabolic subgroup.

\smallskip

$ s_{\min}( \mid \; \mid \chi_{\omega_{E/F}} 1_{\GL_2} \rtimes \lambda') =  \mid \; \mid^{1/2} \chi_{\omega_{E/F}} \otimes
\mid \; \mid^{3/2} \chi_{\omega_{E/F}} \otimes \lambda' + \mid \; \mid^{1/2} \chi_{\omega_{E/F}} \otimes 
\mid \; \mid^{-3/2} \otimes \lambda' + \mid \; \mid^{-3/2} \chi_{\omega_{E/F}} \otimes \mid \; \mid^{-1/2} 
\chi_{\omega_{E/F}}
 \otimes \lambda' +
\mid \; \mid^{-3/2} \chi_{\omega_{E/F}} \otimes \mid \; \mid^{1/2} \chi_{\omega_{E/F}} \otimes \lambda', $

hence by Casselmann square-integrability criterion $ s_{\min}(\beta) =
\mid \; \mid^{1/2} \chi_{\omega_{E/F}} \otimes \mid \; \mid^{3/2} \chi_{\omega_{E/F}} \otimes \lambda'. $ 
By \cite{MR1285969}, Théorème 1.7, $ s_{\min}(\beta^{\textasciicircum}) =
\mid \; \mid^{-1/2} \chi_{\omega_{E/F}} \otimes \mid \; \mid^{-3/2} \chi_{\omega_{E/F}} \otimes \lambda'. $

$ \beta^{\textasciicircum} $ is an irreducible subquotient of $ \mid \; \mid \chi_{\omega_{E/F}} \St_{\GL_2} \rtimes
\lambda'. $ As its restriction is negative and as $ \Lg( \mid \; \mid \chi_{\omega_{E/F}} \St_{\GL_2} ; \lambda') $
 is the only non-tempered subquotient of
$ \mid \; \mid \chi_{\omega_{E/F}} \St_{\GL_2} \rtimes \lambda', \; \beta^{\textasciicircum} $ must equal
$ \Lg( \mid \; \mid \chi_{\omega_{E/F}}
\St_{\GL_2} \rtimes \lambda'). $

 \smallskip

We have seen that $ s_{\min}(\mid \; \mid \chi_{\omega_{E/F}} \St_{\GL_2} \rtimes \lambda') = 
\mid \; \mid^{3/2} \chi_{\omega_{E/F}} \otimes
\mid \; \mid^{1/2} \chi_{\omega_{E/F}} \otimes \lambda' + \mid \; \mid^{3/2} \chi_{\omega_{E/F}} \otimes 
\mid \; \mid^{-1/2} \otimes \lambda' + \mid \; \mid^{-1/2} \chi_{\omega_{E/F}} \otimes \mid \; \mid^{-3/2} 
\chi_{\omega_{E/F}}
 \otimes \lambda' +
\mid \; \mid^{-1/2} \chi_{\omega_{E/F}} \otimes \mid \; \mid^{3/2} \chi_{\omega_{E/F}} \otimes \lambda'. $

\smallskip

So $ s_{\min}(\Lg(\mid \; \mid \chi_{\omega_{E/F}}
 \St_{\GL_2} ; \lambda')) $ must contain at least the
two negative irreducible subquotients $ \mid \; \mid^{-1/2} 
\chi_{\omega_{E/F}} \otimes \mid \; \mid^{-3/2} 
\chi_{\omega_{E/F}}
 \otimes \lambda' $ and $
\mid \; \mid^{-1/2} \chi_{\omega_{E/F}} \otimes \mid \; \mid^{3/2} \chi_{\omega_{E/F}} \otimes \lambda'. $ Hence
$ \beta^{\textasciicircum} \neq \Lg( \mid \; \mid \chi_{\omega_{E/F}} \St_{\GL_2} ; \lambda'). $

\medskip

We obtain that $ \Lg(\mid \; \mid \chi_{\omega_{E/F}} \St_{\GL_2} ; \lambda') = \widehat{
\Lg( \mid \; \mid^{3/2} \chi_{\omega_{E/F}} ; \pi_{1, \chi_{\omega_{E/F}}})}, $ and
$ \delta = \widehat{\Lg(\mid \; \mid^{3/2} \chi_{\omega_{E/F}} ; \mid \; \mid^{1/2} \chi_{\omega_{E/F}} ; 
\lambda')}. $

\bigskip

$ \delta = \widehat{\Lg( \mid \; \mid^{3/2} 
\chi_{\omega_{E/F}} ; \mid \; \mid^{1/2} \chi_{\omega_{E/F}}; \lambda')} $ is square-integrable, hence
unitary. $ \Lg(\mid \; \mid^{3/2} \chi_{\omega_{E/F}} ; \mid \; \mid^{1/2} \chi_{\omega_{E/F}} ;
 \lambda') $ is the dual of a square-integrable representation. It
should be unitary, but we have no proof for it. See \cite{MR2460906}, where the proof for the 
unitarisability of the Aubert dual of a strongly positive square-integrable representation is given for orthogonal and 
symplectic groups. Applying Theorem 1.1 and Remark 4.7 of \cite{MR2652536} to the representation
$ \mid \; \mid^{3/2} \chi_{\omega_{E/F}} \rtimes \pi_{1, \chi_{\omega_{E/F}}} $ we see that $ 
\Lg(\mid \; \mid \chi_{\omega_{E/F}} \St_{\GL_2} ; \lambda') $
and $ \Lg( \mid \; \mid^{3/2} \chi_{\omega_{E/F}} ; \pi_{1, \chi_{\omega_{E/F}}}) $ are non-unitary.
\end{Proof}

\subsection{$ \mid \; \mid \chi_{\omega_{E/F}} \times \chi_{\omega_{E/F}} \rtimes \lambda'$}

\label{1chiomega}

\medskip

In the Grothendieck group of admissible representations of finite length one has

$  \mid \; \mid \chi_{\omega_{E/F}} \times \chi_{\omega_{E/F}} \rtimes \lambda' = \mid \; \mid^{1/2} \chi_{\omega_{E/F}}
\St_{\GL_2} \rtimes \lambda' + \mid \; \mid^{1/2} \chi_{\omega_{E/F}} 1_{\GL_2} \rtimes \lambda'. $

\medskip

We have no proof that $ \mid \; \mid^{1/2} \chi_{\omega_{E/F}}
\St_{\GL_2} \rtimes \lambda' $ and $ \mid \; \mid^{1/2} \chi_{\omega_{E/F}} 1_{\GL_2} \rtimes \lambda' $ are irreducible.
See \cite{MR1658535}, Proposition 6.3, where a proof is given for symplectic and special orthogonal groups and when
the representation of the $ \GL_{2p}- $ part, $ p \geq 1, $ of the inducing representation, is essentially
square-integrable.

\medskip

\begin{Remark}
If we assume that $ \mid \; \mid^{1/2} \chi_{\omega_{E/F}}
\St_{\GL_2} \rtimes \lambda' $ and by Aubert duality $ \mid \; \mid^{1/2} \chi_{\omega_{E/F}} 1_{\GL_2} \rtimes \lambda' $ are irreducible,
then $ \mid \; \mid^{1/2} \chi_{\omega_{E/F}} \St_{\GL_2} \rtimes \lambda' = \Lg(\mid \; \mid^{1/2} 
\chi_{\omega_{E/F}} \St_{\GL_2}
 ;\lambda') $
and $ \mid \; \mid^{1/2} \chi_{\omega_{E/F}} 1_{\GL_2}
\rtimes \lambda' = \Lg(\mid \; \mid \chi_{\omega_{E/F}} ; \chi_{\omega_{E/F}} \rtimes \lambda'),$ and both
subquotients are non-unitary.
\smallskip

Further we are able to prove that $ \mid \; \mid^{\alpha} \chi_{\omega_{E/F}} \St_{\GL_2} \rtimes \lambda' =
\Lg(\mid \; \mid^{\alpha} \chi_{\omega_{E/F}} \St_{\GL_2} ;\lambda') $ and 
$ \mid \; \mid^{\alpha} \chi_{\omega_{E/F}} 1_{\GL_2} \rtimes \lambda' = \Lg(\mid \; \mid^{\alpha_1}  \chi_{\omega_{E/F}} 
; \mid \; \mid^{\alpha_2} \chi_{\omega_{E/F}} \rtimes \lambda'),$ are 
non-unitary for $ 0 < \alpha < 1/2, \; 1/2 < \alpha_1 < 1, \; \alpha_2 = 1 - \alpha_1. $ See Remarks \ref{remark1chiomega}
and \ref{remarkalphachiomegast}.
\end{Remark}

\bigskip

\subsection{$ \mid \; \mid \chi_{1_{F^*}} \times \chi_{1_{F^*}} \rtimes 
\lambda' $}

\label{1chi1F*}

We can not give a complete decompositon of  $ \mid \; \mid \chi_{1_{F^*}} \times \chi_{1_{F^*}} \rtimes 
\lambda' $ into irreducible subquotients. We have the following result:

\medskip

Let $ \chi_{1_{F^*}} \in X_{1_{F^*}}. $ By \cite{Ky} $ \chi_{1_{F^*}} \rtimes \lambda' = \sigma_{1, \chi_{1_{F^*}}} +
 \sigma_{2, \chi_{1_{F^*}}}, $ where $ \sigma_{1, \chi_{1_{F^*}}} $ and $ \sigma_{2, \chi_{1_{F^*}}} $ are irreducible
tempered. The representation $ \mid \; \mid \chi_{1_{F^*}} \times \chi_{1_{F^*}} \rtimes 
\lambda' $ is reducible and we have $  \mid \; \mid \chi_{1_{F^*}} \times \chi_{1_{F^*}} \rtimes \lambda' =
\mid \; \mid^{1/2} \chi_{1_{F^*}} \St_{\GL_2} \rtimes \lambda' + \mid \; \mid^{1/2} \chi_{1_{F^*}} 1_{\GL_2} \rtimes \lambda' =
\mid \; \mid \chi_{1_{F^*}} \rtimes \sigma_{1, \chi_{1_{F^*}}} + \mid \; \mid \chi_{1_{F^*}} \rtimes \sigma_{2, \chi_{1_{F^*}}}
. $ Moreover,

\medskip

$ \mid \; \mid^{1/2} \chi_{1_{F^*}} \St_{\GL_2} \rtimes \lambda' = \Lg(\mid \; \mid^{1/2} \chi_{1_{F^*}} \St_{\GL_2} ;
 \lambda') + \delta + A_1, $

$ \mid \; \mid^{1/2} \chi_{1_{F^*}} 1_{\GL_2} \rtimes \lambda' = \Lg(\mid \; \mid \chi_{1_{F^*}} ; \sigma_{1, \chi_{1_{F^*}}}) +
\Lg(\mid \; \mid \chi_{1_{F^*}} ; \sigma_{2, \chi_{1_{F^*}}}) + A_2, $

$ \mid \; \mid \chi_{1_{F^*}} \rtimes \sigma_{1, \chi_{1_{F^*}}} = \Lg(\mid \; \mid \chi_{1_{F^*}} ;
 \sigma_{1, \chi_{1_{F^*}}}) + \delta + A_2, $

$\mid \; \mid \chi_{1_{F^*}} \rtimes \sigma_{2, \chi_{1_{F^*}}} = \Lg(\mid \; \mid \chi_{1_{F^*}} ; \sigma_{2, \chi_{1_{F^*}}}
) + \Lg(\mid \; \mid^{1/2} 
\chi_{1_{F^*}} \St_{\GL_2} ; \lambda') + A_1, $

\medskip

where $ \delta $ is square-integrable. $ \Lg(\mid \; \mid^{1/2} \chi_{1_{F^*}} \St_{\GL_2} ;
 \lambda'), \; \Lg(\mid \; \mid \chi_{1_{F^*}} ; \sigma_{1, \chi_{1_{F^*}}}) $ and $ \Lg(\mid \; \mid \chi_{1_{F^*}} ;
 \sigma_{2, \chi_{1_{F^*}}}) $ are unitary. $ A_1 $ and $ A_2 $ are either both empty, or $ A_1 $ is equal to $ \delta $
 or $ \delta', $ where $ \delta' $ is square-integrable and $ \delta \neq \delta', $ 
and $ A_2 $ is either equal to $ \Lg(\mid \; \mid \chi_{1_{F^*}} ; \sigma_{2, \chi_{1_{F^*}}}) $ or to
$ \Lg(\mid \; \mid^{1/2} \chi_{1_{F^*}} \St_{\GL_2} ;
 \lambda'). $

\bigskip

We prove this assertion:
$ \Lg( \mid \; \mid^{1/2} \chi_{1_{F^*}} \St_{\GL_2} ; \lambda'), \Lg(\mid \; \mid \chi_{1_{F^*}} ; \sigma_{1, \chi_{1_{F^*}}}
) $ and
$ \Lg(\mid \; \mid \chi_{1_{F^*}} ; \sigma_{2, \chi_{1_{F^*}}}) $ are the only non-tempered irreducible subquotients of
$ \mid \; \mid \chi_{1_{F^*}} \times \chi_{1_{F^*}} \rtimes \lambda' . $  
$ \Lg( \mid \; \mid^{1/2} \chi_{1_{F^*}} \St_{\GL_2} ; \lambda') $ is a subquotient of $ \mid \; \mid^{1/2} 
\chi_{1_{F^*}} \St_{\GL_2} \rtimes \lambda'. $ We consider Jaquet restriction to the minimal parabolic subgroup.

\smallskip

$ s_{\min}(\mid \; \mid^{1/2} \chi_{1_{F^*}} \St_{\GL_2} \rtimes \lambda') = \mid \; \mid \chi_{1_{F^*}} \otimes 
\chi_{1_{F^*}} \otimes \lambda' + \mid \; \mid \chi_{1_{F^*}} \otimes \chi_{1_{F^*}} \otimes \lambda' + 
\chi_{1_{F^*}} \otimes \mid \; \mid^{-1} \chi_{1_{F^**}} \otimes \lambda' +
\chi_{1_{F^*}} \otimes \mid \; \mid \chi_{1_{F^*}} \otimes \lambda' =
2   \mid \; \mid \chi_{1_{F^*}} \otimes 
\chi_{1_{F^*}} \otimes \lambda' + \chi_{1_{F^*}} \otimes \mid \; \mid^{-1} \chi_{1_{F^*}} \otimes \lambda' +
\chi_{1_{F^*}} \otimes \mid \; \mid \chi_{1_{F^*}} \otimes \lambda'. $

\smallskip

By the Casselman square-integrability criterion $ \Lg( \mid \; \mid^{1/2} \chi_{1_{F^*}} \St_{\GL_2} ; \lambda') $ is the
only non-tempered irreducible subquotient of $ \mid \; \mid^{1/2} \chi_{1_{F^*}} \St_{\GL_2} \rtimes \lambda' . $ 
A tempered representation is the subquotient of a representation induced from a 
square-integrable representation. $ \mid \; \mid^{1/2} \chi_{1_{F^*}} \St_{\GL_2} \otimes \lambda' $
is not square-integrable, hence any other subquotient of $ \mid \; \mid^{1/2} \chi_{1_{F^*}} \St_{\GL_2} \rtimes \lambda' $
must be square-integrable. Hence $ \Lg(\mid \; \mid \chi_{1_{F^*}} ; \sigma_{1, \chi_{1_{F^*}}}) $ and 
$ \Lg(\mid \; \mid \chi_{1_{F^*}} ; \sigma_{2, \chi_{1_{F^*}}}) $ are 
subquotients of $ \mid \; \mid^{1/2} \chi_{1_{F^*}} 1_{\GL_2} \rtimes \lambda'. $ Let $ \delta $ be a square-integrable
subquotient of $ \mid \; \mid^{1/2} \chi_{1_{F^*}} \St_{\GL_2} \rtimes \lambda'. $

\smallskip

$ \mid \; \mid \chi_{1_{F^*}} \rtimes \sigma_{1, \chi_{1_{F^*}}} $ and 
$ \mid \; \mid \chi_{1_{F^*}} \rtimes \sigma_{2, \chi_{1_{F^*}}} $ have the same
Jaquet-restrictions: 

\smallskip

$ s_{\min}(\mid \; \mid \chi_{1_{F^*}} \rtimes \sigma_{1, \chi_{1_{F^*}}}) = s_{\min}(\mid \; \mid \chi_{1_{F^*}} \rtimes 
\sigma_{2, \chi_{1_{F^*}}}) 
= \mid \; \mid \chi_{1_{F^*}} \otimes 
\chi_{1_{F^*}} \otimes \lambda' + \mid \; \mid^{-1} \chi_{1_{F^*}} \otimes \chi_{1_{F^*}} \otimes \lambda' + \chi_{1_{F^*}} 
\otimes 
\mid \; \mid \chi_{1_{F^**}}
\otimes \lambda' +
\chi_{1_{F^*}} \otimes \mid \; \mid^{-1} \chi_{1_{F^*}} \otimes \lambda'. $

\smallskip

 We chose $ \sigma_{1, \chi_{1_{F^*}}} $ and $ \sigma_{2, \chi_{1_{F^*}}} $ such that $ \delta $ is a subquotient of
$ \mid \; \mid \chi_{1_{F^*}} \rtimes \sigma_{1, \chi_{1_{F^*}}} $ and $ \Lg( \mid \; \mid^{1/2} \chi_{1_{F^*}} \St_{\GL_2} ; \lambda') $
is a subquotient of $ \mid \; \mid \chi_{1_{F^*}} \rtimes \sigma_{2, \chi_{1_{F^*}}}. $

\smallskip

We have no contradiction that $ \mid \; \mid^{1/2} \chi_{1_{F^*}} \St_{\GL_2} \rtimes \lambda' $ contains a second
irreducible square-integrable subquotient $ \delta' $ that would be a subquotient of $
 \mid \; \mid \chi_{1_{F^*}} \rtimes \sigma_{2, \chi_{1_{F^*}}}, $ and that $
 \mid \; \mid^{1/2} \chi_{1_{F^*}} 1_{\GL_2} \rtimes \lambda' $ either contains $ \Lg(\mid \; \mid \chi_{1_{F^*}} ;
 \sigma_{2, \chi_{1_{F^*}}}) $ with multiplicity 2 or $ \Lg(\mid \; \mid^{1/2} \chi_{1_{F^*}} \St_{\GL_2} \rtimes
 \lambda'). \; \Lg(\mid \; \mid \chi_{1_{F^*}} ;
 \sigma_{2, \chi_{1_{F^*}}}) $ and $ \Lg(\mid \; \mid^{1/2} \chi_{1_{F^*}} \St_{\GL_2} \rtimes
 \lambda') $ would then be subquotients of $ \mid \; \mid \chi_{1_{F^*}} \rtimes \sigma_{1, \chi_{1_{F^*}}}. $

\medskip

Let $ 0 < \alpha < 1. $ By Theorem \ref{alphachisigma12} the representations $ \mid \; \mid^{\alpha} \chi_{1_{F^*}} \rtimes
 \sigma_{1,\chi_{1_{F^*}}} $ and 
$ \mid \; \mid^{\alpha} \chi_{1_{F^*}} 
\rtimes \sigma_{2, \chi_{1_{F^*}}} $ are irreducible, they are equal to $ \Lg(\mid \; \mid^{\alpha} \chi_{1_{F^*}}
; \sigma_{1,\chi_{1_{F^*}}}) $ and to $ \Lg(\mid \; \mid^{\alpha} \chi_{1_{F^*}} 
; \sigma_{2, \chi_{1_{F^*}}}), $ respectively. By Theorem \ref{alphachitau} (2) they are unitary. For $ \alpha = 1, $ by 
\cite{MR0324429} the 
irreducible subquotients of 
$ \mid \; \mid \chi_{1_{F^*}} \rtimes
\sigma_{1, \chi_{1_{F^*}}} $ and of $ \mid \; \mid \chi_{1_{F^*}} \rtimes \sigma_{2, \chi_{1_{F^*}}} $ are unitary.

\subsection{$ \mid \; \mid 1 \times \mid \; \mid^{1/2} \chi_{\omega_{E/F}} \rtimes \lambda' $}

Recall that $ \lambda'(\det) \St_{U(3)} $ is the unique
irreducible square-integrable subquotient and that $ \lambda'(\det) 1_{U(3)} $ is the unique irreducible
non-tempered subquotient of $ \mid \; \mid 1 \rtimes \lambda'. $ Let $ \chi_{\omega_{E/F}} \in X_{\omega_{E/F}}. $
Let $ \pi_{1, \chi_{\omega_{E/F}}} $ be the unique irreducible square-integrable subquotient and let
$ \pi_{2, \chi_{\omega_{E/F}}} $ be the unique irreducible non-tempered subquotient of $ \mid \; \mid^{1/2} 
\chi_{\omega_{E/F}}
\rtimes \lambda' $ \cite{Ky}.

\begin{Theorem}
\label{11/21chiomega}

 The representation $ \mid \; \mid 1 \times \mid \; \mid^{1/2} \chi_{\omega_{E/F}} \rtimes \lambda' $ is reducible, and
we have $ \mid \; \mid 1 \times \mid \; \mid^{1/2} \chi_{\omega_{E/F}} \rtimes \lambda' =
\mid \; \mid^{1/2} \chi_{\omega_{E/F}} \rtimes \lambda'(\det) \St_{U(3)} + \mid \; \mid^{1/2} \chi_{\omega_{E/F}} \rtimes
\lambda'(\det) 1_{U(3)}  =  \mid \; \mid 1 \rtimes  \pi_{1, \chi_{\omega_{E/F}}} +  \mid \; \mid 1 \rtimes 
\pi_{2, \chi_{\omega_{E/F}}}. $ Moreover we have

\medskip

$\mid \; \mid^{1/2} \chi_{\omega_{E/F}} \rtimes \lambda'(\det) \St_{U(3)} = \Lg(\mid \; \mid^{1/2} \chi_{\omega_{E/F}} ;
\lambda'(\det) \St_{U(3)}) + \delta, $

$ \mid \; \mid^{1/2} \chi_{\omega_{E/F}} \rtimes \lambda'(\det) 1_{U(3)}  = \Lg( \mid \; \mid 1 ; 
\mid \; \mid^{1/2} \chi_{\omega_{E/F}} ; \lambda') + \Lg(\mid \; \mid 1 ; \pi_{1, \chi_{\omega_{E/F}}}),$

$ \mid \; \mid 1 \rtimes  \pi_{1, \chi_{\omega_{E/F}}} = \Lg(\mid \; \mid 1 ; \pi_{1, \chi_{\omega_{E/F}}}) + \delta, $

$  \mid \; \mid 1 \rtimes \pi_{2, \chi_{\omega_{E/F}}} = \Lg( \mid \; \mid 1 ; \mid \; \mid^{1/2} \chi_{\omega_{E/F}} ; \lambda')
+  \Lg(\mid \; \mid^{1/2} \chi_{\omega_{E/F}} ; \lambda'(\det) \St_{U(3)}), $

\medskip

where $ \delta $ is square-integrable. $ \delta = \widehat{\Lg( \mid \; \mid 1 ; 
\mid \; \mid^{1/2} \chi_{\omega_{E/F}} ; \lambda')}, $ and $ \Lg(\mid \; \mid^{1/2} 
\chi_{\omega_{E/F}} ; \lambda'(\det) \St_{U(3)})  = \widehat{\Lg(\mid \; \mid 1 ; \pi_{1, \chi_{\omega_{E/F}}})}
.$ The representations $ \Lg(\mid \; \mid^{1/2} \chi_{\omega_{E/F}} ; \lambda'(\det) \St_{U(3)}), \;
 \Lg(\mid \; \mid 1 ; \pi_{1, \chi_{\omega_{E/F}}}), \; \Lg( \mid \; \mid 1 ; \mid \; \mid^{1/2} \chi_{\omega_{E/F}} ;
 \lambda') $ and $ \delta $ are all unitary.

\end{Theorem}

\begin{Proof} 
$ \Lg(\mid \; \mid^{1/2} \chi_{\omega_{E/F}} ; \lambda'(\det) \St_{U(3)}), \; \Lg(\mid \; \mid 1 ; 
\pi_{1, \chi_{\omega_{E/F}}}) $ and $ \Lg( \mid \; \mid 1 ; \mid \; \mid^{1/2} \chi_{\omega_{E/F}} ; \lambda') $ are
all the irreducible non-tempered subquotients of $ \mid \; \mid 1 \times \mid \; \mid^{1/2} \chi_{\omega_{E/F}} \rtimes \lambda'. $

$ \Lg(\mid \; \mid^{1/2} \chi_{\omega_{E/F}} ; \lambda'(\det) \St_{U(3)}) $ is a
subquotient of $ \mid \; \mid^{1/2} \chi_{\omega_{E/F}} \rtimes \lambda'(\det) \St_{U(3)}. $

\medskip

We consider Jaquet-restriction to the minimal parabolic subgroup:

$ s_{\min}(\mid \; \mid^{1/2} \chi_{\omega_{E/F}} \rtimes \lambda'(\det) \St_{U(3)}) = \mid \; \mid 1 \otimes
\mid \; \mid^{1/2} \chi_{\omega_{E/F}} \otimes \lambda' +
\mid \; \mid^{1/2} \chi_{\omega_{E/F}} \otimes \mid \; \mid 1 \otimes \lambda' + \mid \; \mid^{-1/2} \chi_{\omega_{E/F}} \otimes \mid \; \mid 1
\otimes \lambda' + \mid \; \mid 1 \otimes \mid \; \mid^{ - 1/2} \chi_{\omega_{E/F}} \otimes \lambda'. $

\medskip

By Casselman square-integrability criterion $ s_{\min}(\Lg(\mid \; \mid^{1/2} \chi_{\omega_{E/F}} ; 
\lambda'(\det) \St_{U(3)}) ) $ must contain the irreducible subquotient
$ \mid \; \mid^{-1/2} \chi_{\omega_{E/F}} \otimes \mid \; \mid 1
\otimes \lambda'. $

\medskip

$ s_{\min}( \mid \; \mid 1 \rtimes  \pi_{1, \chi_{\omega_{E/F}}}) = \mid \; \mid 1 \otimes
\mid \; \mid^{1/2} \chi_{\omega_{E/F}} \otimes \lambda' +
\mid \; \mid^{1/2} \chi_{\omega_{E/F}} \otimes \mid \; \mid 1 \otimes \lambda' + \mid \; \mid^{-1} 1 \otimes \mid \; 
\mid^{1/2}\chi_{\omega_{E/F}} \otimes \lambda' +  \mid \; \mid^{1/2} \chi_{\omega_{E/F}} \otimes \mid \; \mid^{-1} 1 
\otimes \lambda', $ and

$ s_{\min}( \mid \; \mid 1 \rtimes  \pi_{2, \chi_{\omega_{E/F}}}) = \mid \; \mid 1 \otimes
\mid \; \mid^{-1/2} \chi_{\omega_{E/F}} \otimes \lambda' +
\mid \; \mid^{-1/2} \chi_{\omega_{E/F}} \otimes \mid \; \mid 1 \otimes \lambda' + \mid \; \mid^{-1} 1 \otimes \mid \; 
\mid^{-1/2}\chi_{\omega_{E/F}} \otimes \lambda' +  \mid \; \mid^{-1/2} \chi_{\omega_{E/F}} \otimes \mid \; \mid^{-1} 1 
\otimes \lambda'. $ 

\bigskip

The irreducible subquotient $ \mid \; \mid^{-1/2} \chi_{\omega_{E/F}} \otimes \mid \; \mid 1 \otimes \lambda' $ appears
in $ s_{\min}( \mid \; \mid 1 \rtimes  \pi_{2, \chi_{\omega_{E/F}}}), $ not in $ s_{\min}( \mid \; \mid 1 \rtimes 
\pi_{1, \chi_{\omega_{E/F}}}), $ hence
$ \Lg(\mid \; \mid^{1/2} \chi_{\omega_{E/F}} ; \lambda'(\det) \St_{U(3)}) ) $ is a subquotient of $
\mid \; \mid 1 \rtimes  \pi_{2, \chi_{\omega_{E/F}}}. $

\bigskip

$  \Lg(\mid \; \mid 1 ; \pi_{1, \chi_{\omega_{E/F}}}) $ is the unique irreducible quotient of 
$ \mid \; \mid 1 \rtimes  \pi_{1, \chi_{\omega_{E/F}}}. $

\medskip

$ s_{\min}(\mid \; \mid^{1/2} \chi_{\omega_{E/F}} \rtimes \lambda'(\det) 1_{U(3)}) = \mid \; \mid^{1/2} \chi_{\omega_{E/F}}
 \otimes
\mid \; \mid^{-1} 1 \otimes \lambda' + \mid \; \mid^{-1} 1 \otimes \mid \; \mid^{1/2} \chi_{\omega_{E/F}} \otimes  \lambda' +
 \mid \; \mid^{-1/2} \chi_{\omega_{E/F}} \otimes \mid \; 
\mid^{-1} 1 \otimes \lambda' +  \mid \; \mid^{-1} 1 \otimes \mid \; \mid^{-1/2} \chi_{\omega_{E/F}} 
\otimes \lambda'. $ 

\medskip

Looking at Jaquet modules we see that $ \Lg(\mid \; \mid 1 ; \pi_{1, \chi_{\omega_{E/F}}}) $ is a subquotient of $
 \mid \; \mid^{1/2} \chi_{\omega_{E/F}} \rtimes \lambda'(\det) 1_{U(3)}. $

\bigskip

$ \lambda'(\det) 1_{U(3)} $ is a quotient of $ \mid \; \mid 1 \rtimes \lambda' $ \cite{Ky}, hence $
 \mid \; \mid^{1/2} \chi_{\omega_{E/F}} \rtimes \lambda'(\det) 1_{U(3)} $ is a quotient of $ \mid \; \mid 1 \times 
\mid \; \mid^{1/2} \chi_{\omega_{E/F}} \rtimes \lambda'. \; \pi_{2, \chi_{\omega_{E/F}}} $ is a quotient of $ \mid \; \mid^{1/2}
\chi_{\omega_{E/F}} \rtimes \lambda' $ hence $ \mid \; \mid 1 \rtimes \pi_{2, \chi_{\omega_{E/F}}} $ is also a quotient
of $ \mid \; \mid 1 \times 
\mid \; \mid^{1/2} \chi_{\omega_{E/F}} \rtimes \lambda'. \;
\Lg( \mid \; \mid 1 ; \mid \; \mid^{1/2} \chi_{\omega_{E/F}} ; \lambda') $ is the unique irreducible quotient of
$ \mid \; \mid 1 \times \mid \; \mid^{1/2} \chi_{\omega_{E/F}} \rtimes \lambda', $ hence it is a quotient of
$ \mid \; \mid^{1/2} \chi_{\omega_{E/F}} \rtimes \lambda'(\det) 1_{U(3)} $ and of $ \mid \; \mid 1 \rtimes 
\pi_{2, \chi_{\omega_{E/F}}}. $ It is of multiplicity one.

\bigskip

Each irreducible subquotient in $ s_{\min}(\mid \; \mid^{1/2} \chi_{\omega_{E/F}} \rtimes \lambda'(\det) 1_{U(3)}) $ is 
of multiplicity 1. Hence $ \Lg(\mid \; \mid 1 ; 
\pi_{1, \chi_{\omega_{E/F}}}) $ is of multiplicity 1. We have seen that 
$ s_{\min}(\Lg( \mid \; \mid^{1/2} \chi_{\omega_{E/F}} ; \lambda'(\det) \St_{U(3)})) 
$ contains $ \mid \; \mid^{-1/2} \chi_{\omega_{E/F}} \otimes \mid \; \mid 1 \otimes \lambda', $ with multiplicity 1. $ 
\mid \; \mid^{-1/2} \chi_{\omega_{E/F}} \otimes \mid \; \mid 1 \otimes \lambda' $ does not appear in
$ s_{\min}(\mid \; \mid^{1/2} \chi_{\omega_{E/F}} \rtimes \lambda'(\det) 1_{U(3)})), $ hence
$ \Lg(\mid \; \mid^{1/2} \chi_{\omega_{E/F}} ; \lambda'(\det) \St_{U(3)}) $ equally has multiplicity 1.

\medskip

By Casselmann square-integrability criterion any subquotient of $ \mid \; \mid^{1/2} \chi_{\omega_{E/F}} \rtimes \lambda'(\det) 
\St_{U(3)} $ other than $ \Lg(\mid \; \mid^{1/2} \chi_{\omega_{E/F}} ; \lambda'(\det) \St_{U(3)}) $ is square-integrable.

$ \mid \; \mid^{1/2} \chi_{\omega_{E/F}} \rtimes \lambda'(\det) 
1_{U(3)} $ has the two irreducible subquotients $ \Lg( \mid \; \mid 1 ; \mid \; \mid^{1/2} \chi_{\omega_{E/F}} ; \lambda') $
and $ \Lg(\mid \; \mid 1 ; 
\pi_{1, \chi_{\omega_{E/F}}}). $ By Aubert duality $ \mid \; \mid^{1/2} \chi_{\omega_{E/F}} \rtimes \lambda'(\det) 
\St_{U(3)} $ has exactly one square-integrable irreducible subquotient, denoted by $ \delta. $

Looking at Jaquet modules we see that $ \delta $ is a subquotient $ \mid \; \mid 1 \rtimes \pi_{1, \chi_{\omega_{E/F}}}. $

\medskip

$ \delta = \widehat{\Lg( \mid \; \mid 1 ; 
\mid \; \mid^{1/2} \chi_{\omega_{E/F}} ; \lambda')}, $ and $ \Lg(\mid \; \mid^{1/2} 
\chi_{\omega_{E/F}} ; \lambda'(\det) \St_{U(3)})  = \widehat{\Lg(\mid \; \mid 1 ; \pi_{1, \chi_{\omega_{E/F}}})} .$

$ 1 \times \chi_{\omega_{E/F}} \rtimes \lambda' $ is irreducible by Theorem \ref{chi12} and unitary. For $ 0 < \alpha_1 < 1, 
0 < \alpha_2 < 1/2, $ representations $ \mid \; \mid^{\alpha_1} 1 \times \mid \; \mid^{\alpha_2}
\chi_{\omega_{E/F}} \rtimes \lambda' $ are irreducible by Theorem \ref{alpha12chi12} and unitary by Theorem 
\ref{Lg1chiomega} (1). By [Mi],
all irreducible subquotients
of $ \mid \; \mid 1 \times \mid \; \mid^{1/2} \chi_{\omega_{E/F}} \rtimes \lambda' $ are unitary.
\end{Proof}

\subsection{$ \mid \; \mid 1 \times \chi_{1_{F^*}} \rtimes \lambda' $}

Recall that $ \lambda'(\det) \St_{U(3)} $ is the unique irreducible square-integrable subquotient and that $ \lambda'(\det)
1_{U(3)} $ is the unique irreducible non-tempered subquotient of $ \mid \; \mid 1 \rtimes \lambda' $ (\cite{Ky}).  
Let $ \chi_{1_{F^*}} \in X_{1_{F^*}}. $ By \cite{Ky} $ \chi_{1_{F^*}} \rtimes \lambda' =  \sigma_{1,\chi_{1_{F^*}}} + 
\sigma_{2, \chi_{1_{F^*}}}, $ where $ \sigma_{1, \chi_{1_{F^*}}} $ and $ \sigma_{2, \chi_{1_{F^*}}} $ are irreducible 
tempered.

\smallskip

\begin{Theorem}

\label{11chi1F*}

The representation $ \mid \; \mid 1 \times \chi_{1_{F^*}} \rtimes \lambda' $ is reducible and we have
$ \mid \; \mid 1 \times \chi_{1_{F^*}} \rtimes \lambda' =  \chi_{1_{F^*}} \rtimes \lambda'(\det) \St_{U(3)} +
 \chi_{1_{F^*}} \rtimes \lambda'(\det) 1_{U(3)} = \mid \; \mid 1
\rtimes \sigma_{1, \chi_{1_{F^*}}} + \mid \; \mid 1 \rtimes \sigma_{2, \chi_{1_{F^*}}}. $
Furthermore

\medskip

$ \chi_{1_{F^*}} \rtimes \lambda'(\det) \St_{U(3)} = \tau_5 + \tau_6, $

$ \chi_{1_{F^*}} \rtimes \lambda'(\det) 1_{U(3)} = \Lg( \mid \; \mid 1 ; 
\sigma_{1, \chi_{1_{F^*}}})
+ \Lg(\mid \; \mid 1 ; \sigma_{2, \chi_{1_{F^*}}}), $

$ \mid \; \mid 1 \rtimes \sigma_{1, \chi_{1_{F^*}}} = \Lg( \mid \; \mid 1
; \sigma_{1, \chi_{1_{F^*}}}) + \tau_6, $

$ \mid \; \mid 1 \rtimes \sigma_{2, \chi_{1_{F^*}}} = \Lg(\mid \; \mid 1 ;
 \sigma_{2, \chi_{1_{F^*}}}) + \tau_5, $

\smallskip

where $ \tau_5 $ and $ \tau_6 $ are tempered such that $ \tau_5 = \widehat{\Lg(\mid \; \mid 1 ;
 \sigma_{1, \chi_{1_{F^*}}})} $ and $ \tau_6 = \widehat{\Lg(\mid \; \mid 1 ;
 \sigma_{2, \chi_{1_{F^*}}})}. $ All irreducible subquotients are unitary.

\end{Theorem}

\begin{Proof}
$ \Lg( \mid \; \mid 1 ; \sigma_{1, \chi_{1_{F^*}}}) $ and $ \Lg( \mid \; \mid 1
; \sigma_{2, \chi_{1_{F^*}}}) $ are the only non-tempered subquotients of $ \mid \; \mid 1 \times
\chi_{1_{F^*}} \rtimes \lambda'. $

\medskip
$ \Lg( \mid \; \mid 1 ; 
\sigma_{1, \chi_{1_{F^*}}}) $ is the unique irreducible quotient of $ \mid \; \mid 1
\rtimes \sigma_{1, \chi_{1_{F^*}}}. \; \Lg( \mid \; \mid 1 ; \sigma_{2, \chi_{1_{F^*}}}) $ is the
unique irreducible quotient of $ \mid \; \mid 1 \rtimes \sigma_{2, \chi_{1_{F^*}}}. $

\medskip

$ \chi_{1_{F^*}} \rtimes \lambda'(\det) \St_{U(3)} $ is tempered, hence all irreducible subquotients of $
\chi_{1_{F^*}} \rtimes \lambda'(\det) \St_{U(3)} $ are temperd. Hence 
 $ \Lg( \mid \; \mid 1 ; 
\sigma_{1, \chi_{1_{F^*}}}) $ and $ \Lg( \mid \; \mid 1 ; \sigma_{2, \chi_{1_{F^*}}}) $ are subquotients of $ \chi_{1_{F^*}} 
\rtimes \lambda'(\det) 1_{U(3)}. $

\bigskip

Let $ \tau_5 $ and $ \tau_6 $ be two tempered 
subquotients of $ \chi_{1_{F^*}} \rtimes \lambda'(\det) \St_{U(3)}, $ such that $ \tau_5 $ is a subquotient of
$ \mid \; \mid 1 \rtimes \sigma_{2, \chi_{1_{F^*}}} $ and $ \tau_6 $ is a subquotient of
$ \mid \; \mid 1 \rtimes \sigma_{1, \chi_{1_{F^*}}}. $

\bigskip

We now show that no other irreducible subquotients of $ \mid \; \mid 1 \times \chi_{1_{F^*}}
\rtimes \lambda' $ exist. 
Assume there exists a tempered subquotient $ \tau_7 $ of $ \chi_{1_{F^*}} \rtimes \lambda'(\det) \St_{U(3)}. $
Consider Jaquet-restriction to the minimal parabolic subgroup:

\smallskip

$ s_{\min}(\chi_{1_{F^*}} \rtimes \lambda'(\det) \St_{U(3)}) = \chi_{1_{F^*}} \otimes \mid \; \mid 1 \otimes \lambda' + 
\mid \; \mid 1 \otimes \chi_{1_{F^*}} \otimes \lambda' +
\chi_{1_{F^*}} \otimes \mid \; \mid 1 \otimes \lambda' + \mid \; \mid 1 \otimes \chi_{1_{F^*}} \otimes \lambda' =
 2 \; \chi_{1_{F^*}} \otimes \mid \; \mid 1 \otimes \lambda' + 2 \mid \; \mid 1 \otimes \chi_{1_{F^*}} \otimes \lambda'.
$

\medskip

Hence $ \exists \; i \in \{ 5, 6, 7 \} $ such that $ s_{\min}(\tau_i) $ does not contain the irreducible subquotient $
\chi_{1_{F^*}} \otimes \mid \; \mid 1 \otimes \lambda'. $ The Casselman square-integrability criterion implies that $ \tau_i $
is square-integrable. This can not be the case. Hence $ \tau_7 $ does not 
exist, and $ \tau_5 $ and $ \tau_6 $ are of multiplicity 1. By Aubert duality $ \chi_{1_{F^*}} \rtimes \lambda'(\det) 1_{U(3)}
 $ does not have any subquotients other than
$ \Lg( \mid \; \mid 1 ; 
\sigma_{1, \chi_{1_{F^*}}}) $ and $  \Lg(\mid \; \mid 1 ; \sigma_{2, \chi_{1_{F^*}}}), $ both of
multiplicity 1.

\bigskip

We obtain $ \tau_5 =  \widehat{\Lg(\mid \; \mid 1 ; \sigma_{1, \chi_{1_{F^*}}})}, $
and $ \tau_6 =
\widehat{\Lg(\mid \; \mid 1 ; \sigma_{2, \chi_{1_{F^*}}})}. $

\bigskip

$ 1 \rtimes \sigma_{1,\chi_{1_{F^*}}} $ and $ 1 \rtimes \sigma_{2, \chi_{1_{F^*}}} $
are irreducible by Theorem \ref{chi12}. $ 1 $ and by \cite{Ky} $ \sigma_{1, \chi_{1_{F^*}}} $ and $
 \sigma_{2, \chi_{1_{F^*}}} $
are unitary, hence  $ 1 \rtimes \sigma_{1,\chi_{1_{F^*}}} $ and
 $ 1 \rtimes \sigma_{2, \chi_{1_{F^*}}} $ are unitary. For $ 0 < \alpha < 1, $ 
$ \mid \; \mid^{\alpha} 1 \rtimes
\sigma_{1, \chi_{1_{F^*}}} $ and $  \mid \; \mid^{\alpha} 1 \rtimes
\sigma_{2, \chi_{1_{F^*}}} $ are irreducible by Theorem \ref{alphachisigma12} and unitary by Theorem \ref{alphachitau} (2).
By \cite{MR0324429}, all 
irreducible
subquotients of $ \mid \; \mid 1 \times \chi_{1_{F^*}} \rtimes \lambda' $ are unitary.
\end{Proof}

\subsection{$ \mid \; \mid^{1/2} \chi_{\omega_{E/F}} \times \chi_{1_{F^*}} \rtimes \lambda' $}

Let $ \chi_{\omega_{E/F}} \in X_{\omega_{E/F}}. $ Let $ \pi_{1, \chi_{\omega_{E/F}}} $ be the unique square-integrable irreducible
 subquotient and let 
$ \pi_{2, \chi_{\omega_{E/F}}} $ be the unique non-tempered irreducible subquotient of $ \mid \; \mid^{1/2} 
\chi_{\omega_{E/F}} \rtimes \lambda' $ (\cite{Ky}). 
Let $ \chi_{1_{F^*}} \in X_{1_{F^*}}. $ By \cite{Ky} $ \chi_{1_{F^*}} \rtimes \lambda' =  \sigma_{1,\chi_{1_{F^*}}} + 
\sigma_{2, \chi_{1_{F^*}}}, $ where $ \sigma_{1, \chi_{1_{F^*}}} $ and $ \sigma_{2, \chi_{1_{F^*}}} $ are irreducible 
tempered.

\smallskip

\begin{Theorem}

\label{1/2chiomegachi1F*}

The representation $ \mid \; \mid^{1/2} \chi_{\omega_{E/F}} \times \chi_{1_{F^*}} \rtimes \lambda' $ is reducible and we have
$ \mid \; \mid^{1/2} \chi_{\omega_{E/F}} \times \chi_{1_{F^*}} \rtimes \lambda' =  \chi_{1_{F^*}} \rtimes 
\pi_{1, \chi_{\omega_{E/F}}} + \chi_{1_{F^*}} \rtimes \pi_{2, \chi_{\omega_{E/F}}} = \mid \; \mid^{1/2} \chi_{\omega_{E/F}}
\rtimes \sigma_{1, \chi_{1_{F^*}}} + \mid \; \mid^{1/2} \chi_{\omega_{E/F}} \rtimes \sigma_{2, \chi_{1_{F^*}}}. $
Furthermore

\medskip

$ \chi_{1_{F^*}} \rtimes \pi_{1, \chi_{\omega_{E/F}}} = \tau_7 + \tau_8, $

$ \chi_{1_{F^*}} \rtimes \pi_{2, \chi_{\omega_{E/F}}} = \Lg( \mid \; \mid^{1/2} \chi_{\omega_{E/F}} ; 
\sigma_{1, \chi_{1_{F^*}}})
+ \Lg(\mid \; \mid^{1/2} \chi_{\omega_{E/F}} ; \sigma_{2, \chi_{1_{F^*}}}), $

$ \mid \; \mid^{1/2} \chi_{\omega_{E/F}} \rtimes \sigma_{1, \chi_{1_{F^*}}} = \Lg( \mid \; \mid^{1/2} \chi_{\omega_{E/F}}
; \sigma_{1, \chi_{1_{F^*}}}) + \tau_8, $

$ \mid \; \mid^{1/2} \chi_{\omega_{E/F}} \rtimes \sigma_{2, \chi_{1_{F^*}}} = \Lg(\mid \; \mid^{1/2} \chi_{\omega_{E/F}} ;
 \sigma_{2, \chi_{1_{F^*}}}) + \tau_7, $

\smallskip

where $ \tau_7 $ and $ \tau_8 $ are tempered such that $ \tau_7 = \widehat{\Lg(\mid \; \mid^{1/2} \chi_{\omega_{E/F}} ;
 \sigma_{1, \chi_{1_{F^*}}})} $ and $ \tau_8 = \widehat{\Lg(\mid \; \mid^{1/2} \chi_{\omega_{E/F}} ;
 \sigma_{2, \chi_{1_{F^*}}})}. $ All irreducible subquotients are unitary.

\end{Theorem}

\begin{Proof}
$ \Lg( \mid \; \mid^{1/2} \chi_{\omega_{E/F}} ; \sigma_{1, \chi_{1_{F^*}}}) $ and $ \Lg( \mid \; \mid^{1/2} 
\chi_{\omega_{E/F}}
; \sigma_{2, \chi_{1_{F^*}}}) $ are the only non-tempered subquotients of $ \mid \; \mid^{1/2} \chi_{\omega_{E/F}} \times
\chi_{1_{F^*}} \rtimes \lambda'. $

\medskip
$ \Lg( \mid \; \mid^{1/2} \chi_{\omega_{E/F}} ; 
\sigma_{1, \chi_{1_{F^*}}}) $ is the unique irreducible Langlands-quotient of $ \mid \; \mid^{1/2} \chi_{\omega_{E/F}}
\rtimes \sigma_{1, \chi_{1_{F^*}}}. \; \Lg( \mid \; \mid^{1/2} \chi_{\omega_{E/F}} ; \sigma_{2, \chi_{1_{F^*}}}) $ is the
unique irreducible Langlands-quotient of $ \mid \; \mid^{1/2} \chi_{\omega_{E/F}} \rtimes \sigma_{2, \chi_{1_{F^*}}}. $

\medskip

$ \chi_{1_{F^*}} \rtimes \pi_{1, \chi_{\omega_{E/F}}} $ is tempered, hence all irreducible subquotients of $
\chi_{1_{F^*}} \rtimes \pi_{1, \chi_{\omega_{E/F}}} $ are temperd. Hence 
 $ \Lg( \mid \; \mid^{1/2} \chi_{\omega_{E/F}} ; 
\sigma_{1, \chi_{1_{F^*}}}) $ and $ \Lg( \mid \; \mid^{1/2} 
\chi_{\omega_{E/F}} ; \sigma_{2, \chi_{1_{F^*}}}) $ are subquotients of $ \chi_{1_{F^*}} \rtimes
\pi_{2, \chi_{\omega_{E/F}}}. $

\bigskip

Let $ \tau_7 $ and $ \tau_8 $ be two tempered 
subquotients of $ \chi_{1_{F^*}} \rtimes \pi_{1, \chi_{\omega_{E/F}}}, $ such that $ \tau_7 $ is a subquotient of
$ \mid \; \mid^{1/2} \chi_{\omega_{E/F}} \rtimes \sigma_{2, \chi_{1_{F^*}}} $ and $ \tau_8 $ is a subquotient of
$ \mid \; \mid^{1/2} \chi_{\omega_{E/F}} \rtimes \sigma_{1, \chi_{1_{F^*}}}. $

\bigskip

We now show that no other irreducible subquotients of $ \mid \; \mid^{1/2} \chi_{\omega_{E/F}} \times \chi_{1_{F^*}}
\rtimes \lambda' $ exist. 
Assume there exists a tempered subquotient $ \tau_9 $ of $ \chi_{1_{F^*}} \rtimes \pi_{1, \chi_{\omega_{E/F}}}. $
Consider Jaquet-restriction to the minimal parabolic subgroup:

\smallskip

$ s_{\min}(\chi_{1_{F^*}} \rtimes \pi_{1, \chi_{\omega_{E/F}}}) = \chi_{1_{F^*}} \otimes \mid \; \mid^{1/2}
 \chi_{\omega_{E/F}} \otimes \lambda' + \mid \; \mid^{1/2} \chi_{\omega_{E/F}} \otimes \chi_{1_{F^*}} \otimes \lambda' +
\chi_{1_{F^*}} \otimes \mid \; \mid^{1/2} \chi_{\omega_{E/F}} \otimes \lambda' + \mid \; \mid^{1/2} \chi_{\omega_{E/F}}
\otimes \chi_{1_{F^*}} \otimes \lambda' = 2 \; \chi_{1_{F^*}} \otimes \mid \; \mid^{1/2}
 \chi_{\omega_{E/F}} \otimes \lambda' + 2 \mid \; \mid^{1/2} \chi_{\omega_{E/F}} \otimes \chi_{1_{F^*}} \otimes \lambda'.
$

\medskip

Hence $ \exists \; i \in \{ 7, 8, 9 \} $ such that $ s_{\min}(\tau_i) $ does not contain the irreducible subquotient $
\chi_{1_{F^*}} \otimes \mid \; \mid^{1/2}
 \chi_{\omega_{E/F}} \otimes \lambda'. $ The Casselman square-integrability criterion implies that $ \tau_i $
is square-integrable. This can not be the case. Hence $ \tau_9 $ does not 
exist, and $ \tau_7 $ and $ \tau_8 $ are of multiplicity 1. By Aubert duality $ \chi_{1_{F^*}} \rtimes \pi_{2, 
\chi_{\omega_{E/F}}} $ does not have any subquotients other than
$ \Lg( \mid \; \mid^{1/2} \chi_{\omega_{E/F}} ; 
\sigma_{1, \chi_{1_{F^*}}}) $ and $  \Lg(\mid \; \mid^{1/2} \chi_{\omega_{E/F}} ; \sigma_{2, \chi_{1_{F^*}}}), $ both of
multiplicity 1.

\bigskip

We obtain $ \tau_7 =  \widehat{\Lg(\mid \; \mid^{1/2} \chi_{\omega_{E/F}} ; \sigma_{1, \chi_{1_{F^*}}})}, $
and $ \tau_8 = \widehat{\Lg(\mid \; \mid^{1/2} \chi_{\omega_{E/F}} ; \sigma_{2, \chi_{1_{F^*}}})}. $

\bigskip

$ \chi_{\omega_{E/F}} \rtimes \sigma_{1,\chi_{1_{F^*}}} $ and $ \chi_{\omega_{E/F}} \rtimes \sigma_{2, \chi_{1_{F^*}}} $
are irreducible by Theorem \ref{chi12}. $ \chi_{\omega_{E/F}} $ and by \cite{Ky} $ \sigma_{1, \chi_{1_{F^*}}} $ and $
 \sigma_{2, \chi_{1_{F^*}}} $
are unitary, hence  $ \chi_{\omega_{E/F}} \rtimes \sigma_{1,\chi_{1_{F^*}}} $ and
 $ \chi_{\omega_{E/F}} \rtimes \sigma_{2, \chi_{1_{F^*}}} $ are unitary. For $ 0 < \alpha < 1/2, $ 
$ \mid \; \mid^{\alpha} \chi_{\omega_{E/F}} \rtimes
\sigma_{1, \chi_{1_{F^*}}} $ and $  \mid \; \mid^{\alpha} \chi_{\omega_{E/F}} \rtimes
\sigma_{1, \chi_{1_{F^*}}} $ are irreducible by Theorem \ref{alphachisigma12} and unitary by Theorem \ref{alphachitau} (3).
By \cite{MR0324429}, all 
irreducible
subquotients of $ \mid \; \mid^{1/2} \chi_{\omega_{E/F}} \times \chi_{1_{F^*}} \rtimes \lambda' $ are unitary.
\end{Proof}

\subsection{$ \mid \; \mid^{1/2} \chi_{\omega_{E/F}} \times \mid \; \mid^{1/2} \chi_{\omega_{E/F}}' \rtimes
 \lambda' $}

Let $ \chi_{\omega_{E/F}}, \chi_{\omega_{E/F}}' \in X_{\omega_{E/F}}, $ such that $ \chi_{\omega_{E/F}} \neq
 \chi_{\omega_{E/F}}'. $  Let $ \pi_{1, \chi_{\omega_{E/F}}} $ be the unique square-integrable subquotient
 and let 
$ \pi_{2, \chi_{\omega_{E/F}}} $ be the unique non-tempered irreducible subquotient of $ \mid \; \mid^{1/2} 
\chi_{\omega_{E/F}} \rtimes \lambda'. $
 Let $ \pi_{1, \chi_{\omega_{E/F}'}} $ be the unique square-integrable irreducible
 subquotient and let 
$ \pi_{2, \chi_{\omega_{E/F}'}} $ be the unique non-tempered irreducible subquotient of $ \mid \; \mid^{1/2} 
\chi_{\omega_{E/F}}' \rtimes \lambda' $ \cite{Ky}.

\smallskip

\begin{Theorem}

\label{1/21/2chiomega12}

The representation $ \mid \; \mid^{1/2} \chi_{\omega_{E/F}} \times \mid \; \mid^{1/2} \chi_{\omega_{E/F}}' \rtimes \lambda' $
is reducible. We have
$ \mid \; \mid^{1/2} \chi_{\omega_{E/F}} \times \mid \; \mid^{1/2} \chi_{\omega_{E/F}}' \rtimes \lambda' =
\mid \; \mid^{1/2} \chi_{\omega_{E/F}} \rtimes \pi_{1, \chi_{\omega_{E/F}}'} + \mid \; \mid^{1/2} \chi_{\omega_{E/F}} \rtimes
\pi_{2, \chi_{\omega_{E/F}}'} = \mid \; \mid^{1/2} \chi_{\omega_{E/F}}' \rtimes \pi_{1, \chi_{\omega_{E/F}}} + 
\mid \; \mid^{1/2} \chi_{\omega_{E/F}}' \rtimes
\pi_{2, \chi_{\omega_{E/F}}}. $ Furthermore

\smallskip

$ \mid \; \mid^{1/2} \chi_{\omega_{E/F}} \rtimes \pi_{1, \chi_{\omega_{E/F}}'} = \Lg( \mid \; \mid^{1/2} \chi_{\omega_{E/F}}
; \pi_{1, \chi_{\omega_{E/F}}'}) + \delta, $

$ \mid \; \mid^{1/2} \chi_{\omega_{E/F}} \rtimes \pi_{2, \chi_{\omega_{E/F}}'} = \Lg(\mid \; \mid^{1/2} \chi_{\omega_{E/F}}
; \mid \; \mid^{1/2} \chi_{\omega_{E/F}}' ; \lambda') + \Lg( \mid \; \mid^{1/2} \chi_{\omega_{E/F}}' ; \pi_{1,
 \chi_{\omega_{E/F}}}), $

$\mid \; \mid^{1/2} \chi_{\omega_{E/F}}' \rtimes \pi_{1, \chi_{\omega_{E/F}}} = \Lg(\mid \; \mid^{1/2} \chi_{\omega_{E/F}}' ;
\pi_{1, \chi_{\omega_{E/F}}}) + \delta, $

$ \mid \; \mid^{1/2} \chi_{\omega_{E/F}}' \rtimes \pi_{2, \chi_{\omega_{E/F}}} = \Lg(\mid \; \mid^{1/2} \chi_{\omega_{E/F}}
 ; \mid \; \mid^{1/2} \chi_{\omega_{E/F}} ; \lambda') + \Lg( \mid \; \mid^{1/2} \chi_{\omega_{E/F}} ;
 \pi_{1, \chi_{\omega_{E/F}}}'), $

\smallskip

where $ \delta =  \widehat{\Lg(\mid \; \mid^{1/2} \chi_{\omega_{E/F}} ; \mid \; \mid^{1/2} \chi_{\omega_{E/F}}' ; 
\lambda')} $ is square-integrable. $ \Lg(\mid \; \mid^{1/2} \chi_{\omega_{E/F}}
; \mid \; \mid^{1/2} \chi_{\omega_{E/F}}' ; \lambda'), \Lg( \mid \; \mid^{1/2} \chi_{\omega_{E/F}} ; 
\pi_{1, \chi_{\omega_{E/F}}'}) $ and $ \Lg(\mid \; \mid^{1/2} \chi_{\omega_{E/F}}' ; \pi_{1, \chi_{\omega_{E/F}}}) $
are unitary.

\end{Theorem}

\begin{Proof}
$ \Lg(\mid \; \mid^{1/2} \chi_{\omega_{E/F}}
; \mid \; \mid^{1/2} \chi_{\omega_{E/F}}' ; \lambda'), \Lg( \mid \; \mid^{1/2} \chi_{\omega_{E/F}} ; 
\pi_{1, \chi_{\omega_{E/F}}}') $ and $ \Lg(\mid \; \mid^{1/2} \chi_{\omega_{E/F}}' ; \pi_{1, \chi_{\omega_{E/F}}}) $ are all
the non-tempered irreducible subquotients of $ \mid \; \mid^{1/2} \chi_{\omega_{E/F}} \times \mid \; \mid^{1/2}
\chi_{\omega_{E/F}}' \rtimes \lambda'. $

$ \Lg( \mid \; \mid^{1/2} \chi_{\omega_{E/F}} ; \pi_{1, \chi_{\omega_{E/F}}'}) $ is a subquotient of 
$ \mid \; \mid^{1/2} \chi_{\omega_{E/F}} \rtimes \pi_{1, \chi_{\omega_{E/F}}'}. $

Consider Jaquet-restriction to the minimal parabolic subgroup: 

\medskip

$ s_{\min}( \mid \; \mid^{1/2} \chi_{\omega_{E/F}} \rtimes
\pi_{1, \chi_{\omega_{E/F}}'}) = \mid \; \mid^{1/2} \chi_{\omega_{E/F}} \otimes \mid \; \mid^{1/2} \chi_{\omega_{E/F}}'
\otimes \lambda' + \mid \; \mid^{1/2} \chi_{\omega_{E/F}}' \otimes \mid \; \mid^{1/2} \chi_{\omega_{E/F}}
\otimes \lambda' + \mid \; \mid^{-1/2} \chi_{\omega_{E/F}} \otimes \mid \; \mid^{1/2} \chi_{\omega_{E/F}}'
\otimes \lambda' + \mid \; \mid^{1/2} \chi_{\omega_{E/F}}' \otimes \mid \; \mid^{-1/2} \chi_{\omega_{E/F}}
\otimes \lambda', $

\medskip

$ \Lg( \mid \;
\mid^{1/2} \chi_{\omega_{E/F}} ; \pi_{1, \chi_{\omega_{E/F}}'}) $ is non-tempered, hence $ s_{\min}(\Lg( \mid \;
\mid^{1/2} \chi_{\omega_{E/F}} ; \pi_{1, \chi_{\omega_{E/F}}'})) $ must contain the irreducible subquotient
$ \mid \; \mid^{-1/2} \chi_{\omega_{E/F}} \otimes \mid \; \mid^{1/2} \chi_{\omega_{E/F}}'
\otimes \lambda'. $

\medskip

$ s_{\min}( \mid \; \mid^{1/2} \chi_{\omega_{E/F}}' \rtimes
\pi_{1, \chi_{\omega_{E/F}}}) = \mid \; \mid^{1/2} \chi_{\omega_{E/F}}' \otimes \mid \; \mid^{1/2} \chi_{\omega_{E/F}}
\otimes \lambda' + \mid \; \mid^{1/2} \chi_{\omega_{E/F}} \otimes \mid \; \mid^{1/2} \chi_{\omega_{E/F}}'
\otimes \lambda' + \mid \; \mid^{-1/2} \chi_{\omega_{E/F}}' \otimes \mid \; \mid^{1/2} \chi_{\omega_{E/F}}
\otimes \lambda' + \mid \; \mid^{1/2} \chi_{\omega_{E/F}} \otimes \mid \; \mid^{-1/2} \chi_{\omega_{E/F}}'
\otimes \lambda', $

\medskip

$ s_{\min}( \mid \; \mid^{1/2} \chi_{\omega_{E/F}}' \rtimes
\pi_{2, \chi_{\omega_{E/F}}}) = \mid \; \mid^{1/2} \chi_{\omega_{E/F}}' \otimes \mid \; \mid^{-1/2} \chi_{\omega_{E/F}}
\otimes \lambda' + \mid \; \mid^{-1/2} \chi_{\omega_{E/F}} \otimes \mid \; \mid^{1/2} \chi_{\omega_{E/F}}'
\otimes \lambda' + \mid \; \mid^{-1/2} \chi_{\omega_{E/F}}' \otimes \mid \; \mid^{-1/2} \chi_{\omega_{E/F}}
\otimes \lambda' + \mid \; \mid^{-1/2} \chi_{\omega_{E/F}} \otimes \mid \; \mid^{-1/2} \chi_{\omega_{E/F}}'
\otimes \lambda'. $

\bigskip

The irreducible suquotient $ \mid \; \mid^{-1/2} \chi_{\omega_{E/F}} \otimes \mid \; \mid^{1/2} \chi_{\omega_{E/F}}'
\otimes \lambda' $ appears in $ s_{\min}( \mid \; \mid^{1/2} \chi_{\omega_{E/F}}' \rtimes
\pi_{2, \chi_{\omega_{E/F}}}), $ not in $ s_{\min}( \mid \; \mid^{1/2} \chi_{\omega_{E/F}}' \rtimes
\pi_{1, \chi_{\omega_{E/F}}}). $ Hence $ \Lg( \mid \;
\mid^{1/2} \chi_{\omega_{E/F}} ; \pi_{1, \chi_{\omega_{E/F}}'}) $ is also a subquotient of 
$ \mid \; \mid^{1/2} \chi_{\omega_{E/F}}' \rtimes
\pi_{2, \chi_{\omega_{E/F}}}. $

\bigskip

$ \Lg( \mid \; \mid^{1/2} \chi_{\omega_{E/F}}' ; \pi_{1, \chi_{\omega_{E/F}}}) $ is a subquotient of 
$ \mid \; \mid^{1/2} \chi_{\omega_{E/F}}' \rtimes \pi_{1, \chi_{\omega_{E/F}}}. $ In the same manner as above we find
that $ \Lg( \mid \; \mid^{1/2} \chi_{\omega_{E/F}}' ; \pi_{1, \chi_{\omega_{E/F}}}) $ is also a subquotient of 
$ \mid \; \mid^{1/2} \chi_{\omega_{E/F}} \rtimes \pi_{2, \chi_{\omega_{E/F}}'}. $

\medskip

$ \pi_{2,\chi_{\omega_{E/F}}} $ is a quotient of $ \mid \; \mid^{1/2} \chi_{\omega_{E/F}} \rtimes \lambda', $ hence
$ \mid \; \mid^{1/2} \chi_{\omega_{E/F}}' \rtimes \pi_{2,\chi_{\omega_{E/F}}} $ is a quotient of
$ \mid \; \mid^{1/2} \chi_{\omega_{E/F}} \times \mid \; \mid^{1/2} \chi_{\omega_{E/F}}' \rtimes \lambda'.
\; \pi_{2, \chi_{\omega_{E/F}}'} $ is a quotient of $ \mid \; \mid^{1/2} \chi_{\omega_{E/F}}' \rtimes \lambda', $
hence $ \mid \; \mid^{1/2} \chi_{\omega_{E/F}} \rtimes \pi_{2, \chi_{\omega_{E/F}}'} $ is a quotient of $ 
\mid \; \mid^{1/2} \chi_{\omega_{E/F}} \times \mid \; \mid^{1/2} \chi_{\omega_{E/F}}' \rtimes \lambda'. \;
\Lg(\mid \; \mid^{1/2} \chi_{\omega_{E/F}}
; \mid \; \mid^{1/2} \chi_{\omega_{E/F}}' ; \lambda') $ is the unique irreducible quotient of
$ \mid \; \mid^{1/2} \chi_{\omega_{E/F}} \times \mid \; \mid^{1/2} \chi_{\omega_{E/F}}' \rtimes \lambda', $
hence $ \Lg(\mid \; \mid^{1/2} \chi_{\omega_{E/F}}
; \mid \; \mid^{1/2} \chi_{\omega_{E/F}}' ; \lambda') $ is a quotient of $ \mid \; \mid^{1/2}
 \chi_{\omega_{E/F}}' \rtimes \pi_{2, \chi_{\omega_{E/F}}} $ and of $ \mid \; \mid^{1/2} \chi_{\omega_{E/F}} \rtimes 
\pi_{2, \chi_{\omega_{E/F}}'}. $

\medskip

A tempered representation is the subquotient of a representation induced from a square-integrable representation of a
parabolic subgroup.
Here, for $ i = 0,1,2, \; \Ind_{M_0}^{M_i}(\mid \; \mid^{1/2} \chi_{\omega_{E/F}} \otimes \mid \; \mid^{1/2} \chi_{\omega_{E/F}}'
\otimes \lambda') $ does not contain any square-integrable subquotient. Hence any irreducible subquotient of
$ \mid \; \mid^{1/2} \chi_{\omega_{E/F}} \times \mid \; \mid^{1/2}
\chi_{\omega_{E/F}}' \rtimes \lambda' $ other than
$ \Lg(\mid \; \mid^{1/2} \chi_{\omega_{E/F}}
; \mid \; \mid^{1/2} \chi_{\omega_{E/F}}' ; \lambda'), \Lg( \mid \; \mid^{1/2} \chi_{\omega_{E/F}} ; 
\pi_{1, \chi_{\omega_{E/F}}}') $ and $ \Lg(\mid \; \mid^{1/2} \chi_{\omega_{E/F}}' ; \pi_{1, \chi_{\omega_{E/F}}}) $ must be
square-integrable.

$ s_{\min}(\mid \; \mid^{1/2} \chi_{\omega_{E/F}} \rtimes 
\pi_{1, \chi_{\omega_{E/F}}'}) $ contains only one negative subquotient, $ \mid \; \mid^{-1/2} \chi_{\omega_{E/F}} \otimes 
\mid \; \mid^{1/2} \chi_{\omega_{E/F}}' \otimes \lambda', $ hence
$ \Lg( \mid \; \mid^{1/2} \chi_{\omega_{E/F}} ; 
\pi_{1, \chi_{\omega_{E/F}}'}) $ is the only non-tempered irreducible subquotient of $ 
 \mid \; \mid^{1/2} \chi_{\omega_{E/F}} \rtimes \pi_{1, \chi_{\omega_{E/F}}'}. $ 
Let $ \delta $ denote a square-integrable irreducible subquotient of
$ \mid \; \mid^{1/2} \chi_{\omega_{E/F}} \rtimes 
\pi_{1, \chi_{\omega_{E/F}}'}. $
Looking at Jaquet modules we find that $ \delta $ is also a subquotient of $ \mid \; \mid^{1/2} \chi_{\omega_{E/F}}' 
\rtimes \pi_{1, \chi_{\omega_{E/F}}}. $

\medskip

So far we have seen:

\smallskip

$ \mid \; \mid^{1/2} \chi_{\omega_{E/F}} \rtimes \pi_{1, \chi_{\omega_{E/F}}'} = \Lg( \mid \; \mid^{1/2} \chi_{\omega_{E/F}}
; \pi_{1, \chi_{\omega_{E/F}}'}) + \delta + A_1, $

$ \mid \; \mid^{1/2} \chi_{\omega_{E/F}} \rtimes \pi_{2, \chi_{\omega_{E/F}}'} = \Lg(\mid \; \mid^{1/2} 
\chi_{\omega_{E/F}}
; \mid \; \mid^{1/2} \chi_{\omega_{E/F}}' ; \lambda') + \Lg( \mid \; \mid^{1/2} \chi_{\omega_{E/F}}' ; \pi_{1,
 \chi_{\omega_{E/F}}}) + A_2, $

$\mid \; \mid^{1/2} \chi_{\omega_{E/F}}' \rtimes \pi_{1, \chi_{\omega_{E/F}}} = \Lg(\mid \; \mid^{1/2} \chi_{\omega_{E/F}}' ;
\pi_{1, \chi_{\omega_{E/F}}}) + \delta + A_3, $

$ \mid \; \mid^{1/2} \chi_{\omega_{E/F}}' \rtimes \pi_{2, \chi_{\omega_{E/F}}} = \Lg(\mid \; \mid^{1/2} 
\chi_{\omega_{E/F}}
 ; \mid \; \mid^{1/2} \chi_{\omega_{E/F}}' ; 
\lambda') + \Lg( \mid \; \mid^{1/2} \chi_{\omega_{E/F}}
; \pi_{1, \chi_{\omega_{E/F}}'}) + A_4, $

\smallskip

where $ A_1, A_2, A_3 $ and $ A_4 $ are sums of tempered representations.
We will now show that $ A_1, A_2, A_3 $ and $ A_4 $ are empty.

\medskip

$ s_{\min}(\mid \; \mid^{1/2} \chi_{\omega_{E/F}} \rtimes \pi_{2, \chi_{\omega_{E/F}}}') $ does not contain any 
non-negative subquotients. Hence by the Casselman square-integrability criterion all irreducible subquotients of
$ \mid \; \mid^{1/2} \chi_{\omega_{E/F}} \rtimes \pi_{2, \chi_{\omega_{E/F}}}' $ are non-tempered.
Each subquotient in $ s_{\min}(\mid \; \mid^{1/2} \chi_{\omega_{E/F}} \rtimes \pi_{2, \chi_{\omega_{E/F}}'}) $ is of
multiplicity 1. Hence $ \Lg( \mid \; \mid^{1/2} \chi_{\omega_{E/F}} ; \mid \; \mid^{1/2} \chi_{\omega_{E/F}}' ;
\lambda') $ and $ \Lg( \mid \; \mid^{1/2} \chi_{\omega_{E/F}}' ; \pi_{1, \chi_{\omega_{E/F}}}) $ are of
multiplicity 1 in $  \mid \; \mid^{1/2} \chi_{\omega_{E/F}} \rtimes \pi_{2, \chi_{\omega_{E/F}}'}. $
The irreducible subquotient $ \mid \; \mid^{-1/2} \chi_{\omega_{E/F}}
\otimes \mid \; \mid^{1/2} \chi_{\omega_{E/F}}' \otimes \lambda' $ in $ s_{\min}(\Lg( 
\mid \; \mid^{1/2} \chi_{\omega_{E/F}} ; 
\pi_{1, \chi_{\omega_{E/F}}'})) $ does not appear in $ s_{\min}(\mid \; \mid^{1/2} \chi_{\omega_{E/F}} \rtimes 
\pi_{2, \chi_{\omega_{E/F}}'}). $ Hence $ \Lg( 
\mid \; \mid^{1/2} \chi_{\omega_{E/F}} ; 
\pi_{1, \chi_{\omega_{E/F}}'}) $ is no subquotient of $  \mid \; \mid^{1/2} \chi_{\omega_{E/F}} \rtimes 
\pi_{2, \chi_{\omega_{E/F}}}'. $

Equivalently we obtain
that all irreducible subquotients of $ \mid \; \mid^{1/2} \chi_{\omega_{E/F}}' \rtimes \pi_{2, \chi_{\omega_{E/F}}} $
are non-tempered. $ \Lg( \mid \; \mid^{1/2} \chi_{\omega_{E/F}} ; \mid \; \mid^{1/2} \chi_{\omega_{E/F}}' ;
\lambda') $ and $ \Lg( \mid \; \mid^{1/2} \chi_{\omega_{E/F}} ; \pi_{1, \chi_{\omega_{E/F}}}') $ are of
multiplicity 1 and $ \Lg( \mid \; \mid^{1/2} \chi_{\omega_{E/F}} ; \pi_{1,\chi_{\omega_{E/F}}'}) $ is no subquotient of
$ \mid \; \mid^{1/2} \chi_{\omega_{E/F}}' \rtimes \pi_{2, \chi_{\omega_{E/F}}}. $
By Aubert duality $ \mid \; \mid^{1/2} \chi_{\omega_{E/F}} \rtimes \pi_{1, \chi_{\omega_{E/F}}'} $ and 
$ \mid \; \mid^{1/2} \chi_{\omega_{E/F}}' \rtimes \pi_{1, \chi_{\omega_{E/F}}} $ do not have any other subquotients.

\medskip

We obtain that $ \delta = \widehat{\Lg( \mid \; \mid^{1/2} \chi_{\omega_{E/F}} ; \mid \; \mid^{1/2} \chi_{\omega_{E/F}}' ;
\lambda')} $ and $ \Lg( \mid \; \mid^{1/2} \chi_{\omega_{E/F}} ; \pi_{1, \chi_{\omega_{E/F}}}') 
= \widehat{\Lg( \mid \; \mid^{1/2} \chi_{\omega_{E/F}}' ; \pi_{1,\chi_{\omega_{E/F}}})}. $

\medskip

$ \chi_{\omega_{E/F}} \times \chi_{\omega_{E/F}}' \rtimes \lambda' $ is irreducible by Theorem \ref{chi12} and unitary.
For $ 0 < \alpha_1, \alpha_2 < 1/2, \;  \mid \; \mid^{\alpha_1} \chi_{\omega_{E/F}} \times 
\mid \; \mid^{\alpha_2 } \chi_{\omega_{E/F}}' \rtimes \lambda' $ is irreducible by Theorem \ref{alpha12chi12} and unitary by
Theorem \ref{Lgchiomega12} (1). By \cite{MR0324429}, all irreducible subquotients of $  \mid \; \mid^{1/2} \chi_{\omega_{E/F}} \times
 \mid \; \mid^{1/2} \chi_{\omega_{E/F}}' \rtimes
\lambda' $ are unitary.
\end{Proof}

\bibliographystyle{amsalpha}
\providecommand{\bysame}{\leavevmode\hbox to3em{\hrulefill}\thinspace}
\providecommand{\MR}{\relax\ifhmode\unskip\space\fi MR }
\providecommand{\MRhref}[2]{%
  \href{http://www.ams.org/mathscinet-getitem?mr=#1}{#2}
}
\providecommand{\href}[2]{#2}

\bibliography{biblio}
\nocite{*}
\end{document}